\documentclass[11pt,reqno]{amsart}

\usepackage{amsmath, amsfonts, amssymb, amscd,graphics, graphicx, curves,
graphpap, tikz, color, xcolor, enumitem, hyperref, 
mathabx, mathrsfs}

\usepackage[mathscr]{euscript}
  \usepackage[normalem]{ulem}

\makeatletter
\newcommand{\doublewidetilde}[1]{{%
  \mathpalette\double@widetilde{#1}%
}}

\textwidth= 
148mm
\numberwithin{equation}{section}

\theoremstyle{plain}
\newtheorem{theo}{Theorem}[section]
\newtheorem{lem}[theo]{Lemma}
\newtheorem{prop}[theo]{Proposition}

\newtheorem{cor}[theo]{Corollary}

\newtheorem{lemma}[theo]{Lemma}
\newtheorem{sublemma}[theo]{Sublemma}

\theoremstyle{definition}
\newtheorem{rem}[theo]{Remark}

\newtheorem{definition}[theo]{Definition}

\newenvironment{pf}{\noindent{\it Proof.\,}}{\hfill $\square$\medskip}

\newcommand{\beq}{\begin{equation}}
\newcommand{\eeq}{\end{equation}}
\renewcommand{\a}{\alpha}
\renewcommand{\b}{\beta}

\renewcommand{\d}{\delta}

\newcommand{\f}{\varphi}
\newcommand{\g}{\gamma}
\newcommand{\h}{\eta}

\renewcommand{\k}{\kappa}

\renewcommand{\o}{\omega}

\renewcommand{\r}{\rho}
\newcommand{\s}{\sigma}
\renewcommand{\t}{\tau}

\newcommand{\D}{\Delta}

\renewcommand{\O}{\Omega}

\newcommand{\bR}{\mathbb{R}}

\newcommand{\bN}{\mathbb{N}}

\newcommand{\bB}{\mathbb{B}}

\newcommand{\gc}{\mathfrak{c}}
\newcommand{\gd}{\mathfrak{d}}

\newcommand{\gh}{\mathfrak{h}}

\newcommand{\gm}{\mathfrak{m}}
\newcommand{\gn}{\mathfrak{n}}

\newcommand{\gu}{\mathfrak{u}}

\newcommand{\gC}{\mathfrak{C}}

\newcommand{\gL}{\mathfrak{L}}

\newcommand{\cA}{\mathscr{A}}
\newcommand{\cC}{\mathcal{C}}

\newcommand{\cE}{\mathscr{E}}

\newcommand{\cG}{\mathscr{G}}

\newcommand{\cK}{\mathscr{K}}
\newcommand{\cL}{\mathscr{L}}

\newcommand{\cN}{\mathscr{N}}

\newcommand{\cP}{\mathscr{P}}
\newcommand{\cQ}{\mathscr{Q}}

\newcommand{\cS}{\mathscr{S}}

\newcommand{\p}{\partial}

\renewcommand{\square}{\kern1pt\vbox
{\hrule height 0.6pt\hbox{\vrule width 0.6pt\hskip 3pt
\vbox{\vskip 6pt}\hskip 3pt\vrule width 0.6pt}\hrule height0.6pt}\kern1pt}

\newcommand{\wt}{\widetilde}

\newcommand{\wh}{\widehat}

\newcommand{\bt}{\begin{theo}\ \ }
\newcommand{\et}{\end{theo}}
\newcommand{\bp}{\begin{prop}\ \ }
\newcommand{\ep}{\end{prop}}
\newcommand{\bc}{\begin{cor}\ \ }
\newcommand{\ec}{\end{cor}}
\newcommand{\bl}{\begin{lem}\ \ }
\newcommand{\el}{\end{lem}}
\newcommand{\bd}{\begin{definition}}
\newcommand{\ed}{\end{definition}}
\newcommand{\be}{\begin{equation}}
\newcommand{\ee}{\end{equation}}

\def\<#1,#2>{\langle\,#1,\,#2\,\rangle}
\newcommand{\arr}{\begin{array}{rlll}}
\newcommand{\ea}{\end{array}}
\newcommand{\bea}{\begin{eqnarray}}
\newcommand{\eea}{\end{eqnarray}}
\newcommand{\bean}{\begin{eqnarray*}}
\newcommand{\eean}{\end{eqnarray*}}

\newcommand{\sC}{\mathsf{C}} 
 
\newcommand{\sH}{\mathsf{H}} 
\newcommand{\sQ}{\mathsf{Q}} 
\newcommand{\sfs}{\mathsf{s}}

\renewcommand{\=}{:=}

\newcommand{\dist}{\operatorname{dist}}
\newcommand{\ve}{\varepsilon}

\newcommand{\needle}{\mathcal{N}\hskip-2pt\text{\it eedle}}
\newcommand{\cneedle}{\widecheck{\mathcal{N}}\hskip-2pt\text{\it eedle}}

  \newcommand{\vertiii}[1]{{\left\vert\kern-0.25ex\left\vert\kern-0.25ex\left\vert #1 
    \right\vert\kern-0.25ex\right\vert\kern-0.25ex\right\vert}}

\newcommand{\diam}{{\operatorname{diam}}}

\def\sideremark#1{\ifvmode\leavevmode\fi\vadjust{
\vbox to0pt{\hbox to 0pt{\hskip\hsize\hskip1em
\vbox{\hsize3cm\tiny\raggedright\pretolerance10000
\noindent #1\hfill}\hss}\vbox to8pt{\vfil}\vss}}}

\title[On the Pontryagin Maximum Principle]
{On the Pontryagin Maximum Principle\\ under differential constraints of higher order}
  \author[Franco Cardin, Cristina Giannotti and Andrea Spiro]{Franco Cardin \qquad Cristina Giannotti \qquad  Andrea Spiro}
  
\thanks{{\it Acknowledgments}. This research was partially supported by the Projects MIUR ``Regular and stochastic behaviour in dynamical systems''   by GNFM  of INdAM}

\begin{document}

\begin{abstract}  Exploiting our previous  results on higher order controlled  Lagrangians
in [Nonlinear Anal.  {\bf 207} (2021), 112263],  we 
derive here an  analogue of the classical  first order  Pontryagin Maximum Principle  (PMP)
for   cost  minimising   problems  subjected to  higher order differential  constraints 
$\frac{d^k x^j}{dt^k} = f^j\big(t, x(t), \frac{d x}{dt}(t), \ldots, \frac{d^{k-1}  x}{dt^{k-1}}(t), u(t)\big)$,  $t \in [0,T]$, 
where $u(t)$  is a control curve in a  compact  set $K \subset \bR^m$.  This result and its proof can be  considered as a detailed  illustration of one of the claims of  that previous paper,  namely  that the results of that paper,  originally established in a  smooth differential geometric framework,   yield directly  properties  holding under much weaker and more common assumptions.
In addition, for  further clarifying  our  motivations,  in the last section we display a couple of quick indications on how  the  two-step approach  of this paper  (i.e., a preliminary easy-to-get differential geometric discussion followed by a refining analysis  to    weaken the regularity assumptions) might be fruitfully exploited also  in the context  of control problems governed by partial differential equations or in studies on  the dynamics of controlled mechanical systems.
 \end{abstract}
\maketitle

\setcounter{section}{0}
\section{Introduction}
 In our previous paper \cite{CGS},  we considered the notion of controlled Lagrangians of higher order and, developing  the differential geometric   approach proposed in  \cite{CS},  we   proved  that  -- under certain strong regularity assumptions -- a  generalised version of the classical  Pontryagin Maximum Principle  (PMP) holds for control problems with higher order constraints of Euler-Lagrange type. Roughly speaking, on the one hand the results of \cite{CGS} can be considered  as  generalisations  to   controlled  Lagrangians  of arbitrary order  of   certain facts on  the first order Lagrangians, which were  established by Ioffe and Tihomirov in their elementary  proof of   the PMP   \cite{IT}.  On the other hand, 
 the results in \cite{CGS}  stem from  a fresh   differential geometric approach to  variational problems --  rooted in  Stokes' Theorem -- which  
we think will reveal to be a fruitful addition  to  the traditional Hamiltonian   tool box  of  control theory.  In fact, our differential geometric approach,  involving   controlled Lagrangians rather than controlled Hamiltonians,   admits straightforward generalisations in settings  where  several independent variables are involved  \cite{Li} and/or  lead to  applications  of the  Noether Theorem on differential constraints with symmetries \cite{Ol}.  We therefore consider it much better suited than the traditional Hamiltonian approach  for dealing with   control problems  in   Continuum Mechanics  (where  systems are governed by partial differential equations)   or for  discussions on the dynamics  of controlled systems, where information on symmetries and conservation laws can be exploited.    \par
   \smallskip
In this paper we use  the above mentioned    results of \cite{CGS}  to  derive {\it in a direct way} an  analogue of the classical  PMP 
 for     cost  minimising problems  subjected to  higher order differential  constraints of the form
$\frac{d^k x^j}{dt^k} = f^j\big(t, x(t), \frac{d x}{dt}(t), \ldots, \frac{d^{k-1}  x}{dt^{k-1}}(t), u(t)\big)$,  $t \in [0,T]$, 
with $u(t)$  control curve in a  compact  set $K \subset \bR^m$.  The contents of this paper can be  taken as a detailed  illustration  of a claim  we made in \cite{CGS} (see also \cite{CS}), 
namely that, despite the fact that  those results were  established  under   strong regularity assumptions on solutions and control curves, they  can be  nonetheless used to  directly  derive results  that hold under much weaker  regularity conditions.
\par
\smallskip
The  class of controlled problems  to which  this paper  is devoted is described in detail  as follows. 
 Consider the evolutions  $x = (x^i): [0,T] \to \sQ = \bR^n$  of a controlled dynamical system    on a   fixed  time interval $[0,T]$  
 and 
 the   cost  minimising problem, which is    determined   by the following three ingredients: 
  \begin{itemize}[leftmargin = 10pt]
 \item[--]   a family of  (possibly measurable) curves  $u:[0,T] \to K \subset \bR^m$, $u(t) = (u^a(t))$ with  values in a  compact set $K\subset \bR^m$, playing the role of the control curves for the  system; 
\item[--]  a system of $k$-th order differential constraints 
 of the form
\beq\label{diffe-1-ter} \frac{d ^k x^j}{d t^k} = f^j\left(t, x(t) , \frac{dx}{d t}(t) , \ldots, \frac{d ^{k-1} x}{d t^{k-1}}(t) , u(t)\right)  \eeq
 and  a   set $\cA_{\text{\it init}}$ of standard  $(k-1)$-th order initial conditions $\sfs$  such that  the following holds:   for  each pair $U = (u(t), \sfs)$  there exists   a unique associated  solution    $x^{(U)}(t)$ on $[0, T]$   satisfying \eqref{diffe-1-ter} and the  initial condition $\sfs$;  
 \item[--]  a  $\cC^1$  {\it terminal cost}  function $\mathsf C =  \mathsf C\left(x(T),  \frac{dx}{d t}(T) , \ldots, \frac{d ^{k-1} x}{d t^{k-1}}(T) \right)$, which depends  on the  $(k-1)$-th order jets at the final time $t = T$ of the solutions $x^{(U)}(t)$.
 \end{itemize}
 For a fixed $\sfs_o \in \cA_{\text{\it init}}$,   a  curve $u_o(t)$ in  $K$,   for which   $x_o(t) = x^{( u_o(t), \sfs_o)}(t)$  has  minimal  cost among the controlled curves with identical initial condition $\sfs_o$, is called {\it optimal control}.  If  $k = 1$ and $f = (f^i(t, x, u))$   is a  function    
 on a set of the form $[0, T] \times \O \times K$, with $\O \subset \bR^{n} $ open, which is continuous and   continuously differentiable with respect to    $x$, the described  problem of determining  optimal controls  is one of the  Mayer problems,  to which  the classical  Pontryagin Maximum Principle (PMP) applies (see e.g. \cite{PBGM, AS, BP, Ce, Ga, IT, Su3} and references therein).
Let us  briefly recall it. If $k = 1$,   given a pair $U_o = (u_o(t), \sfs_o)$,  let  us denote by   $x_o(t)  = x^{(U_o)}(t)\in \sQ$    the corresponding controlled curve and by  $p_o(t) = (p_{oi}(t))$, $t \in [0, T]$,     the  unique curve in $\sQ^* \simeq \bR^n$  satisfying the  linear differential  equations 
\beq   \label{coproblem-ter}
\dot p_j  + 
 p_i \frac{\p f^i}{\p x^j} \Big|_{(t, x_o(t), u_o(t))} = 0\quad  \text{with  the terminal condition}\quad p_{j}(T) = - \frac{\p \mathsf C}{\p x^j}\Big|_{x_o(T)}
\eeq
(here and throughout the paper we follow the Einstein convention on summations).
 Further, for any  fixed  $\t_o \in [0, T]$,  let  $ \sH = \sH^{(\sfs_o, u_o, \t_o)}: K \to \bR$  be  the  {\it Pontryagin function} defined by
 \beq  \sH^{(\sfs_o, u_o, \t_o)}(\o) \= p_{oi}(\t_o) f^i(\t_o, x_o(\t_o), \o)\ .\eeq
The classical  PMP states that  {\it   if  $u_o(t)$ is an optimal control for the considered cost problem, then    $ \sH^{(\sfs_o, u_o, \t_o)}(u_o(\t_o)) = \max_{\o \in K}  \sH^{(\sfs_o, u_o, \t_o)}(\o)$ for  almost all choices of  $\t _o\in [0,T]$}.  In many  interesting cases  this famous necessary condition   is  so    restrictive   that    can be used to completely  determine the     optimal controls.  
 \par
 \smallskip
It is quite simple to check  that,   given a pair $U_o = (u_o(t), \sfs_o)$,  at  a fixed choice of   time 
$ t_o \in [0, T]$,    the  equations which give the differential constraints  on the  $x$  and  the associated  system on the  $p$ coincide with the  Euler-Lagrange equations at $t_o$ of the (degenerate) first order Lagrangian  
 $$\cL^{(u_o, t_o)}(t, x, \dot x, p) \=  p_j (\dot x^j - f^j(t, x, u_o(t_o)))
  \ .$$
  The function   $L(t, x, \dot x, p, u) = p_j (\dot x^j - f^j(t, x, u))$,  which gives the Lagrangians $\cL^{(u_o, t_o)}$,  is called the {\it controlled Lagrangian} of  the  considered control problem (\cite{CS, CGS}). We also  recall that for any fixed $\t_o \in [0, T]$, the  Pontryagin function $ \sH = \sH^{(\sfs_o, u_o, \t_o)}$  is maximised exactly  where the function  
 \beq   \cP^{(\sfs_o, u_o, \t_o)}(\o) \= - L\big(\t_o, x_o(\t_o), \dot x_o(\t_o), p_o(\t_o), \o\big)\eeq 
  is maximised. 
 Indeed,    the difference $ \cP^{(\sfs_o, u_o, \t_o)}(\o) - \sH^{(\sfs_o, u_o, \t_o)}(\o)$ is  equal to $- p_{oj}(\t_o)  \dot x^j_o(\t_o)$, a real number  which  is  independent of  $\o$. This remark is at the basis of the so-called {\it Lagrangian version} of the PMP (see e.g. \cite{IT}). 
   \par
 \smallskip
We may now present the higher order   analog  of all this.  Consider a  control problem  with differential constraints  of the form \eqref{diffe-1-ter}  and a cost function $ \sC = \mathsf C\left(x(T), \frac{dx}{d t}(T) , \ldots, \frac{d ^{k-1} x}{d t^{k-1}}(T) \right)$. Assume  that  $f  = (f^i)$ 
satisfies the  following conditions: 
\begin{itemize}[leftmargin = 20pt]
\item[($\a$)]  it  is   $\cC^{k-1}$; 
\item[($\b$)] its partial derivatives  of order $k-1$  are  continuously differentiable  with respect to   each   argument different from $t$ and $u = (u^a)$
\end{itemize}
(note that  $(\a)$ and  $(\b)$ reduce to  the   assumptions of  the classical PMP when $k =  1$). \par
  Then,  
 consider   $n$ auxiliary dual  variables $p = (p_i) \in \sQ^* \simeq \bR^n$  and
  the controlled Lagrangian of order $k$ on the  jets  of curves in $\sQ \times \sQ^* \simeq \bR^{2n}$  defined by 
\beq \label{lag-ter}  L\left(t, x , \frac{dx}{d t} , \ldots, \frac{d ^{k} x}{d t^{k}} , p, u\right) \= p_j\left(\frac{d^k x^j}{dt^k} -  f^j\Big(t, x , \frac{dx}{d t}, \ldots, \frac{d ^{k-1} x}{d t^{k-1}} , u\Big) \right) \ .  
\eeq 
As before, for each  fixed pair $U_o = (u_o(t), \sfs_o)$, given by a control curve $u_o(t)$ and an  initial condition $\sfs_o$, we denote by  $(x_o(t), p_o(t)) \in \sQ \times \sQ^*$   the unique curve  which solves  the Euler-Lagrange equations of the Lagrangian  
$$\cL^{(u_o, t_o)}\left(t, x, \frac{dx}{d t} ,\ldots, \frac{d ^{k} x}{d t^{k}} , p\right) \=  L\left(t, x, \frac{dx}{d t} ,\ldots, \frac{d ^{k} x}{d t^{k}} , p, u_o(t_o)\right)$$
 at each  fixed $t_o \in [0, T]$  (they are explicitly given in    \eqref{diffe-1-ter}  and   \eqref{controlledEL-1bis*})  and  satisfying  the  end-point conditions defined as follows.  
The curve  $x_o(t)$ satisfies the   initial condition  at $t = 0$ given by  $\sfs_o$,  while    $p_o(t)$ satisfies  the   terminal conditions at $t = T$ listed below (here,    $x_{(s)}$ stands for  the $s$-th order derivative   $x_{(s)} \= \frac{d^s x}{d t^s}$ and $\frac{D}{Dt}$ denotes   the  total differential derivative  with the term in  $\frac{d u}{dt}$  removed -- see \eqref{22-bis} for the explicit definition):
\begin{align}
\label{31-ter} & p_i|_{t = T} =\left. \frac{\p \sC}{\p x^i_{(k-1)}}\right|_{
j^{k-1}_{T}(x_o)}\ , \\
\nonumber &\qquad  \vdots \\
& \frac{ d^\ell p_{i}}{dt^\ell}\bigg|_{t = T} =  \left.\left( (-1)^{\ell}\frac{\p \sC}{\p x^i_{(k  -1- \ell)}} +
 \sum_{h = 0}^{ \ell  -1 } (-1)^{\ell + h}  \frac{D^{h}}{ Dt ^{h} } \left( p_m\frac{\p f^m}{\p x^i_{(k - \ell + h )}}\right)   \right)\right|_{
j^{2k-3}_{T}(x_o)}\ ,\\
\nonumber  &\qquad  \vdots\\
 \label{33-ter} & \frac{d^{k-1} p_{i}}{d t^{k-1}}\bigg|_{t = T} =  \left.\left( (-1)^{k-1}\frac{\p \sC}{\p x^i} +
 \sum_{h = 0}^{ k  -2 } (-1)^{k - 1 + h}  \frac{D^{h}}{Dt ^{h} } \left( p_m\frac{\p f^m}{\p x^i_{(h +1)}}\right)   \right)\right|_{
j^{2k-3}_{T}(x_o)}\ .
\end{align}
We remark that, for each curve $x_o(t)$ solving the equations  \eqref{diffe-1-ter}, the corresponding differential problem on the  $p_j(t)$ is meaningful provided that all derivatives of $x_o(t)$ up to order $2 k - 2$ are at least  almost everywhere  defined.  This is trivially true when  $k = 1$ and $u_o(t)$ is measurable, but when $k > 1$  other assumptions are  needed. 
  For  this reason,  we impose the following condition, which  implies the desired property for any $k$ (see  Lemma \ref{lemma32}):
\begin{itemize}[leftmargin = 20pt]
\item[($\g$)]  $u_o(t)$ is measurable and,  if $k  > 1$, it satisfies the additional requirements:
\begin{itemize}[leftmargin = 5pt]
\item it  is piecewise $\cC^{k -1}$; 
\item each derivative $u_{o(\ell)}(t)$, $1 \leq \ell \leq k-1$,   takes values in a fixed compact set $K^{(\ell)} \subset \bR^m$.
\end{itemize}
\end{itemize}
Finally,   for any fixed $\t_o \in [0, T]$,  we define   
 \beq\label{gP-ter}  \sH^{(\sfs_o, u_o, \t_o)}: K \to \bR\ ,\quad  \sH^{(\sfs_o, u_o, \t_o)}(\o) \= p_{oj}(\t_o) f^j\left(\t, x_o(\t_o),  \ldots,x_{o(k-1)}(\t_o) , \o\right)\ .\eeq 
We may now state  the result we are interested in:
 \begin{theo} \label{main1} Assume that $f$ satisfies $(\a)$ and $(\b)$ and  that $u_o(t)$    satisfies  $(\g)$. Then $u_o(t)$ is an optimal control  
 only if  for  almost all  $\t_o \in [0,T]$
 \beq \label{main11} \sH^{(\sfs_o, u_o, \t_o)}(u_o(\t_o)) = \max \big\{\,\sH^{(\sfs_o, u_o, \t_o)}(\o)\ ,\ \o \in K\,\big\}\ .\eeq
 \end{theo}
 If   $f$ is  $\cC^\infty$  and the family of the control curves in which $u(t)$  is allowed to vary is assumed to be the class of the $\cC^\infty$ curves, the proof of \eqref{main11} is very short and elementary -- it is essentially a  direct  consequence of  \cite[Cor.7.7]{CGS}.  \par
 \smallskip
 An alternative way   (but undoubtfully {\it not direct})  to prove  Theorem \ref{main1} demands a  translation  of the whole setting into a corresponding  problem with first order constraints, to which the classical PMP can be  applied. More precisely, one first needs to  introduce $n \times (k-1)$   auxiliary variables, say $y_1 = (y^i_1)$, \ldots, $y_{k-1} = (y^i_{k-1})$,  and translate the   constraints  \eqref{diffe-1-ter}   into  the system of  first order constraints
 \beq\label{diffe-1-quater} 
 \begin{split}
 & \frac{d x^i}{dt} = y_1^i\ , \quad  \frac{d y^i_1}{dt} = y_2^i\ , \quad  \frac{d y^i_2}{dt} = y_3^i\ , \quad \ldots\quad \frac{d y^i_{k-2}}{dt} = y^i_{k-1}\ ,\\
& \frac{d y^i_{k-1}}{d t} = f^i\left(t, x(t) , y_1 , \ldots, y_{k-1} , u(t)\right) \ . 
 \end{split}
 \eeq
Then,  by introducing  another set of additional   variables 
$$\wt p_0 = (\wt p_{0|i})\ ,\quad \wt p_{1} = (\wt p_{1| i})\ , \ldots , \ \wt p_{k-1} = (\wt p_{k-1|i})\ ,$$
 (dual to  the variables $y_0=  (y_0^i = x^i)$, $y_1 = (y^i_1)$, \ldots , $y_{k-1} = (y^i_{k-1})$),   one may considers  an appropriate   Pontryagin function  $\sH'{}^{(\sfs_o, u_o, \t_o)}$ and a set of  necessary conditions  on  optimal controls, which are directly implied   by   the  classical  PMP with first order constraints.    Such conditions  turn out to be different from those of   Theorem  \ref{main1}. In fact,  they involve  a larger number auxiliary variables   in the definition   of the testing  function \eqref{gP-ter}.  Nonetheless, with    little additional effort,   one can come  back  to the  differential constraints on the  auxiliary functions and replace  each of them in an appropriate way. At the end one   reaches a 
 testing  function involving   a minimal  set of auxiliary variables, namely  the above defined function \eqref{main11}. \par
 \smallskip
The main purpose of this paper is  to show that, instead of adopting the above mentioned back and forth type argument  (namely,   introducing  new auxiliary  variables with the purpose  of  completely removing them in a second time)  Theorem \ref{main1} can be proved in a direct way,  by just extending  to  low regularity settings  the simple   differential geometric arguments  -- based on  Stokes' Theorem and of the Principle of Minimal Labour \cite{CS} --  that  hold when  there are no restrictions on the level of  regularity.  In this paper  the extensions  to low regularity settings  are reached  by standard approximation techniques (similar to those  considered  by Gamkrelidze for first order control problems, see e.g. \cite{Ga}) together with a couple of ad hoc  lemmas  for  estimating  variations of costs.  All  arguments are in principle quite simple but  demand  a number of tedious  checks.  Aiming to make the  proof  self-contained  and ready  to  be used in  future works, all proofs are presented    in  great detail.   This is   done also with the hope  the proof can be taken as    a convincing illustration that  an   approach to control problems, structured into 
\begin{itemize}[leftmargin = 10pt]
\item  a preliminary  easy-to-get  differential geometric proof, made 
under strong regularity assumptions, followed by  
\item a second   step  in which the  previous  claims are improved to new  statements that are true under  a minimal set of regularity assumptions, 
\end{itemize}
can be fruitfully exploited in   much more involved  contexts of  control theory.  A couple of examples of settings where such   two-step  approach  might provide interesting results  are given in the last section.   
\par
\begin{rem}  The approach we followed for our proof of  Theorem \ref{main1} made us  aware  that the minimal  (or, more precisely,  very close to the minimal) set of regularity  conditions  to be imposed 
is  that the  function  $f = (f^i)$ satisfies the  conditions $(\a)$ and $(\b)$ and the control curve  $u_o(t)$  satisfies $(\g)$.  Such  conditions become  progressively weaker  if the order of the system is reduced by introducing   auxiliary variables. This phenomenon can be synthesised by saying that 
 {\it    the stronger are the  regularity properties of the differential constraints, the fewer  auxiliary variables need to be considered}. This    property  was  unexpected to us and can  be  taken as an interesting by-product of our  approach.
 \end{rem}
The paper is structured as follows.  After  a preliminary section, in \S\ \ref{definingtriples} we review certain basic notions and   results  given in \cite{CGS}.   In \S\ \ref{key}  we present the two main  lemmas that  allow  using  approximations to  extend the results of \cite{CGS}  to  problems with     weaker regularity assumptions.  The proofs of Theorem \ref{main1} and of  parts of the two main lemmas are given  in \S\ \ref{concluding} and \S\ \ref{theremaining}, respectively.  Suggestions for  further investigations  are  given in  \S \ref{added}.\par
\noindent {\it Acknowledgements.} We are sincerely grateful to the referee for her/his accurate reading and   very nice and useful suggestions  to improve the presentation.  
\par
\vskip  1 cm
\section{Preliminaries}
\label{preliminaries}
\subsection{Notational issues} \label{notationalissues}
The standard coordinates of the  ambient   space $\cQ = \bR^N$  are usually denoted by $q = (q^i)$. When  $\cQ$ is even dimensional  and of   the form  $\cQ = \sQ \times  \sQ^*$ for some $\sQ =  \bR^{n}$,  the   coordinates are  indicated as pairs   $q = (x, p)$,  with  $ x = (x^i) \in \sQ$ and  $p = (p_j) \in \sQ^*$.  The elements of the  controls set  $K \subset \bR^m$ are  $m$-tuples $u = (u^a)$.
Given a $k$-times differentiable curve  $q(t)\in  \cQ$, $t \in I \subset \bR$,  its $r$-th order derivative  is  often indicated with  the short-hand notation
$$q_{(r)} (t)   \=  \left(q^i_{(r)} (t)\right) = \left( \frac{d^r q^i}{d t^r}\bigg|_{t} \right)\ .\qquad \text{We also set}\qquad q_{(0)}(t) \= q(t)\ .$$
Accordingly,   the $r$-th order jet of   a curve  $ \g(t) \=  (t, q(t))  \in \bR \times \cQ$ is denoted   by    
\beq \label{jetcoord} j_t^r(\g) = (t, q(t), q_{(1)}(t), \ldots, q_{(r)}(t)) =  (t, q_{(s)}(t))  \ .\eeq
 The manifold   
of all   $r$-th order jets  of curves, i.e.  the  {\it jet bundle of order $r$} of the  (trivial) bundle $\bR \times \cQ$ over $\bR$,    is denoted by  $J^r(\cQ|\bR)$.  \par
If $g:  J^{r}(\cQ|\bR) \times K\to \bR$ is a smooth  function of pairs $(j^r_t(\g), u)$, formed by  a $r$-th order jet of  a curve  $ \g(t) \=  (t, q(t))$ and  a control $u \in K$, we denote by 
$\frac{D g}{ Dt} $ the smooth real function  on  $J^{r+1}(\cQ|\bR) \times K$ defined by 
\beq  \label{22-bis}
\frac{D g}{D t} \bigg|_{(j^{k+1}_t(\g), u)}\=   \frac{\p g}{\p t}\bigg|_{(j^{r}_t(\g), u)} + \sum_{\ell = 0}^{r}\sum_{i=1}^{N} q^i_{(\ell + 1)} \frac{\p g}{\p q^i_{(\ell)}}\bigg|_{(j^{r}_t(\g), u)} \ .\eeq
Note that   $\frac{D}{Dt}$  differs from  the  total derivative for  the sections of    $\pi: J^{r}(\cQ |\bR) \times K\to \bR$ only because the term  $u^a_{(1)}  \frac{\p}{\p u^a}$   is  missing. We can also say that  the operator $\frac{D}{Dt}$ is the pull-back on $ J^{r}(\cQ|\bR) \times K$ of the classical total derivative operator  $\frac{d}{d t}$  of  $J^{r}(\cQ|\bR)$.
\par
 \subsection{Generalised Mayer problems with constraints of variational type}
 As we mentioned in the Introduction, any $k$-th order control system  
 \eqref{diffe-1-ter}  on curves $x(t)$ can be  considered as a sub-system of  the     Euler-Lagrange  equations  of an appropriate controlled  $k$-th order  Lagrangian   for    curves of the form $t \mapsto (x(t), p(t)) = (x^i(t), p_j(t))$ with  the curve  $t \mapsto p(t)$ in an appropriate  auxiliary space. 
 Despite  of the fact that  \eqref{diffe-1-ter} has order $k$,    the full    Euler-Lagrange system   contains    equations of 
 order $2k$ (see \S \ref{sect51*}). Moreover,   a    solution    $t \mapsto (x(t), p(t))$, $t \in [0, T]$,   of the controlled Euler-Lagrange  equations   is   determined  not only by the  $(k-1)$-th order jet  $\s^{(x)}$ of  $x(t)$ at $t = 0$ and by the control curve  $u(t)$, but also  by the  $(2k-1)$-th order jet $\s^{(p)}$ at $t = 0$  of the curve  $p(t)$.  Since  in the curve $ (x(t), p(t))$ only the part  $x(t)$   is relevant for  the  cost problem,  the part  of the initial datum $\s$ for the curve  $ (x(t), p(t))$,  which is not  determined by $x(t)$,  is    freely specifiable and can be considered  as an additional   ``controlling datum''   for the problem on the curves $(x(t), p(t))$.  
In     \cite{CGS}  we developed   a  theory of     control problems with  (smooth) differential constraints of variational type,   which     not only  works for   the classical first order     Mayer problems  and the cost problems of this paper, but it is  designed to be  applicable   to   other contexts. With  such second aim in mind  and on the basis of the previous observation,  in  \cite{CGS} we were naturally led to consider the following definition (see also \cite{CS}).\par
 \medskip
A {\it  generalised Mayer problem with smooth constraints of variational type  of order $k$} (for short, {\it generalised Mayer problem}) is a  cost minimising problem determined by a  triple  $(\cK, L, C)$  of the following kind.
 \begin{itemize}[itemsep=5pt, leftmargin=10pt]
\item $\cK$ is a {\it set of   control pairs} $U = (u(t), \s)$,   given by:
 \begin{itemize}
\item[a)]    a  
smooth curve 
$u: [0,T]   \to K \subset \bR^m$    with values in a    compact subset $K \subset \bR^m$; 
\item[b)]   a $(2k -1)$-jet   $\s = j^{2k-1}_{t = 0}(\g)$   at $t = 0$ of a smooth curve $t \mapsto \g(t)  = (t, q^i(t))$. 
\end{itemize}
 The jets $\s$  of these  pairs are constrained to be elements of  a fixed    set $\cA_{\text{\it init}}\neq \emptyset$.
\item  A {\it (smooth) controlled Lagrangian $L = L(t, q_{(s)}, u)$  of   order $k$}, i.e. a $\cC^\infty$ function 
$$L : J^k(\cQ|\bR) \times K \to \bR$$
for which the following property holds:  for each control pair  $U = (u(t), \s) \in  \cK$, there exists exactly one curve $\g(t) = (t, q^i(t))$, $t \in [0, T]$,   with  $j^{2k-1}_{t = 0}(\g) = \s$
and satisfying the  
Euler-Lagrange equations of the higher-order Lagrangian $\cL^{(u, t_o)}(t, q_{(\b)}) \= L(t, q_{(\b)}, u(t_o))$  at each time $t_o \in T$, i.e.    the    equations 
 \begin{multline} \label{controlledEL} E_i(L) |_{(j^{2k}(\g(t)), u(t))} \= \\
 = \frac{\p L}{\p q^i} \bigg|_{(j^{2k}(\g(t)), u(t))} + \sum_{\b= 1}^k (-1)^\b \left(\frac{D}{ D t}\right)^\b\left(\frac{\p L}{\p q^i_{(\b)}}\right) \bigg|_{(j^{2k}(\g(t)), u(t))} \hskip -1 cm = 0\ ,\quad i = 1, \ldots, N\ , \end{multline} 
where  $\frac{D}{ D t}$ is the operator defined  in \eqref{22-bis}. This curve   is denoted $\g^{(U)}(t)$.  
\item A {\it (smooth) terminal cost function} $C: J^{2k-1}_{t = T}(\cQ|\bR) \to \bR$. 
 \end{itemize}
  The  equations \eqref{controlledEL}  are called {\it differential constraints} of  the triple $(\cK, L, C)$ and the curves $\g^{(U)}$,  with $U = (u(t), \s)) \in \cK$, 
 are called {\it $\cK$-controlled curves}. \par
A triple  $(\cK, L, C)$ as above   is  called a {\it defining triple}. The problem  of determining  the  $\cK$-controlled curves  $\g^{(U_o)}$, $U_o \in \cK$, for which the  terminal cost  $C(j^{2k-1}_{t = T}(\g^{(U_o)}))$  
 is  minimal,   is called the   {\it  generalised Mayer problem} corresponding to   the triple.  The  pairs  $U_o = (u_o(t), \s_o)$ giving the   cost minimisation curves are called {\it optimal controls}.
\par
\smallskip
In  the next section \S \ref{sect51*}, we  explain  how these notions are well fitted  with the cost   problems  described in the Introduction. 
\par
\smallskip
Throughout the paper  we  restrict our discussion to the   defining triples $(\cK, L, C)$ satisfying the following additional technical hypothesis.  We assume that  
there exists  an  open {\it convex} superset  $\wh K \supsetneq K$ such that,  denoting by   $\wh \cK \supsetneq \cK$ 
  the family  of pairs $U = (u(t), \s)$ with   $u(t)$   in  $\wh K$ but initial condition   $\s $  as in  $ \cK$,   there still exists a   unique solution  to    \eqref{controlledEL}  for any   $ U   \in \wh \cK$. The curves of this larger family   are   called {\it $\wh \cK$-controlled}. 
\par
\medskip
\subsection{Defining triples for  the   cost  problems  of this paper} \label{sect51*} Assume that  all data of  one of the  cost minimisation problem of   the Introduction are   of class $\cC^\infty$.  We claim that in this case it is equivalent to the  generalised Mayer problem given by the following defining triple $(\cK, L, \sC)$.
Consider the  configuration space   $\cQ = \sQ \times \sQ^*$ with $\sQ = \bR^n$ and coordinates    $(x, p) = (x^i, p_j)$,  as specified in \S \ref{notationalissues}. The controlled Lagrangian $L$ and the  cost function $\sC$ are the maps
\begin{multline} \label{lag}  L(t, x, x_{(1)}, \ldots, x_{(k)},  p, u) = p_j \bigg(x^j_{(k)} - f^j(t, x, x_{(1)}, \ldots, x_{(k-1)}, u)\bigg) \qquad \text{and}  \\  C = \sC\left(x(T), x_{(1)}(T), \ldots,x_{(k-1)}(T)\right)\ ,  \end{multline}
where $f = (f^j)$ is the function that gives the  constraints \eqref{diffe-1-ter} and $\sC$ is the  cost of the Introduction, which depends only on the  $(k-1)$-th order jet at $t = T$  of  the part $x(t)$ of the $\g(t) = (t, x(t), p(t)) \in [0, T] \times \sQ \times \sQ^*$.  Note that, even if   $C$ depends only on the jets  of order $k-1$, it can be trivially considered as a function on the space of $(2k-1)$-jets  and   hence is subsumed by  the   general notion  of  terminal  cost function  considered    in \cite{CGS}.  Finally, the set  $\cK$ consists of the pairs 
 $U = (u(t), \s = \big(\mathsf s, \wt{\mathsf s})\big) $  where: (a)     $\mathsf s$ is  an initial condition  in a prescribed    set $\cA_{\text{\it init}} \subset J^{k-1}(\sQ|\bR)|_{t = 0}$  for the curve $x(t)$; (b)    $\wt{\mathsf s}$ is an initial condition (which  can be   freely chosen)  for the curve $p(t)$. \par
 \smallskip
 To show that the problem determined by $(\cK, L, \sC)$ is equivalent to our original cost minimising  problem, we  first observe that the  differential constraints determined by $(\cK, L, \sC)$ are  given by  the controlled Euler-Lagrange equations (which  are of normal type)
 \begin{multline} \label{controlledEL-1*} E_{p_i}(L) |_{(t, j_t^{2k}(\g), u(t))} \= \\
 = \frac{\p L}{\p p_i} \bigg|_{(t, j^{2k}_t(\g), u(t))} + \sum_{\b= 1}^k (-1)^\b \left(\frac{D}{ Dt}\right)^\b\left(\frac{\p L}{\p p_{(\b)i}}\right) \bigg|_{(t, j^{2k}_t(\g), u(t))} = \\
 =   x^i_{(k)} -  f^i(t, x^i, x^i_{(1)}, \ldots, x^i_{(k -1)}, u^a) = 0\ , \end{multline}
 \begin{multline} \label{controlledEL-1bis*} E_{x^i}(L) |_{(t,j^{2k}_t(\g), u(t))} \= \\
 = \frac{\p L}{\p x^i} \bigg|_{(t, j^{2k}_t(\g), u(t))} + \sum_{\b= 1}^k (-1)^\b \left(\frac{D}{D t}\right)^\b\left(\frac{\p L}{\p x^i_{(\b)}}\right) \bigg|_{(t, j^{2k}_t(\g), u(t))} = \\
 =   (-1)^{-k}
\left\{ p_j \frac{\p f^j}{\p x^i} - \frac{D}{ Dt} \left( p_j \frac{\p f^j}{\p x^i_{(1)}}\right)+  \frac{D^2}{Dt^2} \left( p_j \frac{\p f^j}{\p x^i_{(2)}} \right) + \right.\hskip 1.8 cm \\
\left.+ \ldots + (-1)^{k -1}  \frac{D^{k -1}}{Dt^{k -1}} \left( p_j \frac{\p f^j}{\p x^i_{(k -1)}} \right)\right\}\bigg|_{(t, j^{2k-2}_t(\g), u(t))} -  p_{(k)i} = 0\ .  \end{multline}
We immediately see that   \eqref{controlledEL-1*}  coincides with   the  differential constraints \eqref{diffe-1-ter}. This  fact together with   the fact  that   $\sC$ is independent of  $p(t) = (p_j(t))$   implies the claimed  equivalence between our cost minimising   problem and the   problem determined by $(\cK, L, \sC)$ -- at least in the case  of smooth data. \par
In the next sections, we   establish some preliminary results on  generalised Mayer problems. We will come back to  this  specific  one    in \S \ref{concluding}.
\par
\subsection{Differential constraints of normal type}
The smooth higher order constraints  \eqref{controlledEL}  of a generalised Mayer problem  are called   {\it of normal type} if, using  an appropriate number of auxiliary variables, they   can be reduced to a    first order system of the type 
\beq \label{first-3}  \frac{d y^A}{d t} = g^A(t, y^B,  u^a(t)) \ , \qquad 1 \leq A \leq \wt N\ ,\eeq
where  the   $g^A$ are   functions uniquely determined  by the   $L$ of the defining triple.  
For such constraints,  we  denote   by $\wh \cK_{\operatorname{meas}}$ ($\supsetneq \wh \cK$) the   set  of the pairs 
$U =  (u(t), \s)$,  in which  the control curve  $u(t)$ in $\wh K$  is merely measurable. The corresponding   curves $\g^{(U)}$  (they exist by  well known facts on  first order equations) are  called {\it $\wh \cK_{\operatorname{meas}}$-controlled}. 
  \par

Some well-known properties  of  controlled first order differential equations (see e.g. \cite[Ch. 3]{BP})  have  useful direct counterparts  for the higher order  constraints of this kind. We collect them in the next lemma, where
we prove a Gronwall-type result  for  generalised Mayer problems of normal type,  namely  that  two controlled curves are close  whenever their initial conditions   and controls are close. In the statement
 the following notation is  used. For any pair of measurable control curves  $u, u': [0,T] \to \wh K$
 we  denote 
 by $\operatorname{dist}(u, u')$    the distance 
\beq \label{distanza} \operatorname{dist}(u, u') \=  \mu_{\text{Leb}}\big(\left\{\ t \in [0, T]\ :\ u(t) \neq u'(t)\  \right\}\big)\ ,\eeq
where  $ \mu_{\text{Leb}}$ is the Lebesgue measure  on the subsets of $[0, T]$.  
We also use the jets coordinates  \eqref{jetcoord} to  identify  $J^{2k-1}(\cQ|\bR)|_{t = 0}  \simeq \bR^{2k N + 1} $  and we  use the  classical
Euclidean norm of $\bR^{2k N + 1}$ to  define  distances  $|\s - \s'|$ between  initial conditions   $\s, \s' \in \cA_{\text{\it init}}$.\par

\begin{lem} \label{normallemma} Let $(\cK, L, C)$ be a defining triple, giving  differential constraints of order $2k$  of normal type, i.e.  equivalent  to  first order equations
of the form  \eqref{first-3}.  Assume also that the initial conditions in  $\cA_{\text{\it init}}$  are   in bijection with a set $\wt \cA_{\text{\it init}}$ of initial conditions for the problem  \eqref{first-3} by means of a Lipschitz continuous map. \par
Given a  $\wh \cK_{\operatorname{meas}}$-controlled curve  $\g_{o}(t)= \g^{ (U_o)}(t) = (t, q_o(t))$,  with $ U_o = (u_o(t), \s_o) \in \wh \cK_{\operatorname{meas}}$, 
 there exist     constants   $\r, \k, \gC$,  $\gC' >0$ such that 
   for any    $U  = (u(t), \s)$, $U'= (u'(t), \s')  \in  \wh \cK_{\operatorname{meas}}$  with 
   $$|\s - \s_o|,\  |\s' - \s_o| < \r\ ,\qquad \dist(u , u_o), \ \dist(u',  u_o) < \rho\ ,$$ 
    the corresponding  curves
$\g^{(U)}, \g^{(U')} : [0, T] \to [0,T] \times \cQ$ 
are such that 
\beq\label{const}  \| \g^{(U)} - \g^{(U')}\|_{\cC^{2k-1}} \leq \gc\, \operatorname{dist}(u, u') + \k |\s - \s'|\qquad\text{with}\  \ \gc \= 4  \gC e^{2 \gC' T} \ .\eeq
The constants $\r$, $\k$,$\gC$,   $\gC'$  depend only on  
\begin{itemize}[leftmargin = 18 pt]
\item[(a)] the Lipschitz constant of the  bijection between  $\cA_{\text{\it init}}$ and $ \wt \cA_{\text{\it init}}$; 
\item[(b)] the function $g = (g^A)$ and thus the controlled  Lagrangian $L$; 
 \item[(c)] the  choice of a cut-off function $\f: \bR \times \bR^{\wt N} \to \bR$, which is identically equal to $1$ on a relatively compact  neighbourhood $\cN \subset \bR^{\wt N + 1}$ of the  trace of the curve  $\wt \g_o: [0, T] \to [0,T] \times \bR^{\wt N}$,   which solves \eqref{first-3} and  corresponds to  $\g_o(t) = (t, q_o(t))$ in $[0, 1] \times \cQ$.
\end{itemize}
\end{lem}  
\begin{pf} Let $\wt \g_o = (t, y^A_o(t))$, $\cN$ and $\f$ as in (c)  and denote by $\wt \s_o \= (y^A_o(0))$ the initial condition of $\wt \g_o$.
 The curve $\wt \g_o$ is a solution not only to    \eqref{first-3}, but also to  the system 
 \beq \label{first-31bis}  \frac{d y^A}{d t} = h^A(t, y^B, u^a(t)) \qquad  \text{where}\ \ h^A(t, y^B, u^a(t))  \=  \f(t, y^B) \ g^A(t, y^B,  u^a(t)) \ .\eeq
 By construction,  
$h|_{\cN} = g|_{\cN}$ and there are constants   $\gC, \gC' > 0$ (depending  on $g$ and $\f$) such that 
$$\sup_{(t, y, \o) \in \bR^{1 + \wt N}\times \wh K} |h^A(t, y, \o)| \leq \gC\ ,\qquad \sup_{(t, y, \o) \in \bR^{1 + \wt N}\times \wh K}  \left\|\frac{\p h^A}{\p y^B}\bigg|_{(t, y, \o)}\right\| \leq \gC'\ .$$
 By  classical arguments based on Gronwall Lemma (see e.g. \cite[Prop. 3.2.2]{BP}), for any two   solutions  $\wt \g^{(u, \wt \s)}(t) = (t,  y^{(u, \wt \s)}(t)) $ and   $\wt \g^{(u', \wt \s')}(t) = (t,  y^{(u', \wt \s')}(t)) $ of  \eqref{first-31bis}, which are    determined by  pairs $(u, \wt \s)$, $(u', \wt \s')$,  given by  measurable curves $u(t), u'(t) \in \wh K$ and initial conditions    $\wt \s, \wt \s' \in \wt \cA_{\text{\it init}}$, 
  we have that 
 $$\|y^{(u, \wt \s)} - y^{(u', \wt \s')}\|_{\cC^0} \leq \gc\, \operatorname{dist}(u, u') + |\wt \s - \wt \s'|\ \qquad \text{with} \  \gc \= 4  \gC e^{2 \gC' T}\ .$$
Therefore if
 $\operatorname{dist}(u, u_o)$, $\operatorname{dist}(u', u_o)$, $|\wt \s - \wt \s_o|$ and $|\wt \s' - \wt \s_o |$
 are sufficiently small, then   both   
 curves   $\wt{\g}^{(u, \wt \s)}$,  ${\wt \g}^{(u', \wt \s')}$ have trace in $\cN$ and  are solutions  to   \eqref{first-3}.  Since the bijection from  $ \cA_{\text{\it init}}$ to  $\wt \cA_{\text{\it init}}$  is Lipschitz,   there is a  $\k>0$ such that $|\wt \s - \wt \s'| \leq \k |\s - \s'| $ and  the lemma follows.  
\end{pf}
\par
  \begin{rem}[Stability under perturbation] \label{remark33}The  claim of  Lemma \ref{normallemma} has the following  extension, which is later used in the proof of our main result. Consider a one-parameter family of defining triples $(\cK, L^{(\d)}, C^{(\d)})$,  each of them  with the same set of control pairs $\cK$, but with  (smooth) Lagrangians and cost functions, depending on a real parameter $\d \in (0, \d_o]$.  Assume that all of such triples give generalised Mayer problems  of normal type and that the associated  equivalent first order constraints 
   \beq \label{first-31delta}  \frac{d y^A}{d t} =  g^{(\d)A}(t, y^B, u^a(t)) \ , \eeq
   are such that the functions   $g^{(\d)}$ tend uniformly on compacta, together with their  first derivatives in the  $y^B$,  to  a limit map
   $g(t,  y^B, u) = \lim_{\d \to 0} g^{(\d)}(t,  y^B, u)$. This limit map   is clearly 
 continuous and   continuously differentiable  with respect to the  $y^B$. \par
   Consider now a pair  $U_o = (u_o(t), \s_o) \in \wh \cK_{\text{meas}}$ and the uniquely associated solution  $\wt \g_o(t) = (t, y^A_o(t))$ to  
   \beq \label{first-31nodelta}  \frac{d y^A}{d t} =  g^{A}(t, y^B, u^a_o(t)) \ , \eeq
   with initial condition given by the point $\wt \s_o$  corresponding  to $\s_o$.   Considering  a neighbourhood $\cN$ and a cut-off function $\f$ as in (c) of Lemma \ref{normallemma} and  setting   
$$ h^{(\d)}(t, y^B, u^a) \=  \f(t, y^B) \ g^{(\d)A}(t, y^B,  u^a(t))\ ,\qquad h(t, y^B) \=  \f(t, y^B, u^a) \ g^A(t, y^B,  u^a(t)) \ ,$$
we have that also the  $h^{(\d)}$ and their first derivatives with respect to the $y^B$  tend uniformly on compacta to  $h$ and to  its corresponding first derivatives. 
Due to this, for any sufficiently small interval $(0, \d_o]$,  it is possible to select   $\d$-independent constants  $\r, \k, \gC$,  $\gC' >0$   such that    the claim of Lemma \ref{normallemma} holds with such  constants for  any  triple  $(\cK, L^{(\d)}, C^{(\d)})$, $\d \in (0, \d_o]$. 
   \end{rem} 
   \par
 \medskip 
\section{The generalised PMP  for Mayer problems with   smooth differential constraints of variational type}\label{Sect4}
 \label{definingtriples}
 \label{generalisedneedle} 
 Let   $(\cK, L, C)$  be a    defining  triple for a  generalised  Mayer problem  with differential constraints of normal type. Given a control pair 
 $ U_o  = (u_o(t),  \s_o) \in \cK$, with  associated  $ \cK$-controlled curve, 
and a   triple $(\t_o, \o_o, \ve) \in  (0, T) \times K \times (0, + \infty)$ with  $ 0 < \ve  < \min\left\{1, \frac{\t_o}{2}, T- \t_o\right\}$, 
we define  
\beq \label{uk}
u^{(\t_o, \o_o, \ve)}: [0, T] \to K\ ,\qquad u^{(\t_o, \o_o, \ve)}(t) \= \left\{\begin{array}{ll} u_o(t) & \text{if} \ t \in \big[0, \tau_o- {\ve} \big),\\[4pt]
 \omega_o & \text{if}\ t \in \big[\tau_o - {\ve},  \tau_o\big),\\[4pt]
 u_o(t) & \text{if} \ t \in \big[\tau_o, T\big]\ .\\
 \end{array}\right.
 \eeq 
We also select a   constant $0 < \gh <\frac{1}{2}$ (\footnote{From now till almost to the end,   such an  $\gh$ is a fixed number,  say e.g. $\gh = \frac{1}{4}$. Only  at  the very end of  \S \ref{finalproof},  where a $\d$-parameterised family of Mayer problems is taken into account, this constant  $\gh$  will be taking depending on $\d$ and tending to $0$ for $\d \to 0$.}) and  with the  (in general discontinuous) curve \eqref{uk}, we  associate  a  {\it smooth} curve $\widecheck u^{(\t_o, \o_o, \ve)}:  [0, T] \to \wh K$
satisfying the  condition  
\beq \label{theh} \widecheck u^{(\t_o, \o_o, \ve)}(t) = u^{(\t_o, \o_o, \ve)}(t) \qquad \text{for any} \ \  t \notin [\t_o - \ve - \gh \ve^2,  \t_o - \ve] \cup [ \t_o, \t_o + \gh \ve^2]\ .\eeq
We assume that the smoothing  algorithm which determines  the smooth  $\widecheck u^{(\t_o, \o_o, \ve)}(t)$  from the non-smooth $u^{(\t_o, \o_o, \ve)}$  is fixed  (the choice of the  algorithm  does not matter).  \par
 \smallskip
 We call  $u^{(\t_o, \o_o, \ve)}$  the   {\it  needle modification}  of $u_o(t)$ with  {\it peak time} $\t_o$,  {\it ceiling value $\o_o$} and   {\it  width $\ve$}. The associated smooth curve 
  $\widecheck u^{(\t_o, \o_o, \ve)}$ is called the {\it  smoothed needle modification} of  $u^{(\t_o, \o_o, \ve)}$ (see Fig. 1 and Fig. 2).  \par\smallskip
 \centerline{
\begin{tikzpicture}
\draw[<->, line width = 1] (1,3.8) to (1,0.5) to (6.5,0.5);
\draw[<->, line width = 1] (7,3.8) to (7,0.5) to (12.5,0.5);
\draw[fill]  (3, 0.5) circle [radius = 0.06];
\node at  (3, 0.2) {\tiny $\t_o - \ve$};
\draw[fill]  (4.5, 0.5) circle [radius = 0.06];
\node at  (4.5, 0.2) {\tiny $\t_o$};
\draw[fill]  (6.2, 0.5) circle [radius = 0.06];
\node at  (6.2, 0.2) {\small $T$};
\draw[fill]  (9, 0.5) circle [radius = 0.05];
\draw[fill,purple]  (8.5, 0.5) circle [radius = 0.06];
\node[purple] at  (8.1, 0.2) {\tiny $\t_o {-} \ve{-} \gh\ve^2$};
\node[purple] at  (11.2, 0.2) {\tiny $\t_o {+} \gh\ve^2 $};
\draw[fill]  (10.5, 0.5) circle [radius = 0.05];
\draw[fill,purple] (11, 0.5) circle [radius = 0.06];
\draw[fill]  (12.2, 0.5) circle [radius = 0.06];
\node at  (12.2, 0.2) {\small $T$};
\draw[fill]  (1, 3) circle [radius = 0.06];
\node at  (0.8, 3) {\small $\o$};
\draw[fill]  (7, 3) circle [radius = 0.06];
\node at  (6.8, 3) {\small $\o$};
\node[blue] at (2.3,1.7) {\tiny$u_o(t)$};
\node[blue] at (5.3,1.9) {\tiny$u_o(t)$};
\node[black] at (3.7, 2.5) {$u^{(\t_o, \o_o, \ve)}$};
\draw [line width = 0.7, blue] (1, 1) to [out=-20, in=200] (3,1.7)   ; 
\draw [line width = 0.7, blue, dashed](3,1.7)  to (3, 3)   ; 
\draw [line width = 0.7, blue] (3, 3) to  (4.5,3)   ; 
\draw [line width = 0.7, blue, dashed] (4.5,3)  to(4.5, 1.7)   ; 
\draw [line width = 0.7, blue] (4.5, 1.7) to [out=-20, in=150] (6.2,1.5)   ; 
\draw [line width = 0.7, blue] (7, 1) to [out=-20, in=200] (9,1.7)   ; 
\draw [line width = 0.7, blue, dashed](9,1.7)  to (9, 3)   ; 
\draw [line width = 0.7, blue] (9, 3) to  (10.5,3)   ; 
\draw [line width = 0.7, blue, dashed] (10.5,3)  to(10.5, 1.7)   ; 
\draw [line width = 0.7, blue] (10.5, 1.7) to [out=-20, in=150] (12.2,1.5)   ; 
\draw [line width = 1, purple](8.5,1.46)  to  [out=20, in=180] (9.1, 3) ; 
\draw [line width = 1, purple] (10.4,3)  to  [out=0, in=180] (11, 1.64)   ; 
\node[blue] at (8.3,1.7) {\tiny$u_o(t)$};
\node[blue] at (11.3,1.9) {\tiny$u_o(t)$};
\node[black] at (9.7, 2.5) {$\widecheck u^{(\t_o, \o_o, \ve)}$};
 \end{tikzpicture}
 }
   \centerline{\tiny\hskip 1 cm \bf Fig. 1 Needle modification \hskip 2 cm \bf Fig. 2 Smoothed needle modification}

The  (non-smooth and smoothed) needle modifications are  essential ingredients  for the following definition,  in which we  combine  the    classical notion of needle variation, developed by Boltyanski for the original proof of the classical PMP, and the concept   of homotopy variation.   
\par
\begin{definition}  Given a controlled curve  $\g^{( U_o)}(t)$  and a triple $(\t_o, \o_o, \ve_o)$   as above,   consider 
 a continuous map  $\Sigma:  [0, \ve_o] \times  [0, 1] \subset \bR^2 \to \cA_{\text{\it init}}  \subset J^{2k-1}(\cQ|\bR)|_{t = 0}$   such that    $\Sigma(\ve , 0) = \Sigma(0, s) = \s_o$ for any $\ve$ and $s$.
 Moreover, for any $s \in [0,1]$ and $\ve \in [0, \ve_o]$,  let us denote by
  $u^{(\ve, s)}$   (resp.  $\widecheck u^{(\ve, s)}$)  the  control curve in   the convex set $\wh K$ defined by 
\beq  u^{(\ve, s)}(t) = (1-s) u_o(t) + s  u^{(\t_o, \o_o, \ve)}(t)\ ,\qquad s \in [0,1]\ \eeq
 \beq  \left(\ \text{resp.}\  \widecheck u^{(\ve, s)}(t) = (1-s) u_o(t) + s  \widecheck u^{(\t_o, \o_o, \ve)}(t)\ ,\qquad s \in [0,1]\ \ \right)\ .\eeq
The {\it  needle variation} (resp. {\it smoothed needle variation})  {\it of  $\g^{( U_o)}$} associated with  $(\t_o, \o_o$, $\Sigma$, $\ve_o)$   is the  one-parameter family of maps
\beq\label{Cardino} \needle^{(\t_o, \o_o, \Sigma, \ve_o)} (\g^{(U_o)})\= \{\  F^{(\t_o, \o_o, \Sigma)(\ve)}:[0,T] \times [0,1] \to [0,T] \times  \cQ\ , \ \ve \in (0, \ve_o]\ \}\eeq 
\beq\label{Cardino-smooth} \left(\text{resp.} \  \cneedle^{(\t_o, \o_o,  \Sigma, \ve_o)} (\g^{(U_o)})\= \{\  \widecheck F^{(\t_o, \o_o, \Sigma)(\ve)}:[0,T] \times [0,1] \to [0,T] \times  \cQ\ , \ \ve \in (0, \ve_o]\ \}\right),\eeq 
given  by  the homotopies of controlled curves 
 $$ F^{(\t_o, \o_o, \Sigma)(\ve)}(t,s) = \g^{(U^\ve(s))}(t)\ ,\qquad U^\ve ( s) \= \big(u^{(\ve, s)}(t),   \Sigma(\ve, s)\big)$$
  $$ \Bigg(\ \text{resp.}\quad  \widecheck F^{(\t_o, \o_o, \Sigma)(\ve)}(t,s) = \g^{(\widecheck U^\ve(s))}(t)\ ,\qquad \widecheck U^\ve(s) \= \big(\widecheck u^{(\ve, s)}(t), \Sigma(\ve, s)\big)\ \Bigg)$$
 \end{definition}
The  class of   needle variations of  a fixed controlled curve  contains  the following  important subclass, 
which plays a crucial  role in the generalised PMP   established in  \cite{CGS}. \par
   Consider    a  $\cK$-controlled curve   $\g_o= \g^{( U_o)}$ and  a  smoothed needle variation 
 $ \cneedle^{(\t_o, \o_o, \Sigma, \ve_o)} (\g_o)$  as above. We introduce  the following notation.  For each homotopy  $\widecheck F^{  (\t_o, \o_o, \Sigma)(\ve)}$, $\ve \in [0, \ve_o]$,  we denote
\begin{itemize}[leftmargin = 15pt]
\item   by
 $\widecheck F^{  (\t_o, \o_o, \Sigma)(\ve)(2k-1)}$ the  homotopy of the curves in $J^{2k-1}(\cQ|\bR) \times \wh K$, given by the $s$-parameterised family of  maps 
 $t \to ( j^{2k-1}_t(\g^{(\ve, s)}) , \widecheck u^{(\ve, s)}(t))$, 
  made of  the $(2k-1)$-jets $ j^{2k-1}_t(\g^{(\ve, s)})$  of the   curves $\g^{(\ve, s)}(t) \= \widecheck F^{  (\t_o, \o_o, \Sigma)(\ve)}(t ,s)$ and the  control curves $\widecheck u^{(\ve, s)}(t)$; 
 \item   by $\cS^{(\ve)}= \widecheck F^{(\t_o, \o_o, \Sigma)(\ve)(2k-1)}([0,T] \times [0,1]) $ the $2$-dimensional submanifold of $ J^{2k-1}(\cQ|\bR) \times \wh K$ spanned by the traces of
 the curves of the homotopy $\widecheck F^{  (\t_o, \o_o, \Sigma)(\ve)(2k-1)}$; 
 \item by $ Y^{(\ve)} =    Y^{(\ve)i}_{(\ell)}  \frac{\p}{\p q^i_{(\ell)}}  + Y^{(\ve)a} \frac{\p}{\p u^a} $  the field of tangent vectors  of $\cS^{(\ve)}$ defined  by 
\beq\label{vectorfield} Y^{(\ve)}|_{\widecheck F^{(\t_o, \o_o, \Sigma)(\ve)(2k-1)}(t,s)}\= \frac{\p \widecheck F^{ (\t_o, \o_o, \Sigma)(\ve)(2k-1)}}{\p s}\bigg|_{(t,s)}\ .\eeq
\end{itemize}
\smallskip
We are now ready to define   the particular   class of needle variations, which are essential for our proof. In the subsequent  Remark \ref{motivgood}, a short explanation of the main ideas which motivates this definition  is given  (see  \cite{CS, CGS} for a discussion in greater detail). 
\begin{definition}  \label{defgood} A  {\it good}  needle variation  of    $\g_o= \g^{( U_o)}$  is a    smoothed needle variation
 $ \cneedle^{(\t_o, \o_o,  \Sigma, \ve_o)} (\g_o)$ which satisfies the following inequality  for any $\ve \in [0, \ve_o]$: 
\begin{multline} \label{additional}     \int_0^T \Bigg(L\big|_{\widecheck F^{  (\t_o, \o_o,  \Sigma)(\ve)(2k-1)}(t,1)}   -   L\big|_{j^{(2k-1)}_t(\g_o)} \Bigg)dt  
 + \int_0^1 \Bigg(- \frac{\p C}{\p q^i_{(\b)}} Y^{(\ve)i}_{(\b)}\bigg|_{{\widecheck F}^{ (\t_o, \o_o, \Sigma)(\ve)(2k-1)}(T,s)} -  \\
-  \sum_{\a = 1}^k \sum_{\b = 0}^{\a-1} (-1)^{\b}  \frac{d^\b}{dt^\b}  \left(\frac{\p L  }{\p q^i_{(\a)}} \right) 
Y^{(\ve)i}_{{(\a-(\b+1))}}  \big|_{\widecheck F^{ (\t_o, \o_o,  \Sigma)(\ve)(2k-1)}(T,s)} \Bigg) ds + \\
 + \int_0^1  \sum_{\a = 1}^k \sum_{\b = 0}^{\a-1} (-1)^{\b}  \frac{d^\b}{dt^\b}  \left(\frac{\p L  }{\p q^i_{(\a)}} \right) 
Y^{(\ve)i}_{{(\a-(\b+1))}}  \big|_{\widecheck F^{ (\t_o, \o_o,  \Sigma)(\ve)(2k-1)}(0,s)} ds \geq 0
\ . \end{multline}
 With a small abuse of language, we sometimes call {\it good} also the (non-smooth) needle variations,  for which the  associated smoothed ones are good.
\end{definition} 
\begin{rem} \label{motivgood} One of the key ideas  of the approach developed in   \cite{CGS}  (see in particular  Sect. 2 of that paper) is the observation   that the terminal cost of a $\cK$-controlled curve $\g = \g^{(U)}(t)$, $U = (u(t), \s)$, is equal to  the  integral along the curve  $t \mapsto (j^{2k-1}_t(\g), u(t))$
in $J^{2k-1}(\cQ|\bR) \times \cK$ of the  $1$-form $\wt L dt + d \wt C$, where $\wt L$ is an appropriate modification of  $L$ and  $\wt C$ is a smooth extension  of the original cost function $C: J^{2k-1}(\cQ|\bR)_{t = T} \to \bR$ over the whole space $J^{2k-1}(\cQ|\bR) \times \cK$, which   vanishes at the submanifold  $J^{2k-1}(\cQ|\bR)_{t = 0} \times \cK$.   The modification $\wt L$ of $L$  is  built in a way that it produces  the same Euler-Lagrange equations of $L$ but has also the additional property of being  identically equal to $0$ along the   solutions of the equations (for the cost problems   of this paper,  the Lagrangian \eqref{lag}  has already this second property and one can just take  $\wt L = L$). 
 The equality between the terminal cost of $\g$ and  the integral of   $\wt L dt + d \wt C $  along the  curve $(j^{2k-1}_t(\g), u(t))$  is an immediate consequence of  the fact that  $\wt L$ vanishes along the solutions. \par
 This  crucial observation implies that  the difference between the terminal costs  of two homotopic controlled curves  is equal to the integral  of  $-(\wt L dt + d \wt C)$ along two arcs out of  the four counter-clockwise oriented  boundary of the surface $\cS \subset J^{2k-1}(\cQ|\bR) \times \cK$ spanned   by the jets   and the controls of the curves of the homotopy. By Stokes' Theorem, the  integral along those two arcs (that is,  the difference between the   terminal costs)   equals the sum of the integral of the $2$-form  $- d( \wt L dt  + d \wt C) = -  d \wt L \wedge dt $ on  $\cS$ and  the integral of   $\wt L dt + d \wt C$ along the  other  two oriented arcs of the boundary. Let us call such two arcs  the ``vertical part of the  boundary" of the homotopy.  \par
   {\it The   condition  \eqref{additional}  is equivalent  to requiring that the   integral of the $1$-form  $\wt L dt + d \wt C$  along the (oriented) vertical part of the  boundary   is non-positive
  for any  homotopy of  a  good needle variation. }
  \par
  Due to this,  the difference between the terminal costs of two controlled curves related   by an homotopy  of a good needle variation is non-negative only if 
  the integral of   $- d\wt L \wedge dt $  on the corresponding surface $\cS$ is non-negative. 
   This  implies that a necessary condition for a controlled curve to be a solution to the cost problem  is that, for any good needle variation, the limit for $\ve \to 0$ of the integral of the $2$-form  $- d\wt L \wedge dt $
   on the surfaces spanned by homotopies is non-negative. In \cite{CS} and \cite{CGS} we made this necessary condition explicit and  show that it reduces to    the classical  PMP in case of  Mayer  problems with smooth data and first order constraints. In the setting of this  paper, the same necessary  condition gives Theorem   \ref{POMPONE} below.   \par
 We finally remark that  \eqref{additional} is basically a condition on the vector fields $Y^{(\ve)}$ we  associated above  with  the homotopies $\widecheck F^{  (\t_o, \o_o, \Sigma)(\ve)(2k-1)}$.  
The existence  of good needle variations   depends on  the existence of a  sufficient  amount of freedom  for  constructing homotopies  with such a property.  In the next Lemma \ref{existenceofgood}) it is shown that 
such a necessary  freedom is granted  for the cost problems of this paper. Roughly speaking this  is essentially due to the fact that   only the derivatives with respect to the $x_{(s)}^i$-variables are 
 relevant for the terminal cost, while the auxiliary variables $p_i$ and their derivatives are  freely specifiable.  We expect that this is a general property, i.e. that  the existence of good needle variations is  always related with   the existence of  a sufficiently large number of 
 (auxiliary and ineffective) variables. 
\end{rem}
\par
\medskip
Now,  for any control pair $U_o = (u_o(t), \s_o)$,  with corresponding curve  $\g_o= \g^{( U_o)}$, and for any  $\t_o \in (0, T)$,  we define  
  \beq\label{pontr-funct} 
 \cP^{(\s_o, u_o, \t_o)}: K \longrightarrow \bR \ ,\qquad  \cP^{(\s_o, u_o, \t_o)} (\o) \= -   L(j^k_{\t_o}(\g_o), \o)\ . \eeq
 By \cite[Cor.7.7]{CGS} the following holds:
 \begin{theo}[Generalised PMP] \label{POMPONE}  
If    $U_o = (u_o(t), \s_o)$ is  an optimal control, then 
  \beq  \label{712}  \cP^{(\s_o, u_o, \t_o)}(u_o(\t_o))  \geq  \cP^{(\s_o, u_o, \t_o)}(\o) \eeq
   for any  pair  $(\t_o, \o) \in (0, T) \times K$,  for which there is at least one {\rm good} needle variation  of the form  $ \cneedle^{(\t_o, \o_o,\Sigma, \ve_o)} (\g_o)$ for $\g_o$.
  \end{theo}
\par
  \medskip
\section{Two  lemmas for the approximation technique}
\label{key}
 According to  Theorem \ref{POMPONE},  {\it if there is a pair $(\t_o, \o_o) \in (0, T) \times K$ for which  there exists an associated {\it good} needle variation for the  $\cK$-controlled curve $\g_o = \g^{(U_o)}$ and   such that the strict  inequality 
   $ \cP^{(\s_o, u_o, \t_o)} (u_o(\t_o)) < \cP^{(\s_o, u_o, \t_o)}(\o_o)  $ occurs,    then $U_o$ cannot be  an optimal control, i.e.  there   must exist  an alternative  $\cK$-controlled curve  $\g \neq \g_o$ with a strictly smaller terminal    cost}.  The proof    in \cite{CGS} of this is constructive and provides   an explicit construction of    curves with smaller costs. \par In  the next two  lemmas, we    present  such a  construction and we  
  make explicit its  dependence   on the data of the Mayer problem and on  the considered  good  needle variation. 
The first lemma holds for any generalised Mayer problem, with no particular assumptions on the  controlled Lagrangian $L$.
It consists of  two  claims:  The first  says that, for all sufficiently small widths  $\ve$, the jets of the controlled curves in 
 the homotopies of a fixed good  needle variations are  in a prescribed  neighbourhood of the jets of the undeformed curve; The second (and more important) claim 
 gives an estimate for  the terminal costs of the deformed curves of the needle variation, which depends  on the value of  the Pontryagin function 
at the parameter $\o_o$ of the needle variation.\par

The second lemma holds  only when the controlled Lagrangian  and the  good needle variations  have   very special forms and shows that 
certain constants, appearing in the statement of the  first lemma,  actually  depend on much  fewer data. 
Both    lemmas apply to   the  cost minimising problems, on which we  focus in this paper, and  play a crucial role in the proof of our main result. \par
 For the first lemma, 
we   need to start  introducing  some useful notation. Given  a  $ \cK$-controlled curve   $ \g_o \=  \g^{( U_o)}$,   $U_o = (u_o(t), \s_o) \in \cK$, and a relatively  compact  neighbourhood   $\cN \subset J^{2k-1}(\cQ|\bR)$ of the  $(2k-1)$-jets of $\g_o$,   for any    $r \in  \bN$ we denote
 \beq 
 \begin{split}
 & {\vertiii L}_{r,\cN} \= \sup_{\smallmatrix (j^{2k-1}_t(\g), u) \in \cN \times K\\ 0 \leq \ell \leq r \ ,\ \ 0 \leq\ell' \leq \ell\endsmallmatrix} \left\{ \left| \frac{\p^\ell L}{( \p t)^{\ell'} \p q^{i_1}_{(m_{1})}\cdots \p q^{i_{\ell-\ell'}}_{(m_{\ell-\ell'})}}\right|\right \}\ ,\\ 
 &{\vertiii{\frac{\p L}{\p u}}}_{r,\cN} \= 
  \sup_{\smallmatrix (j^{2k-1}_t(\g), u) \in \cN \times K\\ 0 \leq \ell \leq r \ ,\ \ 0 \leq\ell' \leq \ell  \endsmallmatrix} \left\{ \left| \frac{\p^{\ell+1} L}{(\p t)^{\ell'} \p q^{i_2}_{(m_2)} \cdots \p q^{i_{\ell-\ell'}}_{(m_{\ell-\ell'})} \p u^a}\right|\right \}\ ,\\
 & \| C\|_{\cC^1, \cN} \= \sup_{j^{2k-1}_T(\g)\in \cN} \left(\big|C(j^{2k-1}_T(\g))\big| +  \sum_{i , s}\bigg| \frac{\p C}{\p q^i_{(s)}}\bigg|_{(j^{2k-1}_T(\g))}\bigg\|\right)\ .
  \end{split}
\eeq
  Further, given a continuous two-parameters family   of  initial conditions $\Sigma = \Sigma(\ve, s) \in\cA_{\text{\it init}} $,   $(\ve, s)  \in [0, \ve_o] \times  [0, 1]$,   we denote by $  \operatorname{diam}(\Sigma)$ the diameter of the set of all such  initial conditions.  We finally recall that  $\gh $ denotes a fixed constant, smaller than $ \frac{1}{2}$,  which 
 appears in the definition of the smoothed needle modifications. 
  \par
    \begin{lemma} \label{corolletto} Let $ \g_o \=  \g^{( U_o)}$  and  $\cN \subset J^{2k-1}(\cQ|\bR)$ be  a  $ \cK$-controlled curve and a relatively compact neighbourhood of  the $(2k-1)$-jets of $\g_o$ as above and assume that there is  at least one good needle variation  $\wt{\cN eedle} =  \cneedle^{(\t_o, \o_o, \wt \Sigma, \wt \ve_o)} (\g_o)$ for  a given  choice of   $(\t_o, \o_o) \in (0, T) \times K$. 
 Let also   $\r$, $\k$, $\gC$, $\gC'$ be the constants that are determined as  in Lemma \ref{normallemma} by $L$, the Lipschitz bijection  between $\cA_{\text{\it init}}$ and $ \wt \cA_{\text{\it init}}$
  and  a cut-off function with support in the open set $\wt \cN \subset \bR^{\wt N+1}$ corresponding  to  $\cN \subset  J^{2k-1}(\cQ|\bR)$ (so that 
    \eqref{const} holds for the controlled curves, which are sufficiently close to  $\g_o$). We finally assume that the map $\ve \mapsto \diam\{\Sigma(\ve, \cdot), s \in [0,1]\}$  is continuous in the argument $\ve$  and that there exists a constant $K_{(\cN, L)}$, depending on $L$ and $\cN$, satisfying the following condition:   the first order system \eqref{first-3} in normal form, which is equivalent to  \eqref{controlledEL},  is such that
    \beq \label{la42} \sup_{(t,y, u) \in \wt \cN\times K} \left\|\frac{\p g^A}{\p u^a}\right\| < K_{(\cN, L)}  {\vertiii{\frac{\p L}{\p u}}}_{k+1,\cN} \  ,\qquad  \sup_{((t, x, u) \in \wt \cN \times K} \left\|\frac{\p g^A}{\p y^\ell}\right\| < K_{(\cN, L)} {\vertiii{L}}_{k+2,\cN}  \ .\eeq
  Then: 
\begin{itemize}[leftmargin = 18pt]
\item[(1)]    There is   a good needle variation $\cneedle^{(\t_o, \o_o, \Sigma, \ve_o)} (\g_o)$  (obtained  from  $\wt{\cN eedle}  =  \cneedle^{(\t_o, \o_o, \wt \Sigma, \wt \ve_o)} (\g_o)$ by appropriately reducing  the  width  $\wt \ve_o$ to a smaller one  $\ve_o$)  such  that: 
    \begin{itemize}[leftmargin = 15pt]
    \item[(a)]  $\diam \Sigma$ and the  distances between  any two control curves $u(t)$, $u'(t)$,   
 corresponding to two curves in a common  homotopy  $F^{(\t_o, \o_o, \Sigma)(\ve)}$, $\ve \in (0, \ve_o)$,  of the family $\cneedle^{(\t_o, \o_o, \Sigma, \ve_o)} (\g_o)$,  are  less than or equal to $\r$;  
   \item[(b)]    any    $(2k-1)$-th order jet of a control curve in  the homotopies of    $\cneedle^{(\t_o, \o_o, \Sigma, \ve_o)} (\g_o)$  is   in  $\cN$.
   \end{itemize}
\item[(2)]  If     $\k_o =  \cP^{(\s_o, u_o, \t_o)}(\o_o)  - \cP^{(\s_o, u_o, \t_o)}(u_o(\t_o))> 0$, then 
 there exists  a constant  $M$,  depending on  $\cN$,   $\t_o$, $\o_o$, $\Sigma$, ${\vertiii L}_{k+2,\cN}$, ${\vertiii{\frac{\p L}{\p u}}}_{k+1, \cN} $, $\| C\|_{\cC^1, \cN}$
 and the  infinitesimal (\footnote{We use  the short expression  ``infinitesimal''  to mean that   $\lim_{\ve \to 0} V(\ve) = 0$.})
\beq\label{infinitesimal1} V(\ve) \= \frac{1}{\ve}  \int_{\t_o - \ve}^{\t_o } \left( \cP^{(\s_{o},  u_{o}, t)}(u_o(\t_o)) - \cP^{(\s_{o},  u_{o}, t)}(u_o(t))  \right)dt    \eeq
 such that  the  controlled curve $\underline \g = \g^{(\underline U)}$,  determined by the  pair  $\underline U = (\underline u(t), \underline \s) \in \cK_{\text{meas}}$ with   $\underline \s \= \Sigma\left( \frac{1}{2M}, 1\right)$ 
and 
\beq
 \underline u(t) \=  \left\{\begin{array}{ll} u_o(t)&\text{if}\ \  t \in [0, T] \setminus  (\t_o -  \frac{1}{2M}, \t_o]\ ,\\[10pt]
 \o_o &\text{if}\ \   t \in  (\t_o  -  \frac{1}{2M}, \t_o]\ ,
\end{array}
\right.
\eeq
has a  terminal cost  that satisfies 
\beq  \label{325}
   C ^{(\underline \g)}  \leq  C^{(\g_o)} -  \frac{\k_o}{4 M} <  C^{(\g_o)} \ .
\eeq
\end{itemize}
There    also  exists  a  constant  $\widecheck M > 0$, for which  \eqref{325} holds  for  the  {\rm smooth} $\cK$-controlled curve  $\underline {\widecheck \g}$, determined  by the same initial condition and  the   {\rm smoothed   version}   $\widecheck{\underline u}(t)$ of  $\underline u(t)$. 
\end{lemma}
\begin{pf} (1)   Given a   good needle variation 
$\wt{\cN eedle} = \cneedle^{(\t_o, \o_o, \Sigma,  \ve_o)} (\g_o) = \{\widecheck F^{(\t_o, \o_o,\Sigma)(\ve)}, \ve \in [0, \ve_o]\}$,  for any $(\ve,s) \in [0, \ve_o] \times [0,1]$ we denote by  $\g^{(\ve, s)}(t) $ and $\widecheck \g^{(\ve, s)}(t)$  the controlled curves 
$$\g^{(\ve, s)}(t) \=  F^{(\t_o, \o_o, \Sigma)(\ve)}(t,s) \ ,\qquad \widecheck \g^{(\ve, s)}(t) \= \widecheck F^{(\t_o, \o_o,  \Sigma)(\ve)}(t,s)\ .$$
We  also use the short notation
$ \g^{(\ve)}(t) \=  \g^{(\ve, 1)}(t) $ and $\widecheck \g^{(\ve)}(t) \=  \widecheck \g^{(\ve, 1)}(t) $.
By construction 
\begin{itemize}[leftmargin = 15pt]
\item $\g^{(\ve,s)}(t)$ is  the $\wh \cK_{\text{meas}}$-controlled curve determined by  $U^{(\ve, s)} \= (u^{( \ve, s)}(t), \Sigma(\ve, s))$ with 
\beq\label{mod1} u^{(\ve, s)}(t) = (1-s) u_o(t) + s  u^{(\t_o, \o_o, \ve)} (t) \ ,\qquad 
 u^{(\t_o, \o_o, \ve)} (t) \=  \left\{\begin{array}{ll} u_o(t)&\text{if}\ \  t \notin  (\t_o - \ve, \t_o]\ ,\\[10pt]
 \o_o &\text{if}\ \   t \in  (\t_o  -  \ve, \t_o]\ .
\end{array}
\right.
\eeq
\item   $\widecheck \g^{(\ve,s)}(t)$ is  the $\wh \cK$-controlled curve determined by  the same initial condition and the smoothed version $\widecheck u^{(\t_o, \o_o, \ve)} (t) $ of $u^{(\t_o, \o_o, \ve)} (t) $, defined in  \eqref{mod1}.
\end{itemize}
 We now replace the good  needle variation $\wt {\cN eedle}$, given in the statement,  by a new one, in which  the homotopies of  curves are the same, but where the maximum value for  the  parameter $\ve$  is  changed into a new value  $\ve_o \leq \wt \ve_o$.   This  $\ve_o$ is chosen small enough  to make the family of initial conditions $\Sigma = \wt \Sigma|_{[0, \ve_o] \times [0,1]}$ such that   $\diam(\Sigma) <  \r$ and,  for any $(\ve, s)  \in  [0, \ve_o] \times  [0,1]$, 
$$\operatorname{dist}(\widecheck u^{(\ve, s)}, u^{(\ve,s)}) \leq 2 \gh \ve^2 < \r\ ,\qquad \operatorname{dist}(\widecheck u^{(\ve, s)},u_o) \leq \ve + 2 \gh \ve^2 < \r\ , $$
In this way   (a) is satisfied. From  Lemma \ref{normallemma}, by the fact that  $\Sigma(0,s) = \s_o$ for any $s$  and from the Lipschitzian assumption on $\Sigma$, it follows that for a sufficiently small $\ve_o$
\begin{align}
&  \label{neigh}  \|  \widecheck  \g^{(\ve,s)} - \g_o\|_{\cC^{2k-1}} \leq  \gc( \ve + 2 \gh \ve^2 ) + \k |\Sigma(\ve,s) - \Sigma(0,s)|< (2 \gc + \k C) \ve  \ ,
\\
& \label{neigh''}  \| \widecheck \g^{(\ve,s)} - \g^{(\ve,s)}\|_{\cC^{2k-1}} \leq 2 \gh  \gc \ve^2  < \gc \ve^2 \ , 
\end{align}
for a constant $C$ determined by $\Sigma$.  By possibly taking a  smaller $\ve_o$,  also (b) is   satisfied.   \par
(2)  Consider the good needle variation determined  in (1) and its  smoothed version. For any $\ve \in [0, \ve_o]$, we   denote
$\widecheck C^{(\ve)}\= C(j^{2k-1}_{t = T}(\widecheck  \g^{(\ve)}))$ and  $C^{(\ve)}\= C(j^{2k-1}_{t = T}( \g^{(\ve)}))$.
By  \cite[Cor. 6.5]{CGS} we know that 
\begin{multline}  
\label{217}    C^{(0)} -  \widecheck C^{(\ve)}  =  \widecheck C^{(0)} -  \widecheck C^{(\ve)}  = \\
=    \int_0^T \left( \int_0^1  Y^{(\ve)a} \frac{\p \cP^{(\Sigma(\ve,s), \widecheck u^{(\ve,s)}, t)} }{\p u^a} \Bigg|_{\widecheck u^{(\ve,s)}(t)} ds  \right) dt 
- \int_0^T \left( \int_0^1  \frac{\p^2\wh  \mu}{\p t \,\p s}\bigg|_{(t,s)}ds \right)dt
 \end{multline}
where $Y^{(\ve)a}$ are the $\frac{\p}{\p u^a}$-components  of the vector field   \eqref{vectorfield} and 
$\wh \mu:[0, T] \times [0,1] \to \bR$  is  an appropriate    function,  which is  completely determined by the    needle variation considered. For our purposes, there is no need to recall   the   detailed definition  of $\wh \mu$, but only to know  that, since  the  needle variation is good,   {\it  the value of  $\displaystyle \int_0^T \left( \int_0^1  \frac{\p^2\wh  \mu}{\p t \,\p s}\bigg|_{(t,s)}ds \right)dt $ is non-positive} (\cite[Lemma 7.6]{CGS}). In addition, 
the next sublemma, whose proof is quite  technical and  is  postponed to    \S\ \ref{sublemmaproof}, gives an estimate for the first term in   \eqref{217}. 
\begin{sublemma}  \label{lemma37} In the  hypotheses of the lemma,  
there exists a  constant $\gn = \gn_{(\t_o, \cN, L, \frac{\p L}{\p u})}$,  depending  on $\t_o$, $\cN$,   ${\vertiii{L}}_{k+2,\cN}$ and  ${\vertiii{\frac{\p L}{\p u}}}_{k+1,\cN}$, such that  for any $\ve \in [0, \ve_o]$
\begin{multline} \label{329}  \int_0^T \left( \int_0^1  Y^{(\ve)a}\frac{\p \cP^{(\Sigma(\ve,s), \widecheck u^{(\ve,s)}, t)} }{\p u^a} \Bigg|_{\widecheck u^{(\ve,s)}(t)} ds  \right) dt  >\\
>   \int_{\t_o - \ve - \gh \ve^2}^{\t_o + \gh \ve^2} \left(\ \cP^{(\s_{(\ve)}, \widecheck u_{(\ve)}, t)}(\widecheck u^{(\t_o, \o_o, \ve)}(t)) -  \cP^{(\s_o,  u_o, t)}(u_o(t)) \right) dt - \gn \ve^2
\end{multline}
where  $\s_{(\ve)}\= \Sigma(\ve, 1)$ and $\widecheck u_{(\ve)} = \widecheck u^{(\ve, 1)}$
\end{sublemma}
Let us now focus  on the right hand side of \eqref{329}. The first integral  decomposes into
\begin{multline}  \label{pp}  \int_{\t_o - \ve - \gh \ve^2}^{\t_o + \gh \ve^2} \left( \cP^{(\s_{(\ve)}, \widecheck u_{(\ve)}, t)}(\widecheck u^{(\t_o, \o_o, \ve)}(t)) - \cP^{(\s_o, u_o, t)}(u_o(t)) \right) dt= \\
=
 \int_{\t_o - \ve - \gh \ve^2}^{\t_o - \ve} \left( \cP^{(\s_{(\ve)}, \widecheck u_{(\ve)}, t)}(\widecheck u^{(\t_o, \o_o, \ve)}(t)) -  \cP^{(\s_o,  u_o, t)}(u_o(t)) \right) dt  + \\
 +  \int_{\t_o - \ve}^{\t_o } \left(\cP^{(\s_{(\ve)}, \widecheck u_{(\ve)}, t)}(\o_o) -  \cP^{(\s_o,  u_o, t)}(u_o(t)) \right) dt  + \\
 +   \int_{\t_o }^{\t_o + \gh \ve^2} \left(\cP^{(\s_{(\ve)}, \widecheck u_{(\ve)}, t)}(\widecheck  u^{(\t_o, \o_o, \ve)}(t)) -  \cP^{(\s_o,  u_o, t)}(u_o(t)) \right) dt \end{multline}
Since $L$ is continuous, for any $(t, \o) \in [0,T]\times \wh K$, we have that 
$\left|\cP^{(\Sigma(\ve,s), \widecheck u^{(\ve,s)}, t)}(\o)\right| < \gc''_{(\cN, L)}$,  with   $\gc''_{(\cN, L)} \= \| L\|_{\infty, \cN \times K}$, 
and   the sum of the  first and the third terms  in \eqref{pp} is    bounded by $4 \gc''_{(\cN, L)} \gh \ve^2$.  On the other hand, the second   decomposes into 
\begin{multline} \label{41212}
 \int_{\t_o - \ve}^{\t_o } \left(\cP^{(\s_{(\ve)}, \widecheck u_{(\ve)}, t)}(\o_o) -  \cP^{(\s_o,  u_o, t)}(u_o(t)) \right) dt  = \\
 =  \int_{\t_o - \ve}^{\t_o } \left(\cP^{(\s_{(\ve)}, \widecheck u_{(\ve)}, t)}(\o_o) -  \cP^{(\s_o,  u_o, \t_o)}(\o_o) \right) dt  +  \\
 +\underset{= \ve \k_o}{ \underbrace{\int_{\t_o - \ve}^{\t_o } \left(\cP^{(\s_o,  u_o, \t_o)}(\o_o) -  \cP^{(\s_o,  u_o, \t_o)}(u_o(\t_o)) \right) dt}}  + \\
 + \int_{\t_o - \ve}^{\t_o } \left(\cP^{(\s_o,  u_o, \t_o)}(u_o(\t_o)) -  \cP^{(\s_o,  u_o, t)}(u_o(\t_o)) \right) dt + \\
 + \int_{\t_o - \ve}^{\t_o } \left(\cP^{(\s_o,  u_o, t)}(u_o(\t_o)) -  \cP^{(\s_o,  u_o, t)}(u_o(t)) \right) dt 
\end{multline} 
By \eqref{neigh}  the   integrands of the first and the third terms are bounded  by 
\beq\label{2ter} \left| \cP^{(\s_{(\ve)}, \widecheck u_{(\ve)}, t)}(\o_o) - \cP^{(\s_{o},  u_{o}, \t_o)}(\o_o) \right| ,\ \  \left| \cP^{(\s_{o},u_o, \t_o)}(u_o(\t_o)) - \cP^{(\s_{o},  u_{o},t)}(u_o(\t_o)) \right|  <  \gc'''_{(\cN, L, \Sigma)} \ve
\eeq
for some   $\gc'''_{(\cN, L, \Sigma)}$  depending  on ${\vertiii L}_{1,\cN}$ and the constants $\gc$, $\k$, $C$.  On the other hand, by  continuity,   the fourth term is an  infinitesimal  of higher order than $\ve$, that is (see \eqref{infinitesimal1})
\beq \label{2quater} \int_{\t_o - \ve}^{\t_o }   \left(\cP^{(\s_{o},  u_{o}, t)}(u_o(\t_o)) - \cP^{(\s_{o},  u_{o}, t)}(u_o(t)) \right)dt   =  V(\ve)  \ve \ .
\eeq
Combining \eqref{217} with Sublemma \ref{lemma37} and the above discussion,   we obtain that 
\beq  \label{358}   C^{(0)}  -   \widecheck C^{(\ve)}\geq \ve(\k_o -  4 \gc''_{(\cN,L)} \gh \ve - 2 \gc'''_{(\cN, L, \Sigma)}\ve -  |V(\ve)| - \gn \ve) \ .\eeq
On the other hand, by   \eqref{neigh} and the smoothness of  the cost function,  there is  a  constant   $\gm = \gm_{(\Sigma, C)}$,  depending on    $ \| C\|_{\cC^1, \cN}$,   such that  
for all sufficiently small $\ve$
\beq \label{358-bis}|\ \widecheck C^{(\ve)} -  C^{(\ve)} | < \gm \| \widecheck \g^{(\ve)} - \g^{(\ve)}\|_{\cC^{2k-1}} \leq   \gc \gm\ve^2\ .\eeq
Hence, from \eqref{358} and \eqref{358-bis} we get 
\beq  \label{358-ter}  \begin{split} & C^{(\ve)}\leq  C^{(0)}  -    \ve(\k_o -   \gd \ve - |V(\ve)|) ) \ ,\\
& \text{where}\ \ \gd  = \gd_{(\t_o, u_o, \cN, L,  \Sigma, C)}\= (4 \gc''_{(\cN,L)} \gh  + 2 \gc'''_{(\cN, L, \Sigma)} +  \gc \gm + \gn)\ .
\end{split}\eeq
Since the map $\ve \mapsto   \gd \ve + |V(\ve)| $ is  infinitesimal for $\ve \to 0^+$, there exists  $M> 0$ such that 
$$ \k_o -  ( \gd \ve + |V(\ve)|)  > \frac{\k_o}{2}\qquad \text{for all}\ \ve \in \left[0, \frac{1}{M}\right] \ .$$ 
Thus,  setting  $\underline \ve \= \frac{1}{2M}$, we have  
$ C^{(\underline \ve)}  \leq   C^{(0)}  -   \underline \ve \frac{\k_o}{2}  =   C^{(0)}  -     \frac{\k_o}{4 M}$, 
that is  \eqref{325}. 
The last   claim is  proved similarly, using just  \eqref{358} in place of   \eqref{358-ter}.\end{pf}
\par
\smallskip
 We now  present the second advertised  lemma, which gives   a radical improvement of  Lemma \ref{corolletto}  under  additional   assumptions on  $L$ and on the  considered good needle variation. More precisely, we assume that the configuration space has the form  $\cQ = \sQ \times \sQ^*$  for an $n$-dimensional affine space $\sQ = \bR^n$ (the coordinates are thus   pairs $q = (x, p)$ with $x  = (x^i)\in \sQ$ and $p = (p_j) \in \sQ^*$) and  $L$ has  the  form 
\beq \label{specialform}  L =  p_i \big(x^i_{(k)} - f^i\big(t, x^i, x^i_{(1)}, \ldots, x^i_{(k-1)}, u^a\big)\big)\ .\eeq
\begin{lem} \label{corolletto-bis} Assume the  hypotheses of Lemma \ref{corolletto} with the exception of the condition \eqref{la42} and let  $\cQ = \sQ \times \sQ^*$  and $L$ are as above. Moreover, for any $\ve \in [0, \ve_o]$,  let   $\cG(t,s) \= \left(j^{2k-1}_{t}({\widecheck \g}{}^{(\ve, s)}), u^{(\ve,s)}(t)\right)$ and   $Y $  the  field 
 of tangent vectors  to  the surface $\cS \= \cG([0, T] \times [0,1])$ defined by  
$Y|_{\cG(t,s)} = \frac{\p \cG}{\p s}\big|_{(t,s)} $.
\par If the family of  initial conditions $\Sigma$, occurring in the definition of   $\needle^{(\t_o, \o_o, \Sigma, \ve_o)} (\g_o)$,   
is such that 
\begin{multline} \label{3.54}
   \sum_{\d = 1}^{k} \sum_{\h = 0}^{\d-1} (-1)^{\h}  \int_0^1  \frac{d^\h}{dt^\h}  \left(\frac{\p  L  }{\p q^i_{(\d)}} \right) Y^i_{(\d-(\h+1))}\bigg|_{\cG(0 ,s)}
ds = \\
= \sum_{\d = 1}^{k} \sum_{\h = 0}^{\d-1} (-1)^{\h}  \int_0^1 \frac{d^\h}{dt^\h}  \left(\frac{\p  L  }{\p q^i_{(\d)}} \right) Y^i_{(\d-(\h+1))}\bigg|_{\cG(T ,s)}
ds =  0\ ,
\end{multline}  
then  claim (2)  of Lemma \ref{corolletto} holds with    constants 
 $M , \widecheck M$  that depend  just on $\t_o$, $\o_o$, $\cN$, ${\vertiii L}_{1,\cN}$ and  $\| C\|_{\cC^1, \cN}$ and  not on  ${\vertiii{L}}_{k+2, \cN} $ and ${\vertiii{\frac{\p L}{\p u}}}_{k+1, \cN} $. 
\end{lem}
Due to the technicalities in the arguments,  the proof  is  given  later   in \S\ \ref{proofsecondlemma}.\par
\smallskip
  \begin{rem} \label{remarketto} 
 Let  $(\cK, L^{(\d)}, C^{(\d)})$, $\d \in (0, \d_o]$  be  a  one-parameter family of defining triples 
 as  in Remark \ref{remark33}. Assume that all associated generalised    Mayer problems are with smooth data and of normal type and   that, for  $\d$ tending to $0$, the  Lagrangians $L^{(\d)}$  and the cost functions $C^{(\d)}$ tend uniformly on compacta to continuous   functions $L^{(0)}$ and $C^{(0)}$.  Assume also that: 
  \begin{itemize}[leftmargin = 20pt]
 \item[(A)] The  partial  derivatives   
 $ $
 \begin{align*} &\frac{\p^\ell L^{(\d)}}{\p q^{i_1}_{(\a_1)} \ldots \p q^{i_\ell}_{(\a_\ell)}}\ ,\   \frac{\p^\ell L^{(\d)}}{(\p t)^{\ell'}\p q^{i_1}_{(\a_1)} \ldots \p q^{i_{\ell - \ell'}}_{(\a_{\ell - \ell'})}}\qquad \qquad\text{with}\ \  1 \leq  \ell \leq k+2\\\
 \text{and} \qquad & \frac{\p^{\ell+1} L^{(\d)}}{\p u^a \p q^{i_1}_{(a_1)}\ldots \p q^{i_\ell}_{(\a_\ell)}}\ ,\  \frac{\p^{\ell+1} L^{(\d)}}{\p u^a  (\p t)^{\ell'} \p q^{i_1}_{(\a_1)}\ldots \p q^{j_{\ell-\ell'}}_{(\a_{\ell-\ell'})}}\qquad \text{with}\  1 \leq  \ell \leq k+1
 \end{align*}
  tend uniformly on compacta  to the corresponding partial derivatives of $L^{(0)}$; 
  \item[(B)] There is a one-parameter family $U^{(\d)} = (u_o^{(\d)}(t), \s_o^{(\d)}) \in \cK$,   
  whose associated controlled curves $\g^{(\d)}(t) \= \g^{(U^{(\d)})}$  converge in the norm of $\cC^{k-1}([0, T])$  to a curve $\g^{(0)}(t)$  such that: (a) it is a solution to the  differential constraints determined by  $L^{(0)}$, (b) it has $\s_o^{(0)} = \lim_{\d \to 0} \s_o^{(\d)}$ as initial condition;  (c) it  is determined by a measurable control  curve $u^{(0)}_o(t) \in K$ with $u_o^{(0)}(t) = \lim_{\d \to 0} u_o^{(\d)}(t)\ \text{a.\ e.}$; 
\item[(C)] There exists a $\d$-parameterised  family 
 of    initial data maps $\Sigma^{(\d)}(\ve,s)$,  converging  uniformly  on $[0, \ve_o] \times [0,1]$ to  a  limit map $\Sigma^{(0)}(\ve,s)$,  and 
 a corresponding   $\d$-parameterised   family of good needle variations,  determined by   the  maps $\Sigma^{(\d)}$ and a pair $(\t_o, \o_o)$, in which  $\t_o$ is  one of the  points where   $u_o^{(0)}(\t_o) = \lim_{\d \to 0} u_o^{(\d)}(\t_o)$;  
\item[(D)] The real value
  $$\k_o^{(0)} \= \left(\cP^{(\s^{(0)}_{o},  u^{(0)}_{o}, \t_o)}(\o_o) - \cP^{(\s^{(0)}_{o},  u_{o}, \t_o)}(u^{(0)}_o(\t_o)) \right) = \lim_{\d \to 0} \k_o^{(\d)}$$
  is strictly positive; 
    \item[(E)] The functions  $V^{(\d)}(\ve)$, $\ve \in (0, \ve_o]$, which are  defined by  \eqref{infinitesimal1} for each $\d$, tend uniformly on compacta   to  a 
   function $\wt V(\ve): (0, \ve_o] \to \bR$,  which 
 is an infinitesimal for $\ve \to 0$.
  \end{itemize}
   Note that  the  $\d$-parameterised family of control pairs $U^{(\d)}(\ve, s) \= \big(u^{(\d)(\ve, s)}(t),   \Sigma^{(\d)}(\ve, s)\big)$,   which determine the  curves  $\g^{(\d)(\ve, s)}(t)$ of the  needle variation  in  (C), have the following property: {\it for any   $(\ve, s) \in [0, \ve_o] \times [0,1]$, the a.e. limit curve 
   $$ u^{(0)(\ve, s)}(t)  \= \lim_{\d \to 0} u^{(\d)(\ve, s)}(t) $$
      is a needle modification of $u_o(t)$}. We may  therefore consider also the following additional assumptions:
   \begin{itemize}[leftmargin = 20pt]
   \item[(F)] For any $\ve \in [0, \ve_o]$ the controlled curves $\g^{(\d; \ve)}(t) \=  \g^{(\d) (\ve, 1)}(t)$ converge uniformly to the solution of the differential constraints for $\d = 0$, which is determined by the control pair   $ \big(u^{(0)(\ve, 1)}(t),   \Sigma^{(0)}(\ve, 1)\big) $: 
   \item[(G)] There is a relatively compact neighbourhood of the set 
   $$\cE\text{\it nds} \=  \{ j^{2k-1}_T(\g^{(\d) (\ve, s)}), \ (\d, \ve, s) \in [0, \d_o] \times [0, \ve_o] \times [0,1]\  \}\ ,$$  
   given by  the $(2k-1)$-jets  at   $t = T$  of the controlled curves of the needle variations, on which $C^{(\d)}$ tend to $C^{(0)}$ in the $\cC^1$-norm.
   \end{itemize}
   \par
 \smallskip
For each  of the above  Lagrangian $L^{(\d)}$, we  may follow the proof of Lemma \ref{corolletto-bis} and derive the  inequality  \eqref{358-ter} for  any sufficiently small $\ve$,  i.e. the  inequality
  \beq\label{manyineq}  C^{(\d)(\ve)}\leq  C^{( \d)(0))}  -    \ve(\k^{(\d)}_o -   \gd^{(\d)} \ve - |V^{(\d)}(\ve)| )\eeq
    relating the cost $C^{(\d)(0)}$ of  the  controlled curve $\g^{(\d)}(t)$  with   the cost  $C^{(\d)(\ve)}$  of the controlled curves $\g^{(\d; \ve)}(t)$. Notice that, under the  assumptions (A) --  (G),  for $\d$ sufficiently small,  we may assume that the constants  $\gd^{(\d)}$   are independent of $\d$, say    $\gd^{(\d)} = \gd$, so that,  letting $\d \to 0$,   
    \beq  C^{(\d = 0)(\ve)}\leq  C^{(\d = 0)(0)} -    \ve(\k^{(0)}_o -   \gd \ve - |\wt V(\ve)|) )\ , \eeq
    where   $C^{(\d = 0)(\ve)}$  is the terminal cost of the controlled curve $\g^{(\d = 0, \ve)}(t) $.  
  From this, using  the same concluding argument of  the proof of Lemma \ref{corolletto}, we obtain  the existence of a  constant $M > 0$ and   an associated  needle modification     for   the  (merely measurable) limit  curve  $u_o^{(0)}(t)$, such that the corresponding cost   satisfies
      \beq  C^{(\d = 0)(\underline \ve)}\leq  C^{(\d = 0)(0)}  -    \frac{\k^{(0)}_o}{4 M}\ . \eeq
      This  fact  is crucially exploited in  
    the proof of our main result,  given in the last section. \par
We finally observe that, by Lemma \ref{corolletto-bis}, if the  Lagrangians $L^{(\d)}$ have the special form  \eqref{specialform} and all  maps $\Sigma^{(\d)}$  satisfy the condition 
 \eqref{3.54},  the   above conclusion on the costs of  the limit controlled curve $\g^{(0)}_o$  and of its needle modifications holds  also  if  (A) is replaced by the  following weaker assumption:  
   \begin{itemize}[leftmargin = 25pt]
   \item[(A')]  the  partial   derivatives
   $\frac{\p L^{(\d)}}{\p q^{i}_{(\a)}}$, $\frac{\p L^{(\d)}}{\p t}$
   tend uniformly on compacta  to the corresponding partial derivatives of $L^{(0)}$.
  \end{itemize}  
  \end{rem} 
\par 
\medskip

\section{The  proof of    Theorem \ref{main1}}
\label{concluding}
\subsection{A  preliminary ``smooth''  version of the main result} \label{sect51}  As we pointed in \S \ref{sect51*},  when all of its  data are of class $\cC^\infty$,  the cost problem presented in  the Introduction    is equivalent to the  generalised Mayer problem determined by a the  defining triple $(\cK, L, \sC)$ given in that section. We recall that the
 configuration space has the form  $\cQ = \sQ \times \sQ^*$,   $\sQ = \bR^n$,  the controlled Lagrangian $L$ and the  cost function $\sC$ are given in \eqref{lag},  and the set  $\cK$ consists of the pairs 
 $U = (u(t), \s = \big(\mathsf s, \wt{\mathsf s})\big) $  where   (a)     $\mathsf s$  is an initial condition  in a prescribed    set $\cA_{\text{\it init}} \subset J^{k-1}(\sQ|\bR)|_{t = 0}$  for the curve $x(t)$ and  (b)    $\wt{\mathsf s}$ is an initial condition, which  can be  arbitrary,  for the curve $p(t)$. We finally recall that the controlled Euler-Lagrange equations are of normal type and are given in \eqref{controlledEL-1*} and \eqref{controlledEL-1bis*}.
\par
 \medskip
Let us now consider the following subset of $\cK$.   Given  $\sfs \in \cA_{\text{\it init}}$ and a  control curve $u(t) \in K$,  we denote by    $x^{(u, \sfs)}(t)$  the unique   solution to  \eqref{controlledEL-1*} with  
  initial condition $\sfs$. We then denote by  $p^{(u, \sfs)}(t)$   the unique solution to  \eqref{controlledEL-1bis*}  with  $x(t) = x^{(u, \sfs)}(t)$,  that  satisfies  the terminal conditions \eqref{31-ter} -- \eqref{33-ter}.  Finally, we set  $\wt{\mathsf s}^{(u, \mathsf s)}$ to  be the initial jet  $\wt{\mathsf s}^{(u, \mathsf s)}= j^{k-1}_{t = 0}(p^{(u, {\mathsf s})})$ of  $p^{(u, {\mathsf s})}(t)$.
 By construction, the pair
 \beq\label{goodpair} U^{(u, \mathsf s)} \= \left(u(t), \s = (\mathsf s, \wt{\mathsf s}^{(u, \mathsf s)})\right) \ ,\eeq
 has $\g^{(U)}(t) \= \left(t, x^{(u, \sfs)}(t), p^{(u, {\sfs})}(t)\right)$ as  associated  $\cK$-controlled curve.  The pairs \eqref{goodpair} and the corresponding controlled curves are called {\it good}. The subset of the good pairs in $\cK$, $\wh \cK$, $\wh \cK_{\text{meas}}$ are  denoted by  $\cK_{\text{good}}$, $\wh \cK_{\text{good}}$, $\wh \cK_{\text{good meas}}$, respectively.\par
\medskip

Our interest in the good $\cK$-controlled curves   comes from the following  lemma. 
\begin{lem} \label{existenceofgood} Let  $\g_o= \g^{( U_o)}$ be a    $ \cK_{\text{\rm good}}$-controlled curve and 
$ {\needle}^{(\t_o, \o_o,  \Sigma, \ve_o)} (\g_o)$ a needle variation,  
whose associated smoothed needle variation $ {\cneedle}^{(\t_o, \o_o, \Sigma,  \ve_o)} (\g_o)\= \{\widecheck  F^{(\t_o, \o_o, \Sigma)( \ve)}, \ve \in [0, \ve_o]\}$ satisfies the following two conditions:
\begin{itemize}[leftmargin = 25pt]
\item[(a)] all  control pairs $U(s, \ve)$ that determine  the $\wh \cK$-controlled curves $\g^{(s, \ve)} =   \widecheck F^{(\t_o, \o_o,  \Sigma)(\ve)}(\cdot, s)$ are good; 
\item[(b)]  the  initial conditions  for the $x$-components $x^{(s, \ve)}(t)$ of the curves  $\g^{(s, \ve)}$  are constant and    independent of $(s, \ve)  \in [0,1] \times  [0, \ve_o]$. 
\end{itemize}
Then  ${\cneedle}^{(\t_o, \o_o, \Sigma, \ve_o)} (\g_o)$ is a {\rm good} needle variation in the sense of Definition \ref{defgood}. 
\end{lem}
\begin{pf} First of all, we claim that  if a controlled curve $\g(t) \= \left(t, x^{(u, \sfs)}, p^{(u, \sfs)}(t)\right)$ is good, then 
\begin{multline} \label{condi} \left.\left( \frac{\p {\mathsf C}}{\p x^i_{(\b)}} + \sum_{\ell = 0}^{k - \b - 1} (-1)^\ell \frac{d^\ell}{dt ^\ell} \left( \frac{\p (p_m( x^m_{(k)} - f^m))}{\p x^i_{(\b + \ell +1)}}\right)\right)\right|_{
j^{2k-2}_{t = T}(x^{(u, \sfs)})} = 0 \\
 \text{for}\ \ \ 0 \leq \b \leq k-1\ \ \text{and} \ \ 0 \leq i \leq n\ .
\end{multline}
This is a consequence of the fact that,  setting $\r\= k - \b -1$, the conditions   \eqref{condi} become 
 \beq \label{condi1} 
\begin{split} 
&\left.\left( \frac{\p  {\mathsf C}}{\p x^i_{(k - \r -1)}} + \sum_{\ell = 0}^{ \r} (-1)^{\ell}  \frac{d^{\ell}}{dt ^{\ell} } \left( \frac{\p (p_m( x^m_{(k)} - f^m))}{\p x^i_{(k - \r  + \ell)}}\right)\right)\right|_{
j^{2k-2}_{t = T}(x^{(u, \sfs)})} = \\
& = 
\left.\left( \frac{\p  {\mathsf C}}{\p x^i_{(k - \r -1)}} +
 \sum_{\ell = 0}^{ \r  -1 } (-1)^{\ell}  \frac{d^{\ell}}{dt ^{\ell} } \left( p_m\frac{\p f^m}{\p x^i_{(k - \r + \ell )}}\right) -  (-1)^{\r} p_{i(\r)}   \right)\right|_{
j^{2k-2}_{t = T}(x^{(u, \sfs)})} = 0\ , 
\end{split}
\eeq
 and these are precisely   the conditions   \eqref{31-ter} -- \eqref{33-ter}.
  Consider now an arbitrary needle variation  ${\needle}^{(\t_o, \o_o,  \Sigma, \ve_o)} (\g_o)$ of $\g_o \= \g^{(U_o)}$ with $U_o = (u_o(t), \sfs_o)$. By the particular form of the  differential constraints, the  value of the controlled Lagrangian $L$ is $0$ at all  $(k-1)$-th order jets  of the  $\wh \cK$-controlled curves.  Hence, using the short-hand notation $\g^{(s,\ve)(k-1)}(t) \= j^{k-1}_t(\g^{(s,\ve)})$,  we may write 
\beq   \int_0^T \Bigg(L\big|_{\g^{(s,\ve)(k-1)}(t)}   -   L\big|_{\g^{(k-1)}_o(t)} \Bigg)dt   = 0 \qquad \text{for any}\ s \in [0, 1]\ .\eeq
From this and the fact that  the only  non-trivial  derivatives
$\frac{\p L  }{\p q^i_{(\a)}}$ with respect to  the jet coordinates $q_{(\a)} = (x^i_{(\a)}, p_{(\a)j})$,  $\a \geq 1$,  are  those   with   $q^i_{(\a)}  = x^i_{(\a)}$, it follows that \eqref{additional}  is satisfied if and only if 
\begin{multline*}       \int_0^1 \Bigg(-  \sum_{\b = 1}^{k-1}\frac{\p \mathsf C}{\p x^i_{(\b)}} Y^{(\ve)x^i}_{(\b)}\bigg|_{\g^{(s, \ve)(k-1)}(T)} -  \\
 - \sum_{\d = 1}^k \sum_{\h = 0}^{\d-1} (-1)^{\h}  \frac{d^\h}{dt^\h}  \left(\frac{\p L  }{\p x^i_{(\d)}} \right) 
Y^{(\ve)x^i}_{{(\d-(\h+1))}}\big|_{\g^{(s,\ve)(2k-1)}(T)} \Bigg) ds + \\
 + \int_0^1  \sum_{\d = 1}^k \sum_{\h = 0}^{\d-1} (-1)^{\h}  \frac{d^\h}{dt^\h}  \left(\frac{\p L  }{\p x^i_{(\d)}} \right) 
Y^{(\ve)x^i}_{{(\d-(\h+1))}}  \big|_{\g^{(s,\ve)(2k-1)}(0)} ds = 
\end{multline*}
\begin{multline}\label{49}
=-  \int_0^1\left. \sum_{\b = 0}^k\Bigg(\Bigg(  \frac{\p  {\mathsf C}}{\p x^i_{(\b)}} +\sum_{\ell = 0}^{k - \b - 1} (-1)^\ell \frac{d^\ell}{dt ^\ell} \left( \frac{\p (p_m( x^m_{(k)} - f^m))}{\p x^i_{(\b + \ell +1)}}\right)\Bigg)\right|_{
j^{2k-2}_{t = T}(x^{(u_o, \sfs_o)})} {\cdot}\\
{\cdot}Y^{(\ve)x^i}_{(\b)}\bigg|_{\g^{(s,\ve)(k-1)}(T)} \Bigg)ds + \\
+ \int_0^1  \sum_{\d = 1}^k \sum_{\h = 0}^{\d-1} (-1)^{\h}  \frac{d^\h}{dt^\h}  \left(\frac{\p L  }{\p x^i_{(\d)}} \right) 
Y^{(\ve)x^i}_{{(\d-(\h+1))}}  \big|_{\g^{(s,\ve)(k-1)}(0,s)} ds \geq 0
\ . \end{multline}
From \eqref{condi} and the definition of  $Y$,  if  the  needle variation  satisfies (a) and (b),    both  integrals in \eqref{49} are zero  and the inequality is satisfied.
\end{pf}
\begin{rem}\label{remarkone}  By definition, for any control curve $u_o(t)$ in  $\wh K$ and  any $\sfs_o \in \cA_{\text{\it init}}$, there exists a uniquely associated good pair $ U^{(u_o, \sfs_o)} \= \left(u_o(t), \s_o = (\sfs_o, \wt{\sfs_o}^{(u_o, \sfs_o)})\right)$. Then   for any good $\wh \cK$-controlled curve $\g_o$
and any  $(\t_o, \o_o) \in (0,T]\times K$, it is possible to construct  a  smoothed needle variation   $ {\cneedle}^{(\t_o, \o_o,  \Sigma, \ve_o)} (\g_o)$  satisfying  both conditions  of Lemma \ref{existenceofgood}.  This means that {\it for  any good controlled curve $\g_o$ and any $(\t_o, \o_o) \in (0,T]\times K$, there is a good needle variation $ {\cneedle}^{(\t_o, \o_o,  \Sigma, \ve_o)} (\g_o)$  associated with $(\t_o, \o_o)$}. 
\end{rem}
Remark  \ref{remarkone} and Theorem \ref{POMPONE} easily imply the following   $\cC^\infty$  version of Theorem \ref{main1}.
\begin{theo}  \label{firstversion} Let  $f = (f^i)$ and $\sC$  of class $\cC^\infty$ and 
 $U_o = (u_o(t), ({\mathsf s}_o,\wt{\mathsf s}^{(u_o, \mathsf s_o)})) \in \cK_{\text{\rm good}}$  with  $u_o(t)$ smooth. 
If $U_o$  is  an optimal control, then 
 \eqref{main11} holds   for any $(\t_o, \o_o) \in (0, T) \times  K$. \par
More precisely,   for any such $(\t_o, \o_o) $,   there exist  constants $M$,  $\widecheck M$, which  depend  on $f$,  $\sC$,  $\t_o$ and $\o_o$,  such that   if
\beq  \label{bad-bis} \k_o \= \sH^{(u_o, {\mathsf s}_o, \t_o)}(\o_o)  - \sH^{(u_o, {\mathsf s}_o, \t_o)}(u_o(\t_o))   > 0\ ,\eeq
then  there is a  needle modification $\underline u(t)$   of $u_o(t)$,    with associated smoothed needle modification  $\widecheck{\underline u}(t)$,  such that 
\beq 
\begin{split}
   \sC(j^{k-1}_{t = T}(x^{(\underline u, \mathsf s_o)}))  \leq   \sC(j^{k-1}_{t = T}(x^{(u_o, \mathsf s_o)})) -  \frac{\k_o}{4 M} < \sC(j^{k-1}_{t = T}(x^{(u_o, \mathsf s_o)}))  \ , \\
   \sC(j^{k-1}_{t = T}(x^{(\widecheck {\underline u}, \mathsf s_o)}))    \leq   \sC(j^{k-1}_{t = T}(x^{( u_o, \mathsf s_o)})) -  \frac{\k_o}{4 \widecheck M} < \sC(j^{k-1}_{t = T}(x^{(u_o, \mathsf s_o)}))  \ . 
   \end{split}
\eeq
\end{theo}
\begin{pf}  If  $U_o = (u_o(t), ({\mathsf s}_o,\wt{\mathsf s}^{(u_o, \mathsf s_o)}))$ is good,  the  function
\eqref{pontr-funct}  is equal to 
\begin{multline} \cP^{((\sfs_o, \wt \sfs^{(u_o, \sfs_o)}), u_o, \t_o)}(\o) 
=  - p^{(u_o, \sfs_o)}_i(\t_o) \ x^{(u_o, \sfs_o)i}_{(k)}(\t_o) +\\
+  p^{(u_o, \sfs_o)}_i(\t_o) f^i(\t, x^{(u_o, \sfs_o)}(\t_o), x^{(u_o, \sfs_o)}_{(1)}(\t_o), \ldots, x^{(u_o, \sfs_o)}_{(k -1)}(\t_o), \o)= \\
 =  - p^{(u_o, \sfs_o)}_i(\t_o) \ x^{(u_o, \sfs_o)i}_{(k)}(\t_o) + \sH^{(u_o,\sfs_o, \t_o)}(\o)\ .\end{multline}
Thus $\o$ is a maximum point for  $\cP^{((\sfs_o, \wt \sfs^{(u_o, \sfs_o)}), u_o, \t_o)}$  if and only if  it is 
a maximum point for  $\sH^{(u_o,\sfs_o, \t_o)}$. The claim  then follows from Theorem \ref{POMPONE},   Lemma \ref{corolletto}  and Remark \ref{remarkone}. 
\end{pf}
\subsection{The proof of Theorem \ref{main1}} \label{finalproof}
First of all, we assume the following condition, which causes  no loss of generality.  
Let  $\overline{\bB}_{\wt R} \subset \bR^{m}$  
 and $\overline{\bB' }_{\wt R} \subset \bR^{m k}$ be two  closed balls centred at the origin and  of radius $\wt R$,  which contain  the compact sets  $K \subset \bR^m$ and 
$$K^{\leq (k-1)} \= K \times K^{(1)}\times \ldots K^{(k-1)}\subset \bR^{m k}\ ,$$
respectively.  Then, we   set   $\wh K  \= \overline{\bB}_{2\wt R}$, $\wh K'  \= \overline{\bB}_{3\wt R}$ and  $\wh K'' \= \overline{\bB}_{4\wt R}$.  We also assume that {\it  $f$  is   extended to  a map on $\O \times \wh K''$, which still  satisfies   $(\a)$ and $(\b)$}. 
\par
\smallskip
 As a   preliminary step, we  need the following  lemma. \par
\begin{lem}\label{lemma32}  For each pair   $(u_o(t), \mathsf s_o)$,   with $\mathsf s_o \in \cA_{init} $ and  $u_o: [0, T] \to \wh K \subset \bR^m$ satisfying the  condition $(\g)$ of the  Introduction,   there exist: 
\begin{itemize}[leftmargin = 15pt]
\item[--] A  unique  solution  $x^{(u_o, \mathsf s_o)}: [0,T] \to \bR^n$  to     \eqref{controlledEL-1*} with initial condition    $j^{k-1}_{t = 0} (x^{(u_o, \mathsf s_o)}) = \mathsf s_o$.  If 
$k = 1$, this  solution 
is   $\cC^{0}$  with bounded   measurable  first  derivative. If $k \geq 2$, the solution is piecewise $\cC^{2k-2}$.
\item[--] A unique solution $p^{(u_o, \mathsf s_o)}: [0,T] \to \bR^n$   to   \eqref{controlledEL-1bis*}  with terminal conditions \eqref{31-ter} -- \eqref{33-ter}.  This solution is 
of  class  $\cC^{k-1}$ and with bounded   measurable  $k$-th derivative.  
\end{itemize}
\end{lem}
\begin{pf}   Consider the  auxiliary  variables $x^i_\ell$,   $1 \leq \ell \leq k-2$,  $1 \leq i \leq n$,   and   the first order  differential problem
\beq\label{diffe-bis}
\begin{split} & \frac{d x^i}{dt} = x^i_1\ ,\quad \frac{d x^i_1}{d t}= x^i_2\ ,\qquad\ldots \ ,\qquad \frac{d x^i_{k-2}}{dt} = x^i_{k-1}\ ,\\
 & \frac{d x^i_{k-1}}{dt} = f^i(t, x^j, x^j, \ldots, x^j_{k-1}, u^a_o(t)) \ ,
\end{split}
\eeq
with  initial  conditions $(x^i, x^i_\ell)|_{t = 0}$  determined   by the  jet    $\sfs_o = (x^i, x^i_{(\ell)})|_{t = 0}$. This problem is  equivalent to the  system   \eqref{controlledEL-1*} with  initial condition $j^{k-1}_{t = 0} (x^{(u_o, \mathsf s_o)}) =\mathsf s_o$.  Hence the existence and uniqueness of a $\cC^{k-1}$ solution   $x^{(u_o,  \sfs_o)}$   is  a consequence of   a well-known result on first order
differential systems  with  control parameters in normal form (see e.g. \cite[Th. 3.2.1]{BP}).  The $(k-1)$-th derivative  of this solution is absolutely continuous with bounded derivative. Moreover, for $k \geq 2$, on each subinterval on which $u_o(t)$ is $\cC^{k-1}$,   the curve $x_{k-1}(t)$ is  $\cC^{k-1}$. It follows that $x^{(u_o,  \sfs_o)}$  is   piecewise $\cC^{2k-2}$. 
\par
  The existence  and uniqueness of $p^{(u_o, \sfs_o)}$  is   checked  by considering    \eqref{controlledEL-1bis*}   as  a system of equations on the functions $p_j(t)$, depending on  the control curve  $\gu(t) \= (u_o(t), j_t^{2k-2}(x^{(u_o,  \sfs_o)}))$ taking values in $K \times J^{2k-2}(\sQ|\bR)$. Since the curve $u_o(t)$ is bounded and measurable and $x^{(u_o,  \sfs_o)}$ is of class $\cC^{2 k -2}$, the result follows from the above mentioned facts on  systems with  control parameters. 
\end{pf}

We are   now ready to   prove the following  crucial result,  which  implies Theorem \ref{main1}. \par
\begin{theo}  \label{firstversion-bis} Let $u_o: [0, T] \to K$ be a  measurable control curve and $\sfs_o \in \cA_{\text{\it init}}$, as in Theorem \ref{main1}, 
and assume that $(\t_o, \o_o) \in (0, T) \times  K$ is a pair, in which   $\t_o$ satisfies  the following  condition:
\begin{itemize}[leftmargin = 15pt]
\item  if $k = 1$, the time $\t_o$ is  one of the   points for  which 
\beq \label{BPcondition}  \lim_{\ve \to 0^+} \frac{1}{\ve} \int_{\t_o - \ve}^{\t_o } \left| f(t, x^{(\sfs_{o},  u_{o})}(t),  u_o(t)) - f(\t_o, x^{(\sfs_{o},  u_{o})}(\t_o),  u_o(\t_o) )\right|dt = 0\ , \eeq
 i.e.  $\t_o$ is a density (Lebesgue) point of the map $f(t, x(t), u(t))$; 
\item if $k \geq 2$,  $\t_o$ is an inner point of a subinterval $I \subset [0, T]$ on which $u_o(t)$ is  $\cC^{k-1}$. 
\end{itemize}
Then  there is  a constant $M = M_{(f, \sC,\t_o, \o_o)} >0$, depending  on  $f$, $\sC$, $\t_o$,  $\o_o$,  such that   if
\beq  \label{bad-bisbis} \k_o \= \sH^{(u_o, {\mathsf s}_o, \t_o)}(\o_o)  - \sH^{(u_o, {\mathsf s}_o, \t_o)}(u_o(\t_o))   > 0\ ,\eeq
then  there is a   needle modification $u'(t)$   of $u_o(t)$ with  peak time $\t_o$ and ceiling $\o_o$, 
satisfying 
\beq  \label{eccolo}
\begin{split}
   \sC(j^{k-1}_{T}(x^{(u', \sfs_o)}))  \leq   \sC(j^{k-1}_{T}(x^{(u_o, \sfs_o)}))) -  \frac{\k_o}{4 M} < \sC(j^{k-1}_{T}(x^{(u_o, \mathsf s_o)}))) 
   \end{split}
\eeq
and  $U_o = (u_o(t), {\mathsf s}_o)$  cannot be   an optimal control. 
\end{theo}
\begin{pf}  The proof is based on a three-step approximation procedure,   which allows   inferring the 
 theorem   from its previous  ``smooth'' version,   Theorem \ref{firstversion}.  For reader's convenience,  here is an outline of the 
arguments which we are going to use  in  case  $k \geq 2$. 
\begin{itemize}[leftmargin= 20pt]
\item[(I)] First we  introduce a  one-parameter family  of {\it globally} $\cC^{k-1}$ curves $v^{(\h)}:[0, T] \to \wh K$, which tends to the piecewise $\cC^{k-1}$ curve $u_o(t)$ for $\h \to 0$ with respect to the  distance \eqref{distanza}.  This family is constructed in such a way  that: (1) $v^{(\h)}(t)$ coincides with $u_o(t)$ on a neighbourhood of $\t_o$ for any $\h$;  (2)  the associated controlled curves with initial datum $\sfs_o$  tend to the curve $x^{(u_o, \sfs_o)}$ in $\cC^{2k - 2}$ norm. 
\item[(II)] Second, for any  fixed value $\d_1$  for the parameter $\h$ and the corresponding  curve $v^{( \h = \d_1)}(t)$ of step  (I), we consider a  one-parameter family  of  polynomials  (in the $t$-variable)  $v^{(\d_1, \wt \h)}:[0, T] \to \wh K'$,  which   converge to  $v^{(\d_1)}(t)$  in the $\cC^{k-1}$-norm for $ \wt \h \to 0$ and such that, the associated controlled curves $x^{(v^{(\d_1,  \wt  \h)}, \sfs_o)}(t)$ and $p^{(v^{(\d_1,  \wt \h)}, \mathsf s_o)}(t)$, defined in Lemma \ref{lemma32}, converge to the curves  $x^{(v^{(\d_1)}, \sfs_o)}(t)$ and $p^{(v^{(\d_1)}, \mathsf s_o)}(t)$ in the $\cC^{2k-2}$ and $\cC^{k-1}$ norm, respectively. 
\item[(III)]Third,  for any  fixed choice  of  $\d_1$, $\d_2 > 0$,  we consider a one parameter family  of smooth functions $f^{(\d_1,\d_2,  \wt{\wt\h})}:  \O  \times \wh K' \to \bR^n$,  which converges in the  $\cC^{k-1}$ norm to the function $f(t, j^{k-1}_t(x), u)$  for $  \wt{\wt\h} \to 0$. The family is constructed in  such a way  that the  smooth solutions $x^{(\d_1,\d_2,  \wt{\wt\h})}(t)$ to the  constraint given by  $f^{(\d_1, \d_2,  \wt{\wt\h})}$,   the polynomial control curve $v^{(\d_1,  \wt \h = \d_2)}(t)$ and    the initial datum $\sfs_o$, tend in the $\cC^{2k-2}$ norm to the  solution of  the   constraint determined by  $f$,  $v^{(\d_1, \d_2)}(t)$ and   $\sfs_o$.
\end{itemize}
After  these preliminary  constructions,  we  show that: 
\begin{itemize}[leftmargin = 20 pt]
\item[(a)] For any  triple $(\d_1, \d_2, \d_3)$,  the  controlled   curve    $x^{(\d_1, \d_2, \d_3)}(t)$, determined by the {\it smooth}  constraint given by $f^{(\d_1, \d_2, \d_3)}$,  the {\it smooth}  control curve  $v^{(\d_1, \d_2)}(t)$ and   the  initial datum $\sfs_o$ satisfies the hypotheses of Theorem \ref{firstversion}.  This implies the existence of  an appropriate  needle modification  $v^{(\d_1, \d_2,  \underline \ve) }(t)$ of  the control curve  $v^{(\d_1, \d_2)}(t)$ (with $\underline \ve = \underline \ve(\d_3)$ depending on $\d_3$),  which determines a controlled curve with a  smaller terminal cost.  
\item[(b)] We then  show  that the $\d$-parameterised family of  needle modifications  $v^{( \d, \d, \underline \ve(\d))}(t)$ converges  in  the $L^1$ norm  to a  needle modification $u'(t)$ of $u_o(t)$, whose associated controlled curve $x^{(u', \sfs_o)}(t)$ gives a terminal cost  satisfying \eqref{eccolo} 
\end{itemize}
\par
\smallskip
The scheme of the proof for the case $k = 1$  is  similar, but requires a preparatory additional  step.  Before starting with  the whole construction,  we   replace $u_o(t)$ by 
the   curve  $u^{(\d_0)}:[0, T] \to \wh K$, which  is  constant and equal to  the value $u_o(\t_o)$ on the interval $[\t_o - \d_0, \t_o + \d_0]$ and is equal to    $u_o(t)$  at all other points.  Then,   working with the  modified   curve $u^{(\d_0)}$, we perform the analogs of the  three steps (I), (II) and (III). This leads to the   construction of    a  $\d$-parameterised family  of needle modifications of  the smooth control curves, which  converge  in  $L^1$-norm   to  a needle modification  $u'(t)$ of $u_o(t)$ with  associated  controlled curve  $x^{(u', \sfs_o)}$  with terminal cost satisfying \eqref{eccolo}, as desired. 
  \par
\medskip
Let us now proceed with the proof  for the case $k \geq 2$.  Let us  consider   the  following {\it derived system of order $2k-1$} associated with \eqref{controlledEL-1*}.  It is the system of equations which can be obtained from  \eqref{controlledEL-1*}
by differentiating  $k-1$ times  with respect to $t$  and by replacing any   $k$-th order derivative $x_{(k)}^i$ by $f^i(t, j^{k-1}_t(x), u_o(t))$  at all places:  
\begin{equation} \label{firstder-bis} 
\begin{split}
& x^i_{(k)} = f^i\ ,\\
&x_{(k+1)}^i = \frac{\p f^i}{\p t}   + \sum_{r = 0}^{k-2} x^j_{(r+1)}  \frac{\p f^i}{\p x^j_{(r)}}    + f^j \frac{\p f^i}{\p x^j_{(k-1)} } 
 + \frac{\p f^i}{\p u^a} u^a_{o(1)}\ ,\\
& \vdots\\
 & x_{(2 k-1)}^i = \frac{\p^{k -1} f^i}{\p t^{k -1}}   + \ldots\ . 
 \end{split} \end{equation}
Let us synthetically  denote these equations   by 
\beq x^i_{(\ell)} = F^i_{(\ell)}(t, j^{k-1}_t(x), j^{k-1}_t(u_o))\ , \qquad k \leq \ell \leq 2 k -1\ .\eeq
Then, consider 
the analogue of the first order  system \eqref{diffe-bis} that  gives the reduction to the first order of the last line of
\eqref{firstder-bis}. 
Using  the  shorter notation
$y \=  (x^i_{j})$,  this system can be  written as  
\beq \label{sistemino} \dot y = g(t,y(t), u_o(t), u_{o(1)}(t), \ldots u_{o(k-1)}(t))\eeq
where $g$ is  a   map  $g: \wt \O \times  K^{\leq (k-1)} \to \bR^{n(2k-1)}$,  for an appropriate open set $\wt \O \subset  \bR^{n (2 k-1) + 1}$,  which is  uniquely determined by  $f: \O \times K \to \bR^{n}$. By the above described  technical assumptions on $f$, we may assume that  $g$ is  defined on a larger domain $\wt \O  \times \wh K^{\leq(k-1)}$ with  $ \wh K^{\leq(k-1)} \supsetneq K^{\leq (k-1)}$.  Such extension is  continuously differentiable in each variable $y^A$.  \par
Finally,   for any   pair  $(u(t), \mathsf s)$,   let us  denote by $y^{(u, \mathsf s)}: [0, T] \to \bR^{n(2k-1)}$ the   unique  solution to \eqref{sistemino}, which is  controlled by   $u(t)$  and  with   the initial condition  that correspond  to the initial condition  $\sfs$ for  $x(t)$.\par
\smallskip
We may now  construct the one-parameter family   described in (I), using the following 
\begin{lem} \label{lemma42-1} Let $(t_0 = 0, t_1)$, $(t_1, t_2)$, \ldots,   $(t_{P-1}, t_P = T)$ be  the intervals on which $u_o(t)$ is $\cC^{k-1}$. Given  $\h> 0$  there exists a $\cC^{k-1}$ curve  $v^{(\h)}: [0, T] \to \wh K$, which coincides with $u_o(t)$ on the  subintervals $(t_{i-1} + \frac{\h}{4P}, t_i - \frac{\h}{4 P})$ (and, in particular, on some neighbourhood of $\t_o$)  and  such that 
\beq\label{firststep-1}
\|x^{(u_o, \mathsf s_o)}(t) - x^{(v^{(\h)}, \mathsf s_o)}(t)\|_{\cC^{2k-2}} \ ,\ \  \|p^{(u_o, \mathsf s_o)}(t) - p^{(v^{(\h)}, \mathsf s_o)}(t)\|_{\cC^{k-1}}< \h
\ .\eeq 
\end{lem}
\begin{pf} It is  almost immediate to realise   that, for any choice of $\d$,  there exists  a curve  $\wt v(t)$ which  is $\cC^{k-1}$ over the whole domain $[0,T]$ and  coincides with  $u_o(t)$ on the subintervals $(t_{i-1} + \frac{\d}{4P}, t_i - \frac{\d}{4 P})$. This implies  that   all distances  $  \dist(u_o, \wt v) $, $  \dist(u_{o(1)}, \wt v_{(1)}) $, \ldots, $  \dist(u_{o(k-1)}, \wt v_{(k-1)}) $ are less than $\d$.  
Now, 
 by the assumptions on  $f$ and on its derivatives
  (which completely determine the  function $g$ in \eqref{sistemino}), 
  there  exists a  unique solution to the reduced-to the-first-order system \eqref{sistemino}   for the pair $(\wt v, \mathsf s_o)$.  Furthermore, 
by  Lemma \ref{normallemma}, there  are constants 
 $\r >0$  and $\gc$ (depending on $f = (f^i)$ and on a cut-off function $\f$ as described  in the statement of that lemma) such that, if  $\d \leq \r$, then
$$\sup_{t \in [0,T]} | y^{(u_o, \mathsf s_o)}(t) - y^{(\wt v, \mathsf s_o)}(t)| < \gc \dist(u_o,\wt  v) < \gc \d$$
 with   $\dist(\cdot, \cdot) $  given by  \eqref{distanza}.
Selecting a $\d_\h$ such that $\d_\h <   \min\big\{ \frac{\h}{\gc}, \r, \h \big\}$,  we get   
$\sup_{t \in [0,T]} | y^{(u_o, \mathsf s_o)}(t) - y^{(\wt v, \mathsf s_o)}(t)| < \h$, meaning that $v^{(\h)}(t) \= \wt v$ satisfies  the first  upper bound in   \eqref{firststep-1}.   By considering a possibly smaller $\d_\h$ also the second bound  is satisfied. This is because  $p^{(\wt v, \mathsf s_o)}$ is a solution of a system of   controlled differential equations, where the  controls are given  by  the curve $\wt v(t)$  and  the curve of  the $(2k-2)$-jets of    $x^{(\wt v, \mathsf s_o)}(t)$. 
\end{pf}
\par
Let us now fix a control curve $v^{(\d_1)}$ as in the previous lemma. The family of polynomials control curves  in (II) is constructed using the next 
\begin{lem}  \label{lemma42}   Given    $\h > 0$,
 there exists  a polynomial curve $v^{(\d_1, \h)}(t)$ in $\wh K' \supset \wh K$   with 
 \beq \label{att_o}|v^{(\d_1, \h)}(\t_o) - v^{(\d_1)}(\t_o)|  = |v^{(\d_1, \h)}(\t_o) - u_o(\t_o)| < \h\eeq
 and  such that the   solutions $x^{(v^{(\d_1)}, \sfs_o)}$ and $x^{(v^{(\d_1, \h)}, \sfs_o)}$ to the differential problem \eqref{diffe-1-ter} (which is the same of  \eqref{controlledEL-1*}) and the associated curves $p^{(v^{(\d_1)}, \sfs_o)}(t)$  and $p^{(v^{(\d_1, \h)}, \sfs_o)}(t)$, defined in Lemma \ref{lemma32},  satisfy  
\beq \label{4.15}
\|x^{(v^{(\d_1)}, \mathsf s_o)}- x^{(v^{(\d_1, \h)}, \mathsf s_o)}\|_{\cC^{2k-2}} \ , \quad \|p^{(v^{(\d_1)}, \mathsf s_o)} - p^{(v^{(\d_1, \h)}, \mathsf s_o)}\|_{\cC^{k-1}} <  \h\ .
\eeq
\end{lem}
\begin{pf}   
By a well known result on interpolation of continuous functions (see e.g. \cite[Thm. 7.1.6]{Ph}),    we may consider a family of Bernstein polynomials   converging  to $v^{(\d_1)}(t)$  in   $\cC^{k-1}$ norm.  Thus, for any choice of  a sufficiently small $\wt \d > 0$, we may select a polynomial $\wt v^{(\d_1)}(t)$   which satisfies \eqref{att_o}, takes  values in $\wh K'$ and  such that 
\beq \label{delta} \int_{[0,T]} |v^{(\d_1)}_{(\ell)}(t) -   \wt v_{(\ell)}^{(\d_1)}(t)| dt  < \wt \d\ ,\qquad 0 \leq \ell \leq k-1\ .\eeq
Now,  by  \cite[Prop. 3.2.5 (i)]{BP},  there exists a $ \wt \d> 0$ (depending on $\h$)  such that   \eqref{delta} implies that the corresponding solutions $y^{(v^{(\d_1)}, \mathsf s_o)}(t)$ and $y^{( \wt v^{(\d_1)}, \mathsf s_o)}(t)$ of  the reduced-to-the-first-order system \eqref{sistemino}  satisfy the inequality 
$\sup_{t \in [0,T]} | y^{(v^{(\d_1)}, \mathsf s_o)}(t) - y^{( \wt  v^{(\d_1)}, \mathsf s_o)}(t)| < \h$.
 Hence,  the polynomial $v^{(\d_1, \h)}(t) \= \wt v^{(\d_1)}(t)$  is such that   the first upper bound in \eqref{4.15} holds.  Using the same  arguments  for the equations on the $p(t)$,  the second bound of \eqref{4.15} can be  satisfied as well.   
 \end{pf}
\par
It is now the turn to present  the family of functions described in (III). \par
\smallskip
\begin{lem}  \label{lemma43}  Let $v^{(\d_1, \d_2)} :[0,T]  \to \wh K'$ be  one of the polynomials described in Lemma \ref{lemma42}, converging to the $\cC^{k-1}$ control curve $v^{(\d_1)}$.  Let  also  $x^{(v^{(\d_1, \d_2)}, \sfs_o)}(t)$ be the  unique solution  to  \eqref{diffe-1-ter}  determined by  the pair $(v^{(\d_1, \d_2)}(t), \mathsf s_o)$, as discussed  in Lemma \ref{lemma32}.  \par
Then there is a  $\h_o > 0$ such that for  $\h \in (0, \h_o]$ there are   $\cC^\infty$  maps $ f^{(\d_1, \d_2, \h)}:  \O  \times \wh K'' \to \bR^n$  satisfying the following conditions: they converge uniformly on compacta to $f$ together with all partial derivatives up to order $k-1$ for $\h \to 0$   and, for each $\h$,  the unique solution  $x^{(v^{(\d_1, \d_2)}, \sfs_o;\h)}(t) $  to  the differential problem
\beq\label{diffe-hat} x_{(k)} =   f^{(\d_1, \d_2, \h)} \left(t, j^{k-1}_t(x), v^{(\d_1, \d_2)}(t)\right)\ ,\qquad j^{k-1}_{t=0}(x) = \mathsf s_o
\eeq
and the associated curves $p^{(v^{(\d_1, \d_2; \h)}, \mathsf s_o)})(t)$  defined in Lemma \ref{lemma32} satisfy
\beq \label{desired} 
\begin{split}  &  \|x^{(v^{(\d_1, \d_2)}, \sfs_o; \h)} - x^{(v^{(\d_1, \d_2)}, \sfs_o)}\|_{\cC^{2k-2}} < \h\ ,\qquad 
 \|p^{(v^{(\d_1, \d_2; \h)}, \mathsf s_o)} - p^{(v^{(\d_1, \d_2)}, \mathsf s_o)}\|_{\cC^{k-1}} <  \h\ ,\\
& \left| f^{(\d_1, \d_2, \h)} \big(t, j^{k-1}_t(x^{(v^{(\d_1, \d_2)}, \sfs_o; \h)}), v^{(\d_1, \d_2)}(t)\big)  {-} f\big(t, j^{k-1}_t( x^{(v^{(\d_1, \d_2)}, \sfs_o)}), v^{(\d_1, \d_2)}(t)\big) \right| < \h \\
& \hskip 10 cm    \text{for any}\ t \in [0,T]\ .
\end{split}
\eeq
\end{lem}
\begin{pf}   Since $v^{(\d_1, \d_2)}(t)$ is  a polynomial,  it satisfies the condition ($\g$) of the Introduction and  the corresponding solution $x^{(v^{(\d_1, \d_2)}, \sfs_o)}(t)$ to \eqref{diffe-1-ter} 
 determines a    curve  of  jets and controls
\begin{multline*} \g^{(v^{(\d_1, \d_2)}, \sfs_o)(k-1)}(t) \=  \left(t, x^{(v^{(\d_1, \d_2)}, \sfs_o)}(t), x^{(v^{(\d_1, \d_2)}, \sfs_o)}_{(1)}(t), \ldots, x^{(v^{(\d_1, \d_2)}, \sfs_o)}_{(k-1)}(t), v^{(\d_1, \d_2)}(t)\right) \in \\
\in  J^{k-1}(\sQ|\bR)|_{[0,T]} \times \wh K'\ ,\end{multline*}
which is of course  continuous and with compact image. Hence, there is a   $\h_o > 0$ such that the set 
\begin{multline*} \Pi_{\h_o}  \=\bigg \{ (t, \mathsf s_t, u)   \in J^{k-1}(\bR^n|\bR)|_{[0,T]} \times \bR^m\ : \\
 |\mathsf s_t - j^{k-1}_t(x^{(v^{(\d_1, \d_2)}, \mathsf s_o)})| \leq  \h_o \qquad \text{and}\qquad  |u - v^{(\d_1, \d_2)}(t)|\leq  \h_o\ \bigg\}
\end{multline*}
is compact and  with   $u \in \overline{\bigcup _{v \in \overline B_{2 \wt R}} \bB_{\h_o}(v)}  \subset  \wh K'' = \overline{B_{4 \wt R}}$.   By the assumptions ($\a$) and ($\b$) on $f$,   there is  a   constant $\gL$ such that  
\beq \|f\|_{\cC^{k-1}( \Pi_{\h_o})} \ , \quad
  \max_{\smallmatrix 
w_0 + |w| + |\ell| = k, \\ 
0 \leq w_0 +  |  w | \leq k -1,\\
(t, \mathsf s_t, u) \in \Pi_{\h_o}\ 
\endsmallmatrix} \left\{\left| \left. \frac{\p^{w_0 + | w|  + |\ell|} f}{  (\p t)^{w_0} (\p u)^w (\p x^i_{(r)})^\ell}
\right|_{(t, \mathsf s_t, u)}
\right|\right \}  
< \gL\ ,\eeq
where we  denote  
$$\frac{\p^{|w|}}{(\p u)^w} = \frac{\p^{|w|}}{(\p u^1)^{w_1} \ldots (\p  u^m)^{w_m}}\ ,\qquad  \frac{\p^{|\ell|}}{(\p x^i_{(r)})^\ell} =  \frac{\p^{|\ell|}}{(\p x^1_{(0)})^{\ell_{1|0}}  \ldots  
(\p x^n_{(k-1)})^{\ell_{n|k-1}}}\ ,$$
with $w $ and $\ell$  multiindices  $ w = (w_1, \ldots, w_m)$ and $\ell = (\ell_{1|0},\ell_{1|1},  \ldots, \ell_{1|k-1},\ell_{2|0}, \ell_{2|1}, \ldots, 
\ell_{2|k-1}, $ $\ldots,  \ell_{n|0}, \ell_{n|1}, \ldots, \ell_{n|k-1})$.
Moreover, for any    $\d > 0$, there is   $\wh f^{(\d_1, \d_2, \d)} \in \cC^\infty(J^{k-1}(\bR^n|\bR)|_{[0,T]} \times \wh K'')$ such that 
\beq \label{eta}  
\| \wh f^{(\d_1, \d_2, \d)} - f\|_{\cC^{k-1}(\Pi_{\h_o})} 
< \d  \eeq
(it is a  consequence of a classical approximation procedure; see e.g.  \cite[Ch. 15]{Tr}). 
\par
We now want to prove that there exists a constant $ \gC$,  depending  on $k$, $n$, $\gL$ and
 \beq \label{sup} \sup \bigg\{|j^{k-1}_t(x)|\ : \ j^{k-1}_t(x) \in \bB_{\h_o}(j^{k-1}_t( x^{(v^{\d_1, \d_2)}, \sfs_o}))\ ,\ t \in [0,T]\bigg\}\ ,\eeq  with the following property: for  any $\h \in (0, \h_o]$ with $\h_o < \frac{e^{-\gC \max\{\gL,1\}}}{4}$, 
if $\d_\h$ is sufficiently small, 
then 
the corresponding function $f^{(\d_1, \d_2, \h)} \= \wh f^{(\d_1, \d_2, \d_\h)}$ satisfies    \eqref{eta}
and  the associated  solution $x^{(\h)}(t) \= x^{( v^{(\d_1, \d_2)} \mathsf s_o; \d_\h)}(t)$  to the system  \eqref{diffe-hat}    satisfies    \eqref{desired}.\par
  To see this,  consider the two  systems of first order  of the form \eqref{diffe-bis} (obtained by introducing the    auxiliary variables $y^i_r$),  which  correspond to the  derived differential systems of  order $2k-1$ associated with the system \eqref{diffe-1-ter} and the system $x_{(k)} =   f^{(\d_1, \d_2, \d)} \left(t, j^{k-1}_t(x), v^{(\d_1, \d_2)}(t)\right)$.  Then, for each choice of $\d$, let  us denote by 
 \begin{multline*} y^{(v^{(\d_1, \d_2)}, \mathsf s_o)}(t) = (y^i(t), y^i_1(t), \ldots,  y^i_{2k-2}(t)) \ ,\\
  y^{(v^{(\d_1, \d_2)}, \mathsf s_o; \d)}(t) = ( y^{\d i}(t),  y^{\d i}_1(t), \ldots,   y^{\d i}_{2k-2}(t))
  \end{multline*}
the solutions to  such two systems, corresponding to  the solutions 
$x^{(v^{(\d_1, \d_2)}, \mathsf s_o)}(t)$ and  $ x^{(v^{(\d_1, \d_2)}, \mathsf s_o; \d)}(t)$, respectively. 
We synthetically denote the two systems of equations,  of which  they are solutions,  by 
$$\dot{y}^i_{\ell-1} =  g^{i}_\ell(t,  y^j_m, v^{(\d_1, \d_2)}(t))\ ,\qquad \dot{ y}^i_{\ell-1} =  g^{(\d)i}_\ell(t,  y^j_m, v^{(\d_1, \d_2)}(t))$$
with  $1 \leq \ell \leq 2k - 1 $. 
\par
By construction  $g^i_m(t, y^j_r, v^{(\d_1, \d_2)}(t)) =  g^{(\d)i}_m(t, y^j_r, v^{(\d_1, \d_2)}(t)) = y^i_{m}$ for any   $0 \leq m \leq 2k-2$.  On the other hand, the functions $g^i_{2k-1}(t, y^j_r, v^{(\d_1, \d_2)}(t))$ and  $ g^{(\d)i}_{2k-1}(t, y^j_r, v^{(\d_1, \d_2)}(t))$ are in general  different, because they are  given  by the $(k-1)$-th order total derivatives of the functions $f$ and  $f^{(\d_1, \d_2, \d)}$, respectively, evaluated at the points $(t, y^j_m, v^{(\d_1, \d_2)}(t))$. 
\par
The initial values of the curves  $y(t) \= y^{(v^{(\d_1, \d_2)}, \mathsf s_o)}(t) $   and $y^{(\d)}(t) \=  y^{(v^{(\d_1, \d_2)}, \mathsf s_o; \d)}(t)$  are denoted by 
 $\overline y = (\overline y^i,  \overline y^i_1, \ldots,  \overline y^i_{2k-2})$ and  $\overline y^\d = (\overline y^{\d i},  \overline y^{\d i}_1, \ldots,  \overline y^{\d i}_{2k-2})$, respectively.  Note that the components  of  $\overline y $ and $\overline y^\d$,  determined by the derivatives up to order $k -1$ at $t = 0$ of the two curves,  are the same and uniquely determined by $\sfs_o$.  By construction of the derived system of order $2k-1$, the remaining components of   $\overline y $ and $\overline y^\d$ might be different, but  also  such that   $\overline y^\d \to \overline y$ for $\d \to 0$. \par
 \smallskip
Let us  denote 
$ z^i_\ell(t) =  y^{\d i}_\ell(t) - y^i_\ell(t)$ for any $0 \leq \ell \leq 2 k-2$. We observe that,   for any $t \in [0, T]$
\begin{multline*}\dot z^i_r(t) = z^i_{r+1}(t)\ \text{if} \ \ 0 \leq r \leq 2k-3\ ,\\ \dot z^i_{2k-2}(t) =  g^{(\d)i}_{2k-1}(t,  y^\d(t), v^{(\d_1, \d_2)}(t)) - g^i_{2k-1}(t,  y(t), v^{(\d_1, \d_2)}(t))
\end{multline*}
 and therefore
 \begin{align*}
 &  |\dot z^i_r(t)|  \leq \sum_{\ell, i } | z^i_\ell(t)|\ ,\qquad  0 \leq r \leq 2 k-3\ ,
 \end{align*}
 \begin{align*}
& |\dot z^i_{2k-2}(t)| \leq \bigg| g^{(\d)i}_{2k-1}(t,  y^\d(t), v^{(\d_1, \d_2)}(t)) - g^i_{2k-1}(t,   y^\d(t), v^{(\d_1, \d_2)}(t))\bigg|  + \\
&\hskip 3 cm + 
 \bigg| g^i_{2k-1}(t,  y^\d(t), v^{(\d_1, \d_2)}(t)) - g^i_{2k-1}(t,  y(t), v^{(\d_1, \d_2)}(t))\bigg|  
\leq \\
& \hskip 10 cm \leq 
\text{Const} (\d  + \gL\sum_{\ell, i } | z^i_\ell(t)|)\ ,\\
& \frac{d}{dt} \sum _{i,\ell}  |z_\ell^i(t) | \leq \sum_{i,\ell} |\dot z^i_\ell(t) | \leq \gC( \d  + \max\{\gL, 1\}\sum_{i,\ell}  | z^i_\ell(t)|)\ ,
\end{align*}
where  $ \text{Const}$ is a constant,  which depends only on $\gL$ and  \eqref{sup}, and 
  $\gC \=   n (2k-1) \text{Const}$.  \par
\smallskip
Hence if  we  take $\d = \d_\h$  so that
   $\sum_{i,\ell}  | z^i_\ell(0)| \leq \h^2 \leq \h \h_o$ and 
$ \d_\h \leq \frac{\h e^{- \gC { \max\{\gL, 1\}  T}}}{4  \gC T }$,  
then by Gronwall's inequality we obtain that 
\begin{multline*}  |j^{2k-2}_t(x^{(v^{(\d_1, \d_2)}, \sfs_o; \d_\h)}) -  j^{2k-2}_t(x^{(v^{(\d_1,\d_2)}, \mathsf s_o)})|  \leq \sum _{i,\ell}  |z_\ell^i(t) |\leq  \\
\leq (  \d_\h \gC T  + \sum_{i,\ell}  | z^i_\ell(0)|)) e^{ \gC \max\{\gL, 1\} T}  
\leq  \d_\h \gC  T  e^{ \gC\max\{\gL, 1\} T}   + \h \h_o e^{ \gC\max\{\gL, 1\} T} \leq \frac{\h}{2}  < \h \ .\end{multline*}
From this and \eqref{eta},  all three estimates in \eqref{desired} follow.
\end{pf} \par
We are now ready to conclude the proof following the arguments described in (a) and (b) above. 
Using Lemmas \ref{lemma42-1},  \ref{lemma42} and \ref{lemma43},  we  may consider the  families of control curves and functions, parameterised  by a  positive  $\d$ tending to  $0$, 
$$v^{\d}(t)\= v^{(\d, \d)}(t)\ ,\qquad  f^{\d}(t, x^\ell_{(r)}, u) \=   f^{(\d, \d, \d)}(t, x^\ell_{(r)}, u)\ .$$
They  have the following properties:  
\begin{itemize}[leftmargin = 15 pt]
\item each  map $v^{\d}: [0, T] \to  \wh K'$ is  polynomial, it  satisfies  $|v^{\d}(\t_o) - u_o(\t_o)| < \d$ and the corresponding solution $x^{(v^{\d}, \sfs_o)}(t)$ to the equations \eqref{diffe-1-ter}  satisfies 
\beq \label{4.15-ter}
\|x^{(u_o, \mathsf s_o)}- x^{(v^{\d}, \mathsf s_o)}\|_{\cC^{2k-2}}  < \d\ ;\eeq
\item $f^{\d}: \Omega \times \wh K'' \to \bR^n$ is a smooth  function  and  the solution  $x^{( v^{\d}, \mathsf s_o; \d)} (t)$  to  the differential problem
\beq\label{diffe-hat-bisbis} x^i_{(k)}(t) =   f^{\d\, i}\left(t, j^{k-1}_t(x), v^{\d}(t)\right)\ ,\qquad j^{k-1}_{t=0}(x) = \sfs_o
\eeq
 satisfies
\beq \label{desired-bis} 
\begin{split}  & \| x^{( v^{\d}, \sfs_o; \d)}(t) - x^{(u_o, \sfs_o)}(t)\|_{\cC^{2k-2}} < \d\hskip 4 cm  \text{and}\\
& \left|f^{\d}(t, j^{k-1}_t(x^{( v^{\d},\sfs_o; \d)}), v^\d(\t_o))  -  f(t,j^{k-1}_t( x^{(u_o, \mathsf s_o)}), u_o(\t_o)) \right| < \d \ .
\end{split}
\eeq
\end{itemize}
Therefore,  for any sufficiently small $\d, \ve> 0$ we may also consider:  
\begin{itemize}[leftmargin = 15pt]
\item  the  real function  $\sH^{\d}: K \to \bR$ which is defined by 
\beq  \sH^{\d}(\o) \= p^{( v^{\d},\sfs_o; \d)}_i(\t_o) f^{\d\,i}(\t_o, j^{k-1}_{\t_o}(x^{(v^{\d}, \sfs_o; \d)}), \o)\ ,\eeq
where we denote by $p^{( v^{\d}, \sfs_o; \d)}(t)$ the solution to \eqref{controlledEL-1bis*} with  $f$ replaced 
by the smooth  $f^{\d}$ and  determined by the control curve $v^{\d}$ and the initial value $\sfs_o$; 
\item  the needle modification  $v^{\d, \ve} \= v^\d{}^{(\t_o, \o_o, \ve)}$ of the polynomial  curve $v^\d(t)$,  with peak time $\t_o$,  ceiling value $\o_o$ and width $\ve$; 
\item  the needle modification $u^{\ve}_o \= u^{(\t_o, \o_o, \ve)}_o$ of the (merely piecewise  $\cC^{k-1}$)  $u_o(t)$, also  with   peak time $\t_o$,  ceiling value $\o_o$ and width $\ve$.
\end{itemize}
By the  Lemmas  \ref{lemma42-1},  \ref{lemma42}  and \ref{lemma43},  for $\d \to 0$ the functions $f^\d$,  the curves  $ p^{(v^\d,\sfs_o; \d)}(t)$ and the curves of jets $j^{k-1}_t(x^{( v^\d, \sfs_o; \d)})$    tend uniformly on compacta to the map $f$, to  the curve  $p^{( u_o,\sfs_o)}(t)$ and to the curve of jets  $j^{k-1}_t(x^{(u_o, \sfs_o)})$, respectively. 
 Therefore, if we set 
 \beq \label{kappadelta} \k_o^\d \= \sH^{\d}(\o_o)  - \sH^{\d}(v^\d(\t_o))   \eeq
 we have that $\lim_{\d \to 0} \k^\d_o = \k_o   > 0$ and hence there is $\d_o > 0$ such that 
$\k_o^\d   > 0$  for any   $\d \in (0, \d_o]$.\par
\smallskip
 We now observe that  the (restrictions to an appropriate relatively compact neighbourhood of the $k$-jets  of the curve $\g_o(t) = (t, x^{(u_o, \sfs_o)}(t),  p^{(u_o, \sfs_o)}(t))$    of the)  Lagrangians 
\beq \label{finallagrangians} L^{(\d)}(t, x^\ell, \ldots, x^\ell_{(k)}, p_j, u) \= p_j \big(x^j_{(k)} - f^{\d\, j}(t, x^\ell, \ldots , x^\ell_{(k-1)}, u)\big)\eeq 
and the one-parameter families of control pairs $U^{(\d)} = (v^{\d}, (\sfs_o, \wt\sfs_o^{(v^{\d}, \sfs_o)}))$ are such that the conditions (A'), (B) and (D) of  Remark \ref{remarketto} 
are satisfied with $\s^{(0)}_o \= (\sfs_o, \wt \sfs_o^{(u_o, \sfs_o)})$,  $u^{(0)}_o\= u_o$ and  that the limit Lagrangian for $\d \to 0$ is
$$L^{(0)}(t, x^\ell, \ldots, x^\ell_{(k)}, p_j, u) \= p_j \big(x^j_{(k)} - f^{j}(t, x^\ell, \ldots , x^\ell_{(k-1)}, u)\big)\ .$$
Moreover, since  each  Lagrangian \eqref{finallagrangians} is smooth, for each $\d$ we may consider  a   one-parameter family of  {\it good} needle variations $ {\needle}^{(\t_o, \o_o,  \Sigma^{(\d)},  \ve_o)} (\g^\d_o)$  for the pair $(\t_o, \o_o)$   as 
defined in  Lemma \ref{existenceofgood} . By definition of the good needle variations,  the family of   initial data maps  $\Sigma^{(\d)}$ converges uniformly to a limit initial data map $\Sigma^{(0)}$ and satisfies conditions (C) and (F) of Remark \ref{remarketto}. We may also consider an appropriate family of smooth cost functions $C^\d$, which depend on the $(k-1)$-jets at $t = T$ of the controlled curves and converge in the  $\cC^1$ norm to the cost function $C = C^{(0)}$ on an appropriate relatively compact neighbourhood of the end $k-1$-jets of the curves of the needle variations determined by $ {\needle}^{(\t_o, \o_o,  \Sigma^{(\d)},  \ve_o)} (\g^\d_o)$. This would imply that also condition (G) is satisfied. \par
We  finally observe that the  functions $V^{(\d)}$, associated with the above described  Lagrangians, control pairs and initial valued maps,  satisfy also condition (E) of  Remark \ref{remarketto}. This is in fact a direct consequence of the property that, by construction of the polynomial curves $v^{\d}(t)  = v^{(\d, \d)}(t)$,   the    $v^\d(t)$ converge in $\cC^{k-1}$ norm to the function $u_o(t)$  on a fixed interval containing  $\t_o$. Since all conditions of Remark \ref{remarketto} are satisfied, we  infer   the existence of a constant $M > 0$ and a needle modification $u'(t)$   of $u_o(t)$,  with  peak time $\t_o$ and ceiling $\o_o$,  such that  \eqref{eccolo} holds. This concludes the proof with  $k\geq 2$.\par
\medskip
Let us now focus on the case $k = 1$,   i.e. on the situation in  which Theorem \ref{main1} reduces to the   classical PMP.  As the reader will shortly see,   the  approximation procedure  which  we used for  $k > 1$ is  valid also for $k = 1$,  but requires  some nontrivial adjustments. We give here  such adjustments in full detail  mainly with the following purposes:   (a)   showing that, in the classical setting,  our   approximation technique  has  the same power of the standard approach;  (b) paving the way  for  future  developments in different   contexts -- see \S \ref{added}. 
  \par
   For the time being, we  assume  that $f(t, x, u)$ is  continuously differentiable not only with respect to $x$ but also  with respect to $t$ (we  show how to remove this assumption later). As announced above, for each sufficiently small $\d_0$ let us denote by $u^{(\d_0)}: [0, T] \to \wh K$ the control curve defined by 
\beq  u^{(\d_0)}(t) = \left\{ \begin{array}{cc} u_o(\t_o) & \text{if}\ t \in [\t_o - \d_0, \t_o]\ ,\\[10pt]
u_o(t) & \text{otherwise}\ .
\end{array} \right.
\eeq
Note that $\dist(u_o, u^{(\d_0)}) \leq \d_0$. Hence, by the usual circle of ideas,  the corresponding control curve $x^{(u^{(\d_0)}, \sfs_o)}(t)$ and the associated curve $p^{(u^{(\d_0)}, \sfs_o)}(t)$ uniformly converge to the curves $x^{(u_o, \sfs_o)}(t)$ and  $p^{(u_o, \sfs_o)}(t)$, respectively, for $\d_0 \to 0$. So, by continuity of $f$ and $C$,  if we set 
\begin{multline} \k^{(\d_0)} \= \\
= p^{(u^{(\d_0)}, \sfs_o)}_i(\t_o) f^i(\t_o, x^{(u^{(\d_0)}, \sfs_o)}(\t_o), \o_o) -  p^{(u^{(\d_0)}, \sfs_o)}_i(\t_o) f^i(\t_o, x^{(u^{(\d_0)}, \sfs_o)}(\t_o), u^{(\d_0)}(\t_o))  = \\
=  p^{(u^{(\d_0)}, \sfs_o)}_i(\t_o) f^i(\t_o, x^{(u^{(\d_0)}, \sfs_o)}(\t_o), \o_o)  -  p^{(u^{(\d_0)}, \sfs_o)}_i(\t_o) f^i(\t_o, x^{(u^{(\d_0)}, \sfs_o)}(\t_o), u_o(\t_o))\ ,\end{multline} 
we directly obtain that 
$\lim_{\d_0 \to 0}  \k^{(\d_0)} = \k_o$ and  $\lim_{\d \to 0} \sC(x^{(u^{(\d_0)}, \sfs_o)}(T)) =    \sC(x^{(u_o, \sfs_o)}(T))$. 
 \par
\medskip
Now,  let us consider  the following analog of Lemma \ref{lemma42-1}.
\begin{lem} \label{lemma42-1bis} Given the curve $u^{(\d_0)}$, for any $\h> 0$, there exists a continuous  curve  $v^{(\d_0, \h)}: [0, T] \to \wh K$, which coincides with $u^{(\d_0)}(t)$ on
the interval $[\t_o - \d_0, \t_o]$ and  such that 
\begin{multline}\label{firststep-1bis}
  \dist(u^{(\d_0)}, v^{(\d_0, \h)})  < \h\ ,\\
   \  \|x^{(u^{(\d_0)}, \mathsf s_o)} - x^{(v^{(\d_0, \h)}, \mathsf s_o)}\|_{\cC^{0}} \ , \|p^{(u^{(\d_0)}, \mathsf s_o)} - p^{(v^{(\d_0, \h)}, \mathsf s_o)}\|_{\cC^{0}}< \h
\ .\end{multline} 
\end{lem}
\begin{pf} We recall that, by the Lusin Theorem (see e.g. \cite[p.14]{Ad}),  for    any choice of $\wt \d > 0$ there exists  a  $v^{(\wt \d)} \in \cC^0([0,T], \bR^m)$ such that  
\beq \label{Lusin} \sup_{t \in [0,T]} |v^{(\wt \d)}(t)| \leq \sup_{t \in [0,T]} |u_o(t)|  \qquad \text{and}\qquad  \dist(u_o^{(\d_0)}, v^{( \wt \d)}) < \wt \d\ .\eeq
Note that such a curve can be taken  equal to   $u^{(\d_0)}(t)$ on the interval
 $[\t_o - \d_0, \t_o]$. Indeed, one can  start with  a continuous  curve  $\wt v^{(\wt \d)}(t)$   satisfying   the first inequality in \eqref{Lusin} and with  
$\dist( u^{(\d_0)}_o,  \wt v^{(\wt \d)}) < \frac{\wt \d}{2}$.  Then one can modify $\wt v^{(\wt \d)}(t)$ just  on  the closed interval $[\t_o  - \d_0 - \frac{\wt \d}{4}, \t_o + \frac{\wt \d}{4} ]$,  determining a  continuous curve $v^{(\wt \d)}$  that satisfies also   the desired additional requirement. Now, given $\h$, by taking $\wt \d$ sufficiently small, the curve $v^{(\d_0, \h)}(t)  \= v^{(\wt \d)}(t)$ satisfies all required inequalities.
 \end{pf}
  \par
  The same arguments of Lemma \ref{lemma42} imply that, for  any  continuous control curve $v^{(\d_0, \d_1)}(t)$  and any  $\h > 0$
 there exists  a polynomial curve $v^{(\d_0, \d_1, \h)}(t)$ in $\wh K' \supset \wh K $   with 
 \beq \label{att_o-bis}|v^{(\d_0, \d_1, \h)}(\t_o) - v^{(\d_0, \d_1)}(\t_o)|  = |v^{(\d_0, \d_1, \h)}(\t_o) - u_o(\t_o)| < \h\eeq
 and  such that the   solutions $x^{(v^{(\d_0, \d_1)}, \sfs_o)}$ and $x^{(v^{(\d_0, \d_1, \h)}, \sfs_o)}$  and the associated curves $p^{(u^{(\d_0, \d_1)}, \sfs_o)}(t)$  and $p^{(v^{(\d_0, \d_1, \h)}, \sfs_o)}(t)$   satisfy 
\beq \label{4.15-bis}
\|x^{(v^{(\d_0, \d_1)}, \mathsf s_o)}- x^{(v^{(\d_0, \d_1, \h)}, \mathsf s_o)}\|_{\cC^{0}} \ , \quad \|p^{(v^{(\d_0, \d_1)}, \mathsf s_o)}(t) - p^{(v^{(\d_0, \d_1, \h)}, \mathsf s_o)}(t)\|_{\cC^{0}} <  \h\ .
\eeq
This can be considered as  the analog of  (II) for  $k= 1$ and with  control curve  given by $u^{(\d_0)}(t)$ in place of $u_o(t)$. It  is also quite immediate to check  that, for any given polynomial control curve  $v^{(\d_0, \d_1, \d_2)}(t)$ and for any $\h> 0$, there exists  a smooth function $ f^{(\d_0, \d_1, \d_2, \h)}:  \O  \times \wh K'' \to \bR^n$, with  the properties given in Lemma \ref{lemma43} for $k = 1$.  In other words,  the analog of (III)   is also true.  \par
\smallskip
At this point,  if  for any $\d  = \d_{0} > 0$  we  define  
\begin{multline} \label{familyk=1} u_o^\d(t) \= u^{(\d_0)}(t)\ ,\ \  \wh u^\d_o(t) \= v^{(\d, \d)}(t)\ ,\\
  v^\d(t)\= v^{(\d, \d, \d)}(t)\ ,\ \   f^\d(t, x, u) \=   f^{(\d, \d,\d, \d)}(t, x, u)\ ,
  \end{multline}
then, for all  $\d$ sufficiently small, we have:  
\begin{itemize}[leftmargin = 15 pt]
\item each $u^\d_o: [0, T] \to   K$  is a measurable curve which coincides with $u_o(t)$ outside of  $[\t_o - \d, \t_o]$ and it is constant and equal to $u_o(\t_o)$ on such interval; 
\item each $\wh u^\d_o: [0, T] \to  \wh K$ is a continuous curve  which is constant  equal to $u_o(\t_o)$  on $[\t_o - \d, \t_o]$ and with $\dist(\wh u_o^\d, u_o^\d) < \d$; 
\item each   $v^\d: [0, T] \to  \wh K'$ is a polynomial  with   $|v^\d(\t_o) - \wh u^\d_o(\t_o)|  = |v^\d(\t_o) -  u_o(\t_o)| < \d$; 
\item $f^\d: \Omega \times \wh K'' \to \bR^n$ is a smooth  function  and  the solution  $x^{(v^\d, \mathsf s_o; \d)} (t)$  to  
\beq\label{diffe-hat-bis} x^i_{(1)}(t) =   f^{\d\, i}\left(t, x(t), v^\d(t)\right)\ ,\qquad x(0) = \sfs_o
\eeq
 satisfies
\beq \label{desired-bisbis} 
\begin{split}  & \| x^{(v^\d, \sfs_o; \d)}(t) - x^{(\wh u^\d_o, \sfs_o)}(t)\|_{\cC^{0}} ,  \| x^{(v^\d, \sfs_o; \d)}(t) - x^{(u_o, \sfs_o)}(t)\|_{\cC^{0}} < \d\ ,\\
& \left|f^\d(t, x^{(v^\d,\sfs_o; \d)}(t), v^\d(t))  -  f(t, x^{(\wh u^\d_o, \mathsf s_o)}(t), \wh u^\d_o(t)) \right| < \d \ ,\quad   \text{for any}\ t \in [0,T]\ ,\\
&   \left|f^\d(t, x^{(v^\d,\sfs_o; \d)}(t), \o)  -  f(t, x^{(u_o, \mathsf s_o)}(t), \o) \right| < \d\ ,\qquad \text{for any}\ \o \in K\ .
\end{split}
\eeq
\end{itemize}
Therefore,  for any sufficiently small $\d , \ve> 0$ we may  consider:  
\begin{itemize}
\item  the  real function  $\sH^{\d}: K \to \bR$ defined by 
\beq  \sH^{\d}(\o) \= p^{(v^\d,\sfs_o; \d)}_i(\t_o) f^{\d\,i}(\t_o, x^{(v^\d, \sfs_o; \d)}(\t_o), \o)\ ,\eeq
where, as usual,  we denote by $p^{( v^\d, \sfs_o; \d)}(t)$ the solution to \eqref{controlledEL-1bis*} with  $f$ replaced 
by the smooth  $f^\d$ and  determined by the control curve $v^\d$ and the initial value $\sfs_o$; 
\item  the needle modification  $v^{\d, \ve}(t) \= v^\d{}^{(\t_o, \o_o, \ve)}(t) $ of the polynomial  curve $v^\d(t)$,  with peak time $\t_o$,  ceiling value $\o_o$ and width $\ve$; 
\item  the needle modification $ u^{\d, \ve}_o(t) \=  u^{\d (\t_o, \o_o, \ve)}_o(t)$ of  $ u^\d_o(t)$ with   peak time $\t_o$,  ceiling value $\o_o$ and width $\ve$;
\item  the needle modification $u^{\ve}_o(t) \= u^{(\t_o, \o_o, \ve)}_o(t)$ of  $u_o(t)$ with   peak time $\t_o$,  ceiling value $\o_o$ and width $\ve$.
\end{itemize}
{\it Note that for $\d < \ve$, the needle modifications  $ u^{\d, \ve}_o(t)$ and  $u^{\ve}_o(t)$ coincide.}\par
\smallskip
By the  Lemmas \ref{normallemma},  \ref{lemma42-1bis} and the above remarks,  for $\d \to 0$ the functions $f^\d$,  the curves  $x^{(v^\d, \sfs_o; \d)}(t)$ and  $ p^{(v^\d,\sfs_o; \d)}(t)$      tend uniformly on compacta to the map $f$ and to  the curves  $x^{(u_o, \sfs_o)}(t)$ and  $p^{( u_o,\sfs_o)}(t)$, respectively. 
 Therefore, setting $\k_o^\d$ as in \eqref{kappadelta},  
 we have that $\lim_{\d \to 0} \k^\d_o = \k_o   > 0$ and hence there is $\d_o > 0$ such that 
$\k_o^\d   > 0$  for any   $\d \in (0, \d_o]$.\par
\smallskip
As for the previous case, it is straightforward  to check that  the conditions (A'), (B), (C), (D)  and (F) of  Remark \ref{remarketto} are satisfied also for the new $\d$-parameterised family  of control curves and Lagrangians $L^{(\d)}$, determined by the functions $f^\d$.  As before, we can also consider an appropriate family of smooth cost functions $C^{(\d)}$ which converge in the $\cC^1$ norm to $C = C^{(0)}$ on an appropriate relatively compact neighbourhood of the endpoints of the curves of the needle variations. So,  in order to conclude, it remains to prove that also condition (E) is satisfied, i.e. that the functions $V^{(\d)}(\ve)$ tend to some function, which is  an infinitesimal for $\ve \to 0$.  In fact, we claim that  the  $V^{(\d)}(\ve)$  tend uniformly on compacta of $(0, \ve_o]$ to the function $V_o(\ve)$ defined by 
\beq V_o(\ve) = \frac{1}{\ve} \int_{\t_o  -\ve} ^{\t_o} \bigg(p^{(u_o,\sfs_o)}_i(t) f^{i}(t, x^{(u_o, \sfs_o)}(t), u_o(\t_o)) - p^{(u_o,\sfs_o)}_i(t)  f^{i}(t, x^{(u_o, \sfs_o)}(t), u_o(t)) \bigg) dt\eeq
and that such a function is an infinitesimal for $\ve \to 0$, as required.\par
To prove  this,  on a fixed  interval $[\r, \ve_o]$, $\r > 0$, we need to show that 
\begin{multline} \label{theclue}  |V^{(\d)}(\ve) - V_o(\ve)| = \\
=\frac{1}{\ve} \bigg|\int_{\t_o - \ve}^{\t_o} \bigg\{ p^{(v^\d,\sfs_o; \d)}_i(t)  f^{\d\,i}(t, x^{(v^\d, \sfs_o; \d)}(t), v^\d(\t_o)) - p^{(v^\d,\sfs_o; \d)}_i(t)  f^{\d\,i}(t, x^{(v^\d, \sfs_o; \d)}(t), v^\d(t))  - \\
- p^{(u_o,\sfs_o)}_i(t) f^{i}(t, x^{(u_o, \sfs_o)}(t), u_o(\t_o)) + p^{(u_o,\sfs_o)}_i(t)  f^{i}(t, x^{(u_o, \sfs_o)}(t), u_o(t)) \bigg\} dt \bigg|
\end{multline}
is uniformly bounded by some constant $\gC_\d$ tending to $0$ for $\d \to 0$. In this regard, we recall that: 
\begin{itemize}[leftmargin = 15pt]
\item $\wh u^\d_o(\t_o) = u^\d_o(\t_o) =  u_o(\t_o)$ for any $\d$ and $v^{\d}(\t_o)$ tends to $u_o(\t_o)$ for $\d \to 0$; 
\item the maps  $x^{(v^\d, \sfs_o; \d)}(t)$,    $p^{(v^\d,\sfs_o; \d)}(t)$  and $f^{\d}(t, x, u)$ converge   to $x^{(u_o, \sfs_o)}(t)$,  $p^{(u_o,\sfs_o)}(t)$ and $f(t, x, u)$, respectively, in the $\cC^0$ norm. 
\end{itemize}
Due to this,  it suffices  to check  that the function $\D V^{(\d)}: [\r,\ve_o] \to \bR$ defined by 
\begin{multline}
\D V^{(\d)}(\ve) {\=}\frac{1}{\ve} \bigg|\int_{\t_o - \ve}^{\t_o} \bigg\{ p^{(u_o,\sfs_o)}_i(t)  f^{i}(t, x^{(u_o, \sfs_o)}(t), u_o(\t_o)) - p^{(u_o,\sfs_o)}_i(t)  f^{i}(t, x^{(u_o, \sfs_o)}(t), \wh u_o^\d(t))  - \\
- p^{(u_o,\sfs_o)}_i(t) f^{i}(t, x^{(u_o, \sfs_o)}(t), u_o(\t_o)) + p^{(u_o,\sfs_o)}_i(t)  f^{i}(t, x^{(u_o, \sfs_o)}(t), u_o(t)) \bigg\} dt \bigg| = \\
= \frac{1}{\ve} \bigg|\int_{\t_o - \ve}^{\t_o} \bigg\{   p^{(u_o,\sfs_o)}_i(t)  f^{i}(t, x^{(u_o, \sfs_o)}(t), \wh u_o^\d(t))  
 - p^{(u_o,\sfs_o)}_i(t)  f^{i}(t, x^{(u_o, \sfs_o)}(t), u_o(t)) \bigg\} dt \bigg| 
\end{multline}
is such that $\sup_{\ve \in[\r, \ve_o] } \D V^{(\d)}(\ve)$ tends to $0$ for $\d \to 0$.  In order to prove this, we first observe that, since  $ \wh u_o^\d(t)$  differs from   $u_o^\d(t)$ just on a set of measure  less than $\d$, we may replace $ \wh u_o^\d(t)$ by     $u_o^\d(t)$  and  get the inequality 
\begin{multline} \D V^{(\d)}(\ve) \leq  \frac{1}{\ve} \bigg|\int_{\t_o - \ve}^{\t_o} \bigg\{   p^{(u_o,\sfs_o)}_i(t)  f^{i}(t, x^{(u_o, \sfs_o)}(t),  u_o^\d(t))  
 -\\ -  p^{(u_o,\sfs_o)}_i(t)  f^{i}(t, x^{(u_o, \sfs_o)}(t), u_o(t)) \bigg\} dt \bigg|  + 
 \text{Const}\ \d
\end{multline}
for an appropriate constant  $\text{Const}$. Secondly, 
we recall that 
$\big\{\ u_o^{\d}(t) \neq u_o(t)\ \big\} = [ \t_o - \d, \t_o]$
and that $ u_o^\d|_{ [ \t_o - \d, \t_o]} \equiv u_o(\t_o)$. Hence, 
if we set $\wt{\text{Const}}\= \max_{t \in [0, T]} |p^{(u_o,\sfs_o)}(t)|$ and take $\d < \r \leq \ve$,  we get  
\begin{multline*} \D V^{(\d)}(\ve) \leq \frac{1}{\ve} \int_{\t_o - \ve}^{\t_o} \bigg| p^{(u_o,\sfs_o)}_i(t) \bigg(  f^{i}(t, x^{(u_o, \sfs_o)}(t),  u_o^\d(t))  
 -  f^{i}(t, x^{(u_o, \sfs_o)}(t), u_o(t)) \bigg) \bigg|dt +\\
+ \text{Const}\ \d  \leq
\end{multline*}
\begin{multline*}
 \leq \frac{\wt{\text{Const}}}{\ve} \int_{\t_o - \ve}^{\t_o} \bigg| f(t, x^{(u_o, \sfs_o)}(t),  u_o^\d(t))  
 -  f(t, x^{(u_o, \sfs_o)}(t), u_o(t)) \bigg|dt  + \\
 + \text{Const}\ \d  =\\
= \frac{\wt{\text{Const}} }{\ve}\int_{\t_o - \d}^{\t_o} \bigg| f(t, x^{(u_o, \sfs_o)}(t), u_o(\t_o))  
 -  f(t, x^{(u_o, \sfs_o)}(t), u_o(t))  \bigg|dt +   \text{Const}\,\d \leq \\
 \leq \frac{\wt{\text{Const} }}{\ve}\int_{\t_o - \d}^{\t_o} \bigg| f(t, x^{(u_o, \sfs_o)}(t), u_o(\t_o))   -   f(\t_o, x^{(u_o, \sfs_o)}(\t_o), u_o(\t_o)) \bigg|dt + \\
 + \frac{\wt{\text{Const}} }{\ve}\int_{\t_o - \d}^{\t_o} \bigg| f(\t_o, x^{(u_o, \sfs_o)}(\t_o), u_o(\t_o)) 
 -  f(t, x^{(u_o, \sfs_o)}(t), u_o(t))  \bigg|dt  + \\
 +   \text{Const}\,\d\ .
 \end{multline*}
It follows that   
\begin{multline*}
\D V^{(\d)}(\ve)  \leq \frac{\wt{\text{Const}} }{\d}  \int_{\t_o - \d}^{\t_o} \bigg| f(t, x^{(u_o, \sfs_o)}(t), u_o(\t_o))   -   f(\t_o, x^{(u_o, \sfs_o)}(\t_o), u_o(\t_o)) \bigg|dt + \\
 + \frac{\wt{\text{Const}}}{\d}\int_{\t_o - \d}^{\t_o} \bigg| f(\t_o, x^{(u_o, \sfs_o)}(\t_o), u_o(\t_o)) 
 -  f(t, x^{(u_o, \sfs_o)}(t), u_o(t))  \bigg|dt + \\
  +   \text{Const}\,\d  
 \end{multline*}
and the right hand side  can be assumed to be smaller than any desired quantity because of 
 the continuity of $f$ and $x^{(u_o, \sfs_o)}(t)$ and the condition  \eqref{BPcondition}. \par
 \smallskip
This proves that for any choice of  $\r \in (0, \ve_o)$,  the restriction $\D V^{(\d)}|_{[\r, \ve_o]}$ converges uniformly to $0$. Note also that, by a similar  argument,  the assumption \eqref{BPcondition} and the continuity of  $ p^{(u_o,\sfs_o)}(t)$, $ x^{(u_o,\sfs_o)}(t)$ and $f(t, x, u)$ imply that  the function 
\begin{multline} 
V_o(\ve) \= 
\frac{1}{\ve} \int_{\t_o - \ve}^{\t_o} \bigg\{   p^{(u_o,\sfs_o)}_i(t) f^{i}(t, x^{(u_o, \sfs_o)}(t), u_o(\t_o)) - p^{(u_o,\sfs_o)}_i(t)  f^{i}(t, x^{(u_o, \sfs_o)}(t), u_o(t)) \bigg\} dt
\end{multline}
is an infinitesimal for $\ve \to 0$. This concludes the proof that also condition (E) is satisfied.  By Remark \ref{remarketto}, we  conclude  that the theorem holds also  in case  $k = 1$ 
and under the additional assumption that  $f$ is continuously differentiable with respect to $t$.\par
In order to conclude, it is now necessary to remove this  assumption. 
This can be done by noting that the term  \eqref{41212} considered in the proof of Lemma \ref{corolletto}, can be  written as  
\beq  \int_{\t_o - \ve}^{\t_o } \left(\cP^{(\s_{(\ve)}, \widecheck u_{(\ve)}, t)}(\o_o) -  \cP^{(\s_o,  u_o, t)}(u_o(t)) \right) dt  = \ve \k_o  + \ve W(\ve)\eeq
 with 
 \begin{multline} \label{theW} W(\ve) \= \frac{1}{\ve}  \int_{\t_o - \ve}^{\t_o } \left(\cP^{(\s_{(\ve)}, \widecheck u_{(\ve)}, t)}(\o_o) - \cP^{(\s_o,  u_o, \t_o)}(\o_o)  \right) - \\
 -  \left(  \cP^{(\s_o,  u_o, t)}(u_o(t))-  \cP^{(\s_o,  u_o, \t_o)}(u_o(\t_o)) \right) dt\ .
 \end{multline} 
Considering this  new function  instead  of the function $V(\ve)$,  the  family of control curves and differential constraints defined  in \eqref{familyk=1} 
satisfies the following analogs of the   \eqref{manyineq}, which  involve the  $\d$-parametrised family of new functions $W^{(\d)}(\ve)$ instead of the $V^{(\d)}(\ve)$:
 \beq\label{manyineq-bis}  C_\d^{(\ve)}\leq  C^{(0)}_\d  -    \ve(\k^{(\d)}_o -   \wt \gd^{(\d)} \ve - |W^{(\d)}(\ve)| )\eeq
 where  $\wt \gd^{(\d)}$ is now a constant which  {\it does not depend} on  the derivative $\frac{\p f^{\d}}{\p t}$.   The explicit expression of   the function $W^{(\d)}(\ve)$ 
can be directly derived from \eqref{theW}. One finds
\begin{multline} W^{(\d)}(\ve) \=  \\
 \=  \frac{1}{\ve} \int_{\t_o  -\ve} ^{\t_o} \bigg\{p^{(\widecheck v^\d_{(\ve)},\sfs_{(\ve)}; \d)}_i(t) \bigg( x^{(\widecheck v^\d_{(\ve)} , \sfs_{(\ve)}; \d)i }_{(1)}(t) - f^{\d\,i}(t, x^{(\widecheck v^{\d}_{(\ve)}, \sfs_{(\ve)}; \d))}(t), \o_o)\bigg) - \\
 - p^{(v^\d,\sfs_o; \d)}_i(\t_o) \bigg( x^{(v^{\d} , \sfs_o; \d)i }_{(1)}(\t_o) - f^{\d\,i}(\t_o, x^{(v^\d, \sfs_o; \d))}(\t_o), \o_o)\bigg) \bigg\} - \\
-\bigg\{- p^{(v^\d,\sfs_o; \d)}_i(t)  f^{\d,i}(t, x^{(v^{\d}, \sfs_o; \d))}(t), v^\d(t)) + p^{(v^\d,\sfs_o; \d)}_i(\t_o)  f^{\d,i}(\t_o, x^{(v^{\d}, \sfs_o; \d))}(\t_o), v^\d(\t_o)) \bigg\} dt\ , 
 \end{multline} 
 where, according to the  notational conventions of Lemma \ref{corolletto},   $\widecheck v^\d_{(\ve)}$ is the smoothed needle variation of the polynomial  curve $v^\d(t)$ of width $\ve$ and peak time $\t_o$, while  $\sfs_{(\ve)}$ is the corresponding initial datum for the curve $\g(t) = (x(t), p(t))$. This initial datum  $\sfs_{(\ve)}$ is determined so that  the initial  value  for $x(t)$ is $\sfs_o$, while  the initial value for $p(t)$ is  prescribed in order to have the usual terminal conditions at $t = T$.  We also assume that the constant $\gh$, used in the definition of the smoothed needle modifications (see \eqref{theh}) is chosen differently for each value of  $\d$ and in a way that  $\gh = \gh(\d)$ tends to $0$ for $\d \to 0$. 
 In order to conclude, it is now sufficient to show (in analogy with what  we did above for the functions $V^{(\d)}(\ve)$) that the  function $W^{(\d)}(\ve)$ 
converge uniformly on compacta on $(0, \ve_o]$ to a function $W^{(0)})(\ve)$ and that the limit function  $W^{(0)}(\ve)$  is an infinitesimal for $\ve \to 0$. This can be checked directly. More precisely, following the same circle of ideas 
as above, one can see that on any interval $[\r, \ve_o]$, $\r > 0$ the function  $W^{(\d)}$ converges in $\cC^0$ norm   to the function
\begin{multline} W^{(0)}(\ve) \=    \frac{1}{\ve} \int_{\t_o  -\ve} ^{\t_o} \bigg\{p^{(u_{(\ve)},\sfs_{(\ve)})}_i(t) \bigg( x^{(u_{(\ve)} , \sfs_{(\ve)})i }_{(1)}(t) - f^{i}(t, x^{(u_{(\ve)}, \sfs_{(\ve)}))}(t), \o_o)\bigg) - \\
 -p^{(u_o,\sfs_o)}_i(\t_o) \bigg( x^{(u_o , \sfs_o)i }_{(1)}(\t_o) - f^{i}(\t_o, x^{(u_o, \sfs_o))}(\t_o), \o_o)\bigg) \bigg\} - \\
-\bigg\{- p^{(u_o,\sfs_o)}_i(t)  f^{i}(t, x^{(u_o, \sfs_o; \d))}(t), u_o(t))   + p^{(u_o,\sfs_o)}_i(\t_o)  f^{i}(\t_o, x^{(u_o, \sfs_o))}(\t_o), u_o(\t_o)) \bigg\} dt 
 \end{multline} 
and that such a function is an infinitesimal.  For brevity, we omit the details.  
\end{pf}
 \par
   \section{The proofs of  Sublemma \ref{lemma37} and Lemma \ref{corolletto-bis}}  \label{theremaining}
  \subsection{The proof of Sublemma \ref{lemma37}}
  \label{sublemmaproof} Let us denote by  $\o^i_{(\d)}$     the $1$-forms     on  the manifold $J^{2k+1}(\cQ|\bR) \times \bR^M$ defined by 
\beq \label{holonomic}
\o^i_{(\d)} \= dq^i_{(\d)} -  q^i_{(\d+1)} dt\ ,\qquad
  \d = 0, \ldots, 2k\ .
\eeq
Using these $1$-forms, we can introduce  the  {\it controlled Poincar\'e-Cartan form}   $\b^{PC}$  associated with the controlled Lagrangian $L$  (\cite[Sect. 5]{CGS})
\beq \label{310} \b^{PC} \=    L dt  + \sum_{\d = 1}^{k} \sum_{\h = 0}^{\d-1} (-1)^{\h}  \frac{d^\h}{dt^\h}  \left(\frac{\p  L  }{\p q^i_{(\d)}} \right) \o^i_{(\d-(\h+1))}  \ .
\eeq
By  basic facts on    variationally equivalent  $1$-forms (see e.g. \cite[Prop. A2]{Sp} and \cite[Proof of Lemma 5.2]{CGS}), the exterior  differential    $d \b^{PC}$   has the form  
\beq \label{39} \begin{split}
 d  \b^{PC} {=} & E(L)_i  \o^i_{(0)} \wedge dt 
 +  \frac{\p  L}{\p u^a} d u^a \wedge dt  +\\
& +   \text{\rm  linear combinations of   wedges  of pairs of   
  $1$-forms of the kind}\ \eqref{holonomic}\ ,
\end{split}\eeq
where the functions $E(L)_i$ are the controlled Euler-Lagrange expressions defined in \eqref{controlledEL}. 
Consider   the (smooth) map
\beq \label{350} \cG = \cG_{(\ve)}: [0, T] \times [0,1] \to J^{2k+1}( \cQ|\bR) \times \bR^M\ ,\quad \cG(t,s) \= \left(j^{2k+1}_{t}({\widecheck \g}{}^{(\ve, s)}), u^{(\ve,s)}(t)\right)\ .\eeq
and  the  fields  of tangent vectors of the surface $\cS \= \cG([0, T] \times [0,1])$ defined by  
\beq \label{4.20} X|_{\cG(t,s)} \= \frac{\p \cG}{\p t}\bigg|_{(t,s)} = \cG_*\left(\frac{\p}{\p t}\bigg|_{(t,s)}\right)\ ,\qquad Y|_{\cG(t,s)} = \frac{\p \cG}{\p s}\bigg|_{(t,s)}  = \cG_*\left(\frac{\p}{\p s}\bigg|_{(t,s)}\right)\ .\eeq
By construction,  each vector  $X|_{\cG(t,s)}$ has the first components that are  tangent to the curve of jets, determined by a $ \cK$-controlled curve ${\widecheck \g}^{(\ve, s)}(t) = (t, q^{(\ve,s)}(t))$. In particular, the  $\frac{\p}{\p t}$-component  of  $X|_{\cG(t,s)}$   is  $1$ for any $(t,s)$.  By the same reason   the  $\frac{\p}{\p t}$-component of $Y|_{\cG(t,s)}$ is identically $0$. 
Due to this and the vanishing  of the $1$-forms \eqref{holonomic} on the tangent vectors of curves of  jets given  by  curves in $\cQ \times \bR$, we have that  for any $(t, s) \in [0, T] \times [0,1]$
\beq
\begin{split} \label{348} &\b^{PC}(X|_{\cG(t,s)}) =      L\big|_{j^{k}_t(\cG(\cdot,s))} \ ,  \\
&d \b^{PC}(X|_{\cG(t,s)}, Y|_{\cG(t,s)})  = -  \frac{\p L}{\p u^a}\bigg|_{j^{k}_t(\cG(\cdot,s))}  \hskip - 0.5 cm Y^a|_{\cG(t,s)}\ ,\ \  \text{with}\  Y^a|_{\cG(t,s)}\=  du^a(Y|_{\cG(t,s)}) \ . 
\end{split}
\eeq
We also have that 
\beq \label{magic} 
 \b^{PC}(Y|_{\cG(t ,s)}) = \sum_{\d = 1}^{k} \sum_{\h = 0}^{\d-1} (-1)^{\h}  \frac{d^\h}{dt^\h}  \left(\frac{\p  L  }{\p q^i_{(\d)}} \right) Y^i_{(\d-(\h+1))}\bigg|_{\cG(t ,s)} \ ,
\eeq
where    $Y^i_{(\a)}$,   $0 \leq \a \leq k-1$,   are the  $\frac{\p}{\p q^i_{(\a)}}$-components of  $Y$. 
From  \eqref{348}, \eqref{magic},  the definition of  $\cP^{(\Sigma(\ve,s), \widecheck u^{(\ve,s)}, t)}$    and    Stokes' Theorem, we have   
\begin{multline*} 
 \int_0^T \left( \int_0^1  Y^a \frac{\p\cP^{(\Sigma(\ve,s), \widecheck u^{(\ve,s)}, t)} }{\p u^a} \Bigg|_{\widecheck u^{(\ve,s)}(t)} ds  \right) dt  = \\
= - \hskip - 1 cm \iint_{[\t_o - \ve - \gh \ve^2, \t_o + \gh \ve^2] \times [0,1]}  \frac{\p L}{\p u^a}\bigg|_{\cG(t,s)}  Y^a dt  ds = \hskip - 1 cm  \iint_{[\t_o - \ve - \gh \ve^2, \t_o + \gh \ve^2] \times [0,1]}  d \b^{PC}(X|_{\cG(t,s)}, Y|_{\cG(t,s)})  dt ds =  
\end{multline*}
\begin{multline} \label{ohoh}
=  \int_{\t_o - \ve - \gh \ve^2}^{\t_o + \gh \ve^2}  \bigg(  L|_{(j^{k}_t(\g_o), u^{(\t_o, \o_o, 0)}(t))} -  L|_{(j^{k}_t(\widecheck \g^{(\ve)}), u^{(\t_o, \o_o, \ve)}(t))} \bigg) dt +
\\
 +   \sum_{\d = 1}^{k} \sum_{\h = 0}^{\d-1} (-1)^{\h}  \int_0^1 \Bigg( \frac{d^\h}{dt^\h}  \left(\frac{\p  L  }{\p q^i_{(\d)}} \right) Y^i_{(\d-(\h+1))}\bigg|_{\cG(\t_o +  \gh \ve^2 ,s)}
- \\
-  \frac{d^\h}{dt^\h}  \left(\frac{\p  L  }{\p q^i_{(\d)}} \right) Y^i_{(\d-(\h+1))}\bigg|_{\cG(\t_o  - \ve -   \gh \ve^2 ,s)}\Bigg)ds = \\
=   \int_{\t_o - \ve - \gh \ve^2}^{\t_o + \gh \ve^2}   \left(\ \cP^{(\s_{(\ve)}, \widecheck u_{(\ve)}, t)}(u^{(\t_o, \o_o, \ve)}(t)) -  \cP^{(\s_o,  u_o, t)}(u_o(t)) \right)  dt +
 \\
 +   \sum_{\d = 1}^{k} \sum_{\h = 0}^{\d-1} (-1)^{\h}  \int_0^1 \Bigg( \frac{d^\h}{dt^\h}  \left(\frac{\p  L  }{\p q^i_{(\d)}} \right) Y^i_{(\d-(\h+1))}\bigg|_{\cG(\t_o +  \gh \ve^2 ,s)}
- \\
\hskip 5 cm -  \frac{d^\h}{dt^\h}  \left(\frac{\p  L  }{\p q^i_{(\d)}} \right) Y^i_{(\d-(\h+1))}\bigg|_{\cG(\t_o  - \ve -   \gh \ve^2 ,s)}\Bigg)ds .  \end{multline}
 The claim is therefore proven if we can show that  for any $s \in [0,1]$ the absolute value 
$$
\left|  \frac{d^\h}{dt^\h}  \left(\frac{\p  L  }{\p q^i_{(\d)}} \right) Y^i_{(\d-(\h+1))}\bigg|_{\cG(\t_o +  \gh \ve^2 ,s)}
-  \frac{d^\h}{dt^\h}  \left(\frac{\p  L  }{\p q^i_{(\d)}} \right) Y^i_{(\d-(\h+1))}\bigg|_{\cG(\t_o  - \ve -   \gh \ve^2 ,s)}\right| $$
  is  bounded above by   $\ve^2$  times   a constant  depending on  $\t_o$, $\cN$,   ${\vertiii{L}}_{k+2,\cN}$ and  ${\vertiii{\frac{\p L}{\p u}}}_{k+1,\cN}$. To check  this,  we first observe that 
  \begin{multline*}\left| \frac{d^\h}{dt^\h}  \left(\frac{\p  L  }{\p q^i_{(\d)}} \right) Y^i_{(\d-(\h+1))}\bigg|_{\cG(\t_o +  \gh \ve^2 ,s)}
-  \frac{d^\h}{dt^\h}  \left(\frac{\p  L  }{\p q^i_{(\d)}} \right) Y^i_{(\d-(\h+1))}\bigg|_{\cG(\t_o  - \ve -   \gh \ve^2 ,s)}\right|  \leq \\
\leq \left| \bigg(\frac{d^\h}{dt^\h}  \left(\frac{\p  L  }{\p q^i_{(\d)}} \right) \bigg|_{\cG(\t_o  +  \gh \ve^2 ,s)}
-  \frac{d^\h}{dt^\h}  \left(\frac{\p  L  }{\p q^i_{(\d)}} \right)  \bigg|_{\cG(\t_o  - \ve -   \gh \ve^2 ,s)}\bigg)Y^i_{(\d-(\h+1))}\bigg|_{\cG(\t_o  +   \gh \ve^2 ,s)}\right| + \\
+ \left| \frac{d^\h}{dt^\h}  \left(\frac{\p  L  }{\p q^i_{(\d)}} \right)  \bigg|_{\cG(\t_o  - \ve -   \gh \ve^2 ,s)}\bigg(Y^i_{(\d-(\h+1))}\bigg|_{\cG(\t_o +  \gh \ve^2 ,s)}
-   Y^i_{(\d-(\h+1))}\bigg|_{\cG(\t_o  - \ve -   \gh \ve^2 ,s)}\bigg)\right|
\end{multline*}
We now recall that for any $(t_o, s_o) \in [0, T] \times [0,1]$   
\beq \label{eheh} 
\begin{split}  
 &Y^i_{(\a)}\big|_{\cG(t_o ,s_o)} = \frac{\p }{\p s} \bigg|_{(t_o, s_o)}(q^{(\ve, s)i}_{(\a)}(t))\ ,\\
 &   \frac{\p}{\p t} \bigg|_{t_o} (Y^i_{(\a)} |_{\cG(\cdot ,s_o)})= \frac{\p^2 }{\p t \p s}\bigg|_{(t_o, s_o)} (q^{(\ve, s)i}_{(\a)}(t) ) = \frac{\p }{\p s}\bigg|_{(t_o,s_o)} (q^{(\ve, s)i}_{(\a+1)}(t)) = 
 Y^i_{(\a+1)}\big|_{\cG(t_o ,s_o)} \ ,
 \end{split}\eeq
where  $q^{(\ve, s)i}_{(\a)}(t_o)$ stands for the $q^i_{(\a)}$-component of the jet $j^{2k -1}_{t= t_o}(\widecheck \g^{(\ve, s)})$ of the $\wh \cK$-controlled curve $\widecheck \g^{(\ve, s)}$.
Combining \eqref{eheh},  the differentiability of $L$ with respect to  $u$, the explicit expressions of the Euler-Lagrange equations (which are obtained by taking at most  $k+1$ derivatives of $L$ with respect to  the jets coordinates)  and a straightforward  generalisation  of a classical fact on solutions to controlled differential equations (see e.g. \cite[Thm. 3.2.6]{BP} and the proof of Lemma \ref{normallemma}), 
 one can check that  for any $ 0 \leq \b \leq k$ and any $(t_o, s_o) \in [\t_o - \ve - \gh \ve^2, \t_o + \gh \ve^2] \times [0,1]$
\begin{multline*} 
 \left| Y^i_{(\b)}\big|_{\cG(t_o ,s_o)} \right| \leq \\
\leq(\ve + 2 \gh \ve^2) \bigg(\sup_{
 {\scriptscriptstyle [\t_o - \ve - \gh \ve^2, \t_o + \gh \ve^2]}
}  \Big|\widecheck u^{(\t_o, \o_o, \ve)}(t) -  u_o(t) \Big| \bigg)  e^{(\ve + 2 \gh \ve^2) K_{(\cN, L)}{\vertiii{L}}_{k+2, \cN}}  K_{(\cN, L)}{\vertiii{\frac{\p L}{\p u}}}_{k+1,\cN} \hskip - 0.5 cm  \leq\\
\leq (\ve + 2 \gh \ve^2) \diam(\wh K) e^{(\ve + 2 \gh \ve^2) K_{(\cN, L)}{\vertiii{L}}_{k+2, \cN}}K_{(\cN, L)}  {\vertiii{\frac{\p L}{\p u}}}_{k+1,\cN}   \ , 
\end{multline*}
where $\diam(\wh K)$  is the diameter of the relatively compact set $\wh K \subset \bR^M$.
Consequently, for any $ 0 \leq \a \leq k-1$ 
\begin{multline} \label{342}
 \left| Y^i_{(\a)}\big|_{\cG(t_o ,s_o)} -    Y^i_{(\a)}\bigg|_{\cG(\t_o  - \ve -   \gh \ve^2 ,s_o)}\right| = \bigg|\int_{t = \t_o  - \ve -   \gh \ve^2}^{t_o}
 \frac{\p}{\p t}\bigg|_t  Y^i_{(\a)}\big|_{\cG(\cdot ,s_o)}  dt\bigg|  = \\
=  \bigg|  \int_{t = \t_o  - \ve -   \gh \ve^2}^{t_o}
 Y^i_{(\a+1)}\big|_{\cG(t ,s_o)}dt  \bigg| 
\leq \\
\leq (\ve + 2 \gh \ve^2)^2\diam(\wh K) e^{(\ve + 2 \gh \ve^2) K_{(\cN, L)}{\vertiii{L}}_{k+2, \cN}} K_{(\cN, L)} {\vertiii{\frac{\p L}{\p u}}}_{k+1,\cN}   \ . 
\end{multline}
Hence,  for any $0 \leq \a \leq k-1$, $0 \leq \b \leq k$, 
\begin{multline} \label{344} \left| \frac{d^\h}{dt^\h}  \left(\frac{\p  L  }{\p q^i_{(\b)}} \right)  \bigg|_{\cG(\t_o  - \ve -   \gh \ve^2 ,s)} \right|\cdot
 \left|\bigg(Y^i_{(\a)}\bigg|_{\cG(\t_o +  \gh \ve^2 ,s)}
-   Y^i_{(\a)}\bigg|_{\cG(\t_o  - \ve -   \gh \ve^2 ,s)}\bigg)\right| \leq\\
\leq(\ve + 2 \gh \ve^2)^2\diam(\wh K) {\vertiii{L}}_{k+1, \cN} e^{(\ve + 2 \gh \ve^2) K_{(\cN, L)}{\vertiii{L}}_{k+2, \cN}} K_{(\cN, L)} {\vertiii{\frac{\p L}{\p u}}}_{k+1,\cN}  .
\end{multline}
A similar line of arguments  
yields to the  estimate
\begin{multline} \label{345}
\left| \bigg(\frac{d^\h}{dt^\h}  \left(\frac{\p  L  }{\p q^i_{(\d)}} \right) \bigg|_{\cG(\t_o  +  \gh \ve^2 ,s)}
-  \frac{d^\h}{dt^\h}  \left(\frac{\p  L  }{\p q^i_{(\d)}} \right)  \bigg|_{\cG(\t_o  - \ve -   \gh \ve^2 ,s)}\bigg)\right| \cdot \left|Y^i_{(\d-(\h+1))}\bigg|_{\cG(\t_o  +   \gh \ve^2 ,s)}\right| \leq\\
\leq 2 (\ve + 2 \gh \ve^2)^2\diam(\wh K) {\vertiii{L}}_{k+1, \cN} e^{(\ve + 2 \gh \ve^2) K_{(\cN, L)}{\vertiii{L}}_{k+2, \cN}} K_{(\cN, L)} {\vertiii{\frac{\p L}{\p u}}}_{k+1,\cN}  .
\end{multline}
From \eqref{ohoh}, \eqref{344}, \eqref{345}, the conclusion follows. \hfill\qed
\par
\medskip
\subsection{The proof of Lemma  \ref{corolletto-bis} } \label{proofsecondlemma}
It  suffices to prove that  the   constant  $\gn_{(\t_o, \cN, L, \frac{\p L}{\p u})}$ of  Sublemma \ref{lemma37}   is independent on ${\vertiii{L}}_{k+2, \cN} $ and  ${\vertiii{\frac{\p L}{\p u}}}_{k+1, \cN} $.  By \eqref{ohoh},  this  is proven if  we can show that, under the assumption  \eqref{3.54}, then 
\begin{multline} \label{355}
  \sum_{\d = 1}^{k} \sum_{\h = 0}^{\d-1} (-1)^{\h}  \int_0^1 \Bigg( \frac{d^\h}{dt^\h}  \left(\frac{\p  L  }{\p q^i_{(\d)}} \right) Y^i_{(\d-(\h+1))}\bigg|_{\cG(\t_o +  \gh \ve^2 ,s)}
- \\
\hskip 5 cm -  \frac{d^\h}{dt^\h}  \left(\frac{\p  L  }{\p q^i_{(\d)}} \right) Y^i_{(\d-(\h+1))}\bigg|_{\cG(\t_o  - \ve -   \gh \ve^2 ,s)}\Bigg)ds = 0\ .
\end{multline}
For this,   we  observe that,  by Stokes' Theorem, \eqref{348} and \eqref{magic}, 
\begin{multline}\label{356}
\sum_{\d = 1}^{k} \sum_{\h = 0}^{\d-1} (-1)^{\h}  \int_0^1 \Bigg( \frac{d^\h}{dt^\h}  \left(\frac{\p  L  }{\p q^i_{(\d)}} \right) Y^i_{(\d-(\h+1))}\bigg|_{\cG(\t_o - \ve -  \gh \ve^2 ,s)} - \\
-\sum_{\d = 1}^{k} \sum_{\h = 0}^{\d-1} (-1)^{\h}  \int_0^1 \Bigg( \frac{d^\h}{dt^\h}  \left(\frac{\p  L  }{\p q^i_{(\d)}} \right) Y^i_{(\d-(\h+1))}\bigg|_{\cG(0 ,s)} = \\
= \hskip - 1 cm \iint_{[0, \t_o - \ve - \gh \ve^2] \times [0,1]}   \hskip - 1 cmd \b^{PC}(X|_{\cG(t,s)}, Y|_{\cG(t,s)})  dt ds - \\
-   \int^{\t_o - \ve - \gh \ve^2}_0  \bigg(  L|_{(j^{k}_t(\g_o), u_o(t))} -  L|_{(j^{k}_t(\widecheck \g^{(\ve)}), u^{(\t_o, \o_o, \ve)}(t))} \bigg) dt  
= - \hskip - 1 cm \iint_{[0, - \t_o - \ve - \gh \ve^2] \times [0,1]}  \frac{\p L}{\p u^a}\bigg|_{\cG(t,s)}  Y^a dt  ds- \\
-   \int^{\t_o - \ve - \gh \ve^2}_0  \bigg(  L|_{(j^{k}_t(\g_o), u_o(t))} -  L|_{(j^{k}_t(\widecheck \g^{(\ve)}), u^{(\t_o, \o_o, \ve)}(t))} \bigg) dt \ ,
 \end{multline}
 where $\b^{PC}$ is the $1$-form  \eqref{39}.  
We now recall that in the region $[0,  \t_o - \ve - \gh \ve^2] \times [0, 1]$ the components $Y^a$ of the vector field $Y$ are identically   $0$. Moreover, if $L$ has the form \eqref{specialform},  the controlled Euler-Lagrange equations imply that the value of $L$ is $0$ at all jets of each $\wh \cK$-controlled curve $\g$, so   that 
$$\int^{\t_o - \ve - \gh \ve^2}_0  \bigg(  L|_{(j^{k}_t(\g_o), u_o(t))} -  L|_{(j^{k}_t(\widecheck \g^{(\ve)}), u^{(\t_o, \o_o, \ve)}(t))} \bigg) dt = 0\ .$$
 Hence  \eqref{356} and \eqref{3.54} imply
\begin{multline} \sum_{\d = 1}^{k} \sum_{\h = 0}^{\d-1} (-1)^{\h}  \int_0^1 \Bigg( \frac{d^\h}{dt^\h}  \left(\frac{\p  L  }{\p q^i_{(\d)}} \right) Y^i_{(\d-(\h+1))}\bigg|_{\cG(\t_o - \ve -  \gh \ve^2 ,s)} = \\
= \sum_{\d = 1}^{k} \sum_{\h = 0}^{\d-1} (-1)^{\h}  \int_0^1 \Bigg( \frac{d^\h}{dt^\h}  \left(\frac{\p  L  }{\p q^i_{(\d)}} \right) Y^i_{(\d-(\h+1))}\bigg|_{\cG(0 ,s)} = 0\ .
\end{multline}
A  similar argument yields  
\begin{multline} \sum_{\d = 1}^{k} \sum_{\h = 0}^{\d-1} (-1)^{\h}  \int_0^1 \Bigg( \frac{d^\h}{dt^\h}  \left(\frac{\p  L  }{\p q^i_{(\d)}} \right) Y^i_{(\d-(\h+1))}\bigg|_{\cG(\t_o +  \gh \ve^2 ,s)} = \\
= \sum_{\d = 1}^{k} \sum_{\h = 0}^{\d-1} (-1)^{\h}  \int_0^1 \Bigg( \frac{d^\h}{dt^\h}  \left(\frac{\p  L  }{\p q^i_{(\d)}} \right) Y^i_{(\d-(\h+1))}\bigg|_{\cG(T ,s)} = 0\ .
\end{multline}
 From this \eqref{355} follows. 
\par
\smallskip
\section{Suggested   investigations}
\label{added}
As we mentioned in the Introduction, here we want  to point out  some problems of Control Theory  where  our  two-steps approach (= a preliminary analysis  based on classical  results  of Differential Geometry, followed by   arguments  devoted to  reduce the regularity assumptions)  has good chances to produce new  results or  to enlighten some particular aspects of the dynamics of controlled systems.  The  discussion is intentionally very sketchy,  because its  purpose  is merely  to provide   suggestions and motivations for future studies. \par
\smallskip
\subsection{Maximum Principles  in  Continuum Dynamics}
Consider the following toy problem.  Let $\cE$ be an (unbounded) elastic continuum,  whose elements are described by  just one space-variable, denoted by $s \in \bR$,  that evolves in the time $t$. 
The deformations of  such continuum  are represented by functions $x(t,s)$ of the time and space variables and are assumed  to satisfy a hyperbolic equations of the form 
\beq  \label{hyper} \frac{\p^2 x}{\p t^2} - \frac{\p^2 x}{\p s^2} = f( t, s, u(t,s))\ , \eeq
where $f: \bR^3 \to \bR$ is a fixed smooth function and  $u(t,s)$  is a control map  with values in a compact set $K \subset \bR$.  The composed map  $f(t, s, u(t,s))$ might be physically interpreted  as  a (density of a) {\it dead load}  attached at the  elements of the continuum   and   varying   in time. Note also that, when $f(s, t, u)$ is  linear in $u$, the equation  \eqref{hyper}  is in the class of  controlled hyperbolic equations,  which is intensively  studied in the theory of control problems governed by partial differential equations (see \cite{Li, To}).\par 
\smallskip
 Following our usual two-step approach, let us at first restrict the discussion  of this toy problem to  deformations $x(s,t)$ and control maps $u(s,t)$ of class $\cC^\infty$ and satisfying all needed assumptions (as,  for instance, rapidly decreasing properties for $s \to 0$) that  may guarantee  that all subsequent  arguments  are  meaningful. 
\par
Consider  the following  problem: {\it given an initial  condition  for $x(t,s)$ at $t = 0$
\beq \label{initcond} x( 0, s) = \f(s)\ , \qquad\ \frac{\p x}{\p t}\bigg|_{(0, s)} = \psi(s) ,\eeq
 look for  a  load  $u_o(t,s)$ such that the corresponding solution to \eqref{hyper}   satisfying \eqref{initcond}
 minimises the integral at $t = T$ (= the terminal cost)} 
 \beq \label{cost-multiple} {\mathsf C}(x_o(T, s)) = \int_\bR \ell(x(T, s)) d s\eeq
 where $\ell(x)$ is a prescribed smooth real  function.
This can be classified as  a control problem, whose  optimal controls are the loads $u_o(t,s)$ satisfying the above minimising requirement.  Inspired by the discussions in \cite{CS, CGS}  about the  control problems involving  just one independent variable, it is natural  to start studying this new type of control problem by considering the {\it  controlled Lagrangian density} on the $2$-jets of maps 
$(t, s) \mapsto (x(t,s), p(t,s))$, 
 defined by 
\beq \cL^{(u(t,s))} (t, s, x, x_t, x_{tt}, x_s, x_{ss}, p)\= p(x_{tt} - x_{ss} - f\left(t,s, u(t, s)\right)) +\left(\ell(x) + t \frac{\p \ell(x)}{\p x} x_t\right) .\eeq
One can  check that, for any fixed choice of the control function $u(t,s)$, 
\begin{itemize}[leftmargin = 15pt]
\item[(1)] The Euler-Lagrange equations determined by $ \cL^{(u(t,s))}$ give a system of two partial differential equations,  the first equal to  \eqref{hyper}, the second equal to  the  hyperbolic equation on $p(t,s)$
\beq   \frac{\p^2 p}{\p t^2} -  \frac{\p^2 p}{\p s^2} =0 \ ; \eeq
\item[(2)] If the pair $(x(t,s), p(t,s))$ is a solution to  the Euler-Lagrange equations in (1), then 
\begin{multline*}\iint_{S = \{0 \leq t \leq T, s \in \bR\}} \hskip - 20 pt \cL^{(u(t,s))} \left(t, s, x(t,s),  \frac{\p x(t,s)}{\p t} , \frac{\p^2 x(t,s)}{\p t^2} ,  \frac{\p x(t,s)}{\p s} ,  \frac{\p^2 x(t,s)}{\p s^2} , p(t,s)\right)dt ds  = \\ =  {\mathsf C}(x(T, s)) \ .\end{multline*}
\end{itemize}
All this shows  that the new setting   is extremely   close to what is considered in \cite{CS, CGS} for control problems  with  differential constraints involving just one independent variable. We are confident  that the same line of arguments considered  there (and in particular the ``road map'' presented  in  \cite[Sect. 2.2]{CGS}) can be followed
for this toy problem and many other  control problems with constraints given by {\it partial differential equations}. This would lead to   analogs of the PMP  (compare, for instance,    \cite{Li, BM}) under strong regularity assumptions, results which   can be considered as  the {\it first step of differential-geometric type} of the  approach we are promoting.  The  direct proof of Theorem \ref{main1}  given in this paper can be then considered as  guiding line for extending the results  of the ``first step''  to reach results under low regularity assumptions. 
\par
\medskip
\subsection{Dynamics of controlled systems with  higher order constraints}
Consider a dynamical system which is subjected to a second order differential constraint in normal form and independent on time, that is  of the form 
\beq\label{diffe-1-Noether} \frac{d ^2 x^j}{d t^2} = f^j\left(x(t) ,u(t)\right)  \ ,\qquad 1 \leq j \leq n\ .\eeq
with control curve  $u(t) = (u^a(t))$ taking values in a relatively compact set $K \subset \bR^m$. Following the first step of our  two-step approach, let us at  first  assume that all data satisfy strong regularity assumptions
(i.e. assume that $f$ is smooth, $u(t)$ varies in the class of smooth curves, $K$ has smooth boundary,  etc.), so that the most common differential geometric tools might be used. Let us also 
denote by $\cL(x, x_{(1)}, x_{(2)}, p, u)$  the second order controlled Lagrangian \eqref{specialform} associated with this control problem: 
$$\cL(x, x_{(1)}, x_{(2)}, p, u) \= p_j\left(x^j_{(2)} - f^j\left(x ,u(t)\right) \right)\ .$$
We remark that,  for any fixed choice of a control curve $u_o(t)$, the Euler-Lagrange equations of $\cL(x, x_{(1)}, x_{(2)}, p, u_o(t)) $  for  $x$ and $p$ coincide with 
the Euler-Lagrange equations of  the equivalent Lagrangian (their difference is a null Lagrangian)
\beq \label{eqLag} \wt \cL(x, x_{(1)}, x_{(2)}, p,p_{(1)} , u_o(t)) \= p_{(1)j} x^j_{(1)} - p_j f^j\left(x ,u_o(t)\right) \ .\eeq
 If we consider the coordinates $(\wt x^i, \wt p_j)$ related with the $(x^i, p_j)$ by 
\beq  x^i = \frac{1}{\sqrt{2}}\wt x^i + \frac{1}{\sqrt{2}} \wt p_i\ ,\qquad p_j =\frac{1}{\sqrt{2}} \wt x_j - \frac{1}{\sqrt{2}}\wt p_j\ ,\eeq
the new Lagrangian \eqref{eqLag} takes  a very  familiar form, namely
\begin{multline}  \label{eqLag*} \wt \cL(x, x_{(1)}, x_{(2)}, p,p_{(1)} , u_o(t)) \= \frac{1}{2}  \sum_{i = 1}^n\left( \left(\wt x_{(1)i}\right)^2 - \frac{1}{2}  \left(\wt p_{(1)i} \right)^2\right) + V(\wt x, \wt p, u_o(t))  \\
\text{where} \ V(\wt x, \wt p, u) \= - \left(\frac{1}{\sqrt{2}} \wt x_j - \frac{1}{\sqrt{2}}\wt p_j\right) f^j\left( \frac{1}{\sqrt{2}}\wt x + \frac{1}{\sqrt{2}} \wt p\ ,  \,u\right) \ .\end{multline}
This is a Lagrangian that describes the dynamics on a Lorentzian $2$-manifold  of a system subjected to  force with  time-dependent potential  $V(\wt x, \wt p, u_o(t))$. 
It is therefore possible  to use  a variety of  well known mathematical physics tools  to study the dynamics of such controlled systems. For instance, for any given choice of a smooth $u_o(t)$,   studying symmetries of  $V$ and using  Noether Theorem (\cite{Ol}),   all conservation laws  that are satisfied (or, more interesting,  violated)  can be explicitly determined. In particular, 
in the time intervals on which $u_o(t)$ is constant (recall that, in several classical   settings, the optimal control $u_o(t)$ is constant a.e.) an appropriate non-positively defined energy is conserved by the corresponding controlled evolution. \par
\smallskip
Furthermore, if we 
denote by    $(Q^I) \= \left(\begin{array}{c} x^i \\
p_j \end{array}  \right)$, we may observe that the Hessian $\frac{\p \wt \cL}{\p Q^I_{(1)} \p Q^J_{(1)}}$ is non-degenerate, a property that allows a formulation of the differential constraints into  a Hamiltonian formulation in  the phase space spanned  by the  coordinates  $Q = (Q^I) = (x^i, p_j)$ and their duals $ P = (P_K)$. In a sense, this  would be a ``true Hamiltonian presentation'' of the  constraints of the controlled system, very much different from the traditional Pontryagin's Hamiltonian type presentation.   We think that it would be  quite important to get a clear view 
of  the relations between these two distinct Hamiltonian type presentations  of the differential constraints and of their dependences on the needle variations. As usual,  answers  to any question  in this topic can at first be obtained via differential geometric tools under strong regularity assumptions.  Secondly one can  extend the   results to  the  lowest possible regularity assumptions following the ideas  of this paper. Similar investigations might -- and, in our opinion,  should -- be made for controlled  systems with differential constraints of order  higher than two and/or by means of the alternative   presentations of Lagrangian type,  which are discussed  in \cite{CGS}.

\bigskip
\bigskip
\font\smallsmc = cmcsc8
\font\smalltt = cmtt8
\font\smallit = cmti8
\hbox{\parindent=0pt\parskip=0pt
\vbox{\baselineskip 9.5 pt \hsize=3.1truein
\obeylines
{\smallsmc 
Franco Cardin
Dipartimento di Matematica 
``Tullio  Levi-Civita''
Universit\`a degli Studi di Padova
Via Trieste, 63
I-35121 Padua
ITALY
\ 
\ 
\ 
}\medskip
{\smallit E-mail}\/: {\smalltt cardin@math.unipd.it 
\ 
}
}
\hskip 0.0truecm
\vbox{\baselineskip 9.5 pt \hsize=3.7truein
\obeylines
{\smallsmc
Cristina Giannotti \& Andrea Spiro
Scuola di Scienze e Tecnologie
Universit\`a di Camerino
Via Madonna delle Carceri
I-62032 Camerino (Macerata)
ITALY
\ 
}\medskip
{\smallit E-mail}\/: {\smalltt cristina.giannotti@unicam.it
\smallit E-mail}\/: {\smalltt andrea.spiro@unicam.it}
}
}
\end{document}